\title{The Nakayama automorphism of the \\ almost Calabi-Yau algebras \\ associated to $SU(3)$ modular invariants}
\author{
        David E. Evans and Mathew Pugh \\ \\
        School of Mathematics, \\
        Cardiff University, \\
        Senghennydd Road, \\
        Cardiff, CF24 4AG, \\
        Wales, U.K.
}
\date{\today}

\documentclass[12pt]{article}
\usepackage{amssymb}
\usepackage{graphicx}
\usepackage[usenames]{color}

\newtheorem{Def}{Definition}[section]
\newtheorem{Prop}[Def]{Proposition}

\newtheorem{Cor}[Def]{Corollary}
\newtheorem{Thm}[Def]{Theorem}

\begin{document}
\maketitle

\begin{abstract}
We determine the Nakayama automorphism of the almost Calabi-Yau algebra $A$ associated to the braided subfactors or nimrep graphs associated to each $SU(3)$ modular invariant. We use this to determine a resolution of $A$ as an $A$-$A$ bimodule, which will yield a projective resolution of $A$.
\end{abstract}

\section{Introduction} \label{sect:intro}

The $SU(2)$ and $SU(3)$ modular invariant partition functions were classified by Cappelli, Itzykson and Zuber \cite{cappelli/itzykson/zuber:1987ii} and Gannon \cite{gannon:1994} respectively.
The $SU(2)$ theory is closely related to the preprojective algebras of Coxeter-Dynkin quivers. The object of this paper is to study the analogous finite dimensional superpotential algebras associated to the $SU(3)$ invariants.

The classical McKay correspondence relates finite subgroups $\Gamma$ of $SU(2)$ with the algebraic geometry of the quotient Kleinian singularities $\mathbb{C}^2/\Gamma$ \cite{reid:2002} but also with the classification of $SU(2)$ modular invariants \cite{cappelli/itzykson/zuber:1987ii, zuber:2002}, classification of subfactors of index less than 4 \cite{ocneanu:1988, izumi:1991, kawahigashi:1995, bion-nadal:1991, izumi:1994}, and quantum subgroups of $SU(2)$ \cite{ocneanu:2000ii, ocneanu:2002, xu:1998, bockenhauer/evans:1999i, bockenhauer/evans:1999ii, bockenhauer/evans/kawahigashi:1999, bockenhauer/evans/kawahigashi:2000}. The study of quotient singularities and their resolution has been assisted with the study of the structure of certain noncommutative algebras.
Minimal resolutions of Kleinian singularities can be described via the moduli space of representations of the preprojective algebra associated to the action of $\Gamma$. This leads to general programme to understand singularities via a noncommutative algebra $A$, often called a noncommutative resolution, whose centre corresponds to the coordinate ring of the singularity \cite{vandenBergh:2004}. The algebra should be finitely generated over its centre, and the desired favourable resolution is the moduli space of representations of $A$, whose category of finitely generated modules is derived equivalent to the category of coherent sheaves of the resolution.
In the case of a quotient singularity $\mathbb{C}^3/\Gamma$ for a finite subgroup $\Gamma$ of $SU(3)$, the corresponding noncommutative algebra $A$ is a Calabi-Yau algebra of dimension 3.

Calabi-Yau algebras arise naturally in the study of Calabi-Yau manifolds, providing a noncommutative version of conventional Calabi-Yau geometry.
An algebra $A$ is Calabi-Yau of dimension $n$ if the bounded derived category of the abelian category of finite dimensional $A$-modules is a Calabi-Yau category of dimension $n$. In this case the global dimension of $A$ is $n$ \cite{bocklandt:2008}.
The derived category of coherent sheaves over an $n$-dimensional Calabi-Yau manifold is a Calabi-Yau category of dimension $n$ and they appear naturally in the study of boundary conditions of the $B$-model in superstring theory over the manifold. For more on Calabi-Yau algebras, see e.g. \cite{bocklandt:2008, ginzburg:2006}.

In \cite[Remark 4.5.7]{ginzburg:2006} Ginzburg introduced, in his terminology, $q$-deformed Calabi-Yau algebras. In the case where $q$ is not a root of unity, these algebras are Calabi-Yau algebras of dimension 3.
In this paper we will study these algebras in the case where $q$ is a root of unity, which are the $SU(3)$ generalizations of preprojective algebras \cite{gelfand/ponomarev:1979} for the Coxeter-Dynkin diagrams $ADE$. We will call these algebras \emph{almost Calabi-Yau algebras}.

In Section \ref{sect:Preliminaries} we bring together the strands needed from subfactor theory, modular tensor categories and their modules, planar algebras and categorical approaches.
We begin in Section \ref{sect:verlinde-sector} by describing our notation for the representation theory of $SU(n)$ at level $k$. In Section \ref{sect:Hecke} we recall the generalized Temperley-Lieb algebras for $SU(n)$, which are representations of the Hecke algebra. Then in Section \ref{sect:VerlindeAlgebra} we review the description of the Verlinde algebra and fusion rules for quantum $SU(n)$ in terms of endomorphisms of a hyperfinite type $\mathrm{III}_1$ factor $N$. This system of endomorphisms has the structure of a modular tensor category. For a braided inclusion $N \subset M$ we obtain a module category which produces a nimrep, or non-negative matrix integer representation of the fusion rules, which is associated to the $SU(n)$ modular invariants. In Section \ref{sect:subfactors} we give a subfactor description of the Verlinde algebra, where the space of intertwiners is identified with the algebra of paths on the nimrep graph.
In Section \ref{sect:SU(2)categorical} we move to the categorical picture, and give a diagrammatic and categorical description of the Verlinde algebra and fusion rules for quantum $SU(2)$. This is based on the diagrammatic representation of the Temperley-Lieb algebra of Kauffman \cite{kauffman:1987}. A categorical approach to Temperley-Lieb algebras was given in \cite{turaev:1994, yamagami:2003, cooper:2007}.
In Section \ref{sect:nimreps} we describe $SU(2)$ module categories in terms of preprojective algebras, following the ideas of Cooper \cite{cooper:2007}. We relate this in Section \ref{sect:BGPA-GHJ} to braided subfactors using the Goodman-de la Harpe-Jones construction \cite{goodman/de_la_harpe/jones:1989} and its manifestation in the bipartite graph planar algebra construction of Jones \cite{jones:2000}. We give the analogous description for $SU(3)$, where the diagrammatic and categorical description of Section \ref{sect:SU(3)categorical} is based on the $SU(3)$-Temperley-Lieb algebra and its diagrammatic representation constructed in \cite{evans/pugh:2009iii} using the $A_2$-webs of Kuperberg \cite{kuperberg:1996}.
In Section \ref{sect:nimrepsSU(3)} we describe the module categories in terms of (almost) Calabi-Yau algebras, following the ideas of Cooper \cite{cooper:2007}.
The corresponding nimreps, the $SU(3)$ $\mathcal{ADE}$ graphs, are the graphs associated to the $SU(3)$ modular invariants.
In Section \ref{sect:A2GPA-GHJ} we relate the construction of these almost Calabi-Yau algebras to braided subfactors using the $SU(3)$ Goodman-de la Harpe-Jones construction \cite{evans/pugh:2009ii} and its manifestation in the $SU(3)$-graph planar algebra construction \cite{evans/pugh:2009iv}.

Then in Section \ref{sect:Hilbert-SU(3)} we compute the Hilbert series of dimensions associated to these almost Calabi-Yau algebras.
The McKay graph of $SU(3)$ is built out of closed paths of length 3, which corresponds to the fact that the fundamental representation $\rho$ of $SU(3)$ satisfies $\rho \otimes \rho \otimes \rho \ni \mathbf{1}$.
One can build an Ocneanu cell system $W$ on the McKay graph of a subgroup of $SU(3)$ or an $\mathcal{ADE}$ graph $\mathcal{G}$, which attaches a complex number to each closed path of length three on the edges of $\mathcal{G}$ \cite{ocneanu:2000ii}. These yield relations on paths of equal length, and one obtains a superpotential algebra $A = A(\mathcal{G},W)$ by taking the quotient by the ideal generated by these relations.
For the $\mathcal{ADE}$ graphs $\mathcal{G}$, we take potentials built on the cell system $W$ computed in \cite{evans/pugh:2009i}, and study the Hilbert series of dimensions of the corresponding quotient algebras $A(\mathcal{G},W)$, which are almost Calabi-Yau algebras.
The preprojective algebras for $SU(2)$ braided subfactors, and the almost Calabi-Yau algebras for $SU(3)$ braided subfactors, are given by the image under a functor $F$ of the quotient $S$ of the tensor algebra generated by the fundamental sector by a tensor ideal which makes $S$ symmetric, that is,
\begin{equation}
F(S) = A(\mathcal{G},W),
\end{equation}
where the functor $F$ is essentially the module category arising from the braided subfactor.

If $H_n$ is the matrix of dimensions of paths of length $n$ in some graph $\mathcal{G}$ in the quotient algebra $A = A(\mathcal{G},W)$, with the indices of the matrix labeled by the vertices, then the matrix Hilbert series $H_A$ of the algebra $A$ is defined as $H_A(t) = \sum H_n t^n$. Let $\Delta_{\mathcal{G}}$ be the adjacency matrix of $\mathcal{G}$.
Then if $\mathcal{G}$ is the McKay graph of a subgroup of $SU(3)$ then $A = A(\mathcal{G},W)$ is a Calabi-Yau algebra of dimension 3 \cite[Theorem 4.4.6]{ginzburg:2006}, and by \cite[Theorem 4.6]{bocklandt:2008} its Hilbert series is given by:
\begin{equation} \label{eqn:H(t)-CYd}
H_A(t) = \frac{1}{1 - \Delta_{\mathcal{G}} t + \Delta_{\mathcal{G}}^T t^2 - t^3}.
\end{equation}
We prove in Theorem \ref{thm:SU(3)Hilbert} that if $\mathcal{G}$ is a finite $SU(3)$ $\mathcal{ADE}$ graph which carries a cell system $W$, and thus yields a braided subfactor \cite{evans/pugh:2009ii}, then
\begin{equation} \label{eqn:H(t)-qCYd}
H_A (t) = \frac{1 - P t^h}{1 - \Delta_{\mathcal{G}} t + \Delta_{\mathcal{G}}^T t^2 - t^3},
\end{equation}
where $P$ is a permutation matrix corresponding to a $\mathbb{Z}_3$ symmetry of the graph, and $h = k+3$ is the Coxeter number of $\mathcal{G}$, where $k$ is the level of $SU(3)$. The permutation matrix $P$ corresponds precisely to the Nakayama permutation for $A$.
This result was mentioned without proof in \cite{evans/pugh:2009v}.
In \cite[Proposition 3.14]{brenner/butler/king:2002} the Hilbert series was given for a $(p,q)$-Koszul algebra (or almost Koszul algebra), of which (\ref{eqn:H(t)-qCYd}) is a particular case, where $A = A(\mathcal{G},W)$ is a $(h-3,3)$-Koszul algebra, see Section \ref{sect:nimrepsSU(3)}.

The dual $A^{\ast} = \mathrm{Hom}(A, \mathbb{C})$ is an $A$-$A$ bimodule, not usually identified with ${}_A A_A$ or ${}_1 A_1$ with standard right and left actions but with ${}_1 A_{\beta}$ with standard left action and the standard right action twisted by an automorphism $\beta$, the Nakayama automorphism \cite{yamagata:1996}. The Nakayama automorphism measures how far away $A$ is from being symmetric. In the case of a preprojective algebra $\Pi$ of a Dynkin quiver, this Nakayama automorphism is identified (up to a sign) with an involution on the underlying Dynkin diagram, which is trivial in all cases, except for the Dynkin diagrams $A_n$, $D_{2n+1}$, $E_6$ where it is the unique non-trivial involution \cite{erdmann/snashall:1998i, erdmann/snashall:1998ii}. In the case of $SU(3)$, the Nakayama automorphism is identified in Theorem \ref{Thm:A^=1_A_beta} with the $\mathbb{Z}_3$ symmetry of the underlying $\mathcal{ADE}$ graph given by the permutation $P$ in (\ref{eqn:H(t)-qCYd}).
This answers a question we posed in \cite[p.411]{evans/pugh:2009v}.
The Hilbert series (\ref{eqn:H(t)-qCYd}) for $A$ is a key ingredient in our determination of the Nakayama automorphism for the $SU(3)$ $\mathcal{ADE}$ graphs.
We also use the Ocneanu cells $W(\triangle)$ computed in \cite{evans/pugh:2009i}, and exploit their $\mathbb{Z}_3$-invariance and in most cases the non-vanishing of (certain linear combinations of) cells which appear in the determination of a basis for $A$. It does not appear that the non-vanishing of these linear combinations can be deduced merely from the existence of cells in \cite{ocneanu:2000ii}.

In Section \ref{sect:nakayama} we obtain the first part of a resolution of $A$ as an $A$-$A$ bimodule
$$A \otimes_R A[3] \rightarrow A \otimes_R V^{\ast} \otimes_R A[1] \rightarrow A \otimes_R V \otimes_R A \rightarrow A \otimes_R A \rightarrow A \rightarrow 0,$$
where $R$, $V$ are the $A$-$A$ bimodules generated by the vertices, edges of $\mathcal{G}$ respectively, the $A$-$A$ bimodule $V^{\ast}$ is the dual space of $V$, and $B[m]$ denotes the graded space $B$ but with grading shifted by $m$.
The algebra $A$ is a Calabi-Yau algebra if and only if the kernel of the leftmost map is zero, as in \cite{ginzburg:2006, bocklandt:2008, broomhead:2008}.
In our case, this kernel $\Omega^4(A)$ is non-zero and is determined in Theorem \ref{thm:Omega^4(A)}.
We show that $\Omega^4(A)$ is isomorphic as an $A$-$A$ bimodule to ${}_1 A_{\beta^{-1}}$, and thus obtain a finite resolution of $A$ as an $A$-$A$ bimodule
$$0 \rightarrow {}_1 A_{\beta^{-1}}[h] \rightarrow A \otimes_R A[3] \rightarrow A \otimes_R V^{\ast} \otimes_R A[1] \rightarrow A \otimes_R V \otimes_R A \rightarrow A \otimes_R A \rightarrow A \rightarrow 0.$$
This \emph{almost Calabi-Yau} condition should be compared with the Calabi-Yau condition expressed above.

This resolution of $A$ will yield a projective resolution of $A$ as an $A$-$A$ bimodule.
The objective in deriving this resolution is to provide a basis for the computation of the Hochschild (co)homology and cyclic homology of the algebras $A(\mathcal{G},W)$ for the $SU(3)$ $\mathcal{ADE}$ graphs. Beginning with a pair $(\mathcal{G},W)$ given by a cell system $W$ on an $SU(3)$ $\mathcal{ADE}$ graph $\mathcal{G}$, we construct a subfactor $N \subset M$ which yields a nimrep which recovers the graph $\mathcal{G}$ as described in Section \ref{sect:VerlindeAlgebra}. Then we can construct the algebra $A(\mathcal{G},W)$ whose Hochschild (co)homology and cyclic homology only depends on the original pair $(\mathcal{G},W)$, or equivalently, on the subfactor $N \subset M$. Thus the Hochschild (co)homology and cyclic homology of $A$ should be regarded as invariants for the subfactor $N \subset M$.

\section{Preliminaries} \label{sect:Preliminaries}

In this section we bring together the strands needed from subfactor theory, modular tensor categories and their modules, planar algebras and categorical approaches as outlined in the Introduction.

\subsection{Representations of $SU(n)$ and $SU(n)_k$} \label{sect:verlinde-sector}

Here we describe our notation for the representation theory of $SU(n)$ at level $k \leq \infty$.
Every irreducible representation $\lambda_m$ of $SU(n)$ is classified by a signature, or highest weight, $m = (m_1,m_2,\ldots,m_{n-1})$, where $m_i$ are integers such that $m_1 \geq m_2 \geq \cdots \geq m_{n-1} \geq 0$, for $i=1,\ldots,n-1$.
A signature $m$ can be represented by a Young diagram with at most $n-1$ rows, and $m_i$ boxes in the $i^{\mathrm{th}}$ row, $i = 1,\ldots,n-1$.
The irreducible positive energy representations of the loop group of $SU(n)$ at level $k$, or $SU(n)_k$, are described by the irreducible representations of $SU(n)$ whose signature has at most $k$ columns.

For $SU(2)$ the signatures are just the integers $k \geq 0$. The fundamental representation is $\rho = \lambda_1$, and the irreducible representations of $SU(2)$ satisfy the fusion rules $\lambda_m \otimes \rho = \lambda_{m-1} \oplus \lambda_{m+1}$ for $m \geq 1$, and $\lambda_0 \otimes \rho = \rho$.
The fusion graph for $SU(2)$ is the infinite Dynkin diagram $A_{\infty}$.
It is well known that the $k^{\mathrm{th}}$ symmetric product of $\mathbb{C}^2$ gives the irreducible level $k$ representation.
The irreducible representations of $SU(2)_k$ satisfy the same fusion rules as those for $SU(2)$, only now $\lambda_m$ is also understood to be zero if $m > k$.

For $SU(3)$ the signatures are pairs $(m_1,m_2)$ with $m_1 \geq m_2 \geq 0$. We will replace the pair $(m_1,m_2)$ by the Dynkin labels $(p,l) = (m_2-m_1,m_1)$, where now $p,l \geq 0$. The conjugate representation of $\lambda_{(p,l)}$ is $\overline{\lambda}_{(p,l)} = \lambda_{(l,p)}$. The fundamental representation $\rho = \lambda_{(1,0)}$ generates every irreducible representation of $SU(3)$ with its conjugate $\overline{\rho}$, and the irreducible representations of $SU(3)$ satisfy the fusion rules
\begin{equation} \label{eqn:SU(3)fusion_rules}
\lambda_{(p,l)} \otimes \rho = \lambda_{(p,l-1)} \oplus \lambda_{(p-1,l+1)} \oplus \lambda_{(p+1,l)}, \;\; \lambda_{(p,l)} \otimes \overline{\rho} = \lambda_{(p-1,l)} \oplus \lambda_{(p+1,l-1)} \oplus \lambda_{(p,l+1)},
\end{equation}
where $\lambda_{(r,s)}$ is understood to be zero if $r < 0$ or $s < 0$.
The fusion graph is the infinite graph $\mathcal{A}^{(\infty)}$ (see \cite[Figure 4]{evans/pugh:2009i}).
The irreducible representations of $SU(3)_k$ satisfy the same fusion rules as those for $SU(3)$, only now $\lambda_{(p,l)}$ is also understood to be zero if $p+l>k$, and the fusion graph is the truncated graph $\mathcal{A}^{(k+3)}$.

\subsection{The generalized Temperley-Lieb algebras} \label{sect:Hecke}

Let $M_n = \mathrm{End}(\mathbb{C}^n)$. By Weyl duality, the fixed point algebra of $\otimes^m M_n$ under the product adjoint action of $SU(n)$ is generated by a representation $\sigma \rightarrow g_{\sigma}$ on $\otimes^m \mathbb{C}^n$ of the group ring of the symmetric, or permutation, group $S_m$. This algebra is generated by unitary operators $g_j$, $j=1,\ldots,m-1$, which represent transpositions $(j,j+1)$, satisfying the relations
\begin{eqnarray}
g_j^2 & = & 1, \label{eqn:Hecke-1} \\
g_i g_j & = & g_j g_i, \quad |i-j|>1, \label{eqn:Hecke-2} \\
g_i g_{i+1} g_i & = & g_{i+1} g_i g_{i+1}, \label{eqn:braiding_relation}
\end{eqnarray}
and the vanishing of the antisymmetrizer
\begin{equation}\label{SU(N)condition}
\sum_{\sigma \in S_m} \mathrm{sgn}(\sigma) g_{\sigma} = 0.
\end{equation}
Writing $g_j = 1 - U_j$, these unitary generators and relations are equivalent to the self-adjoint generators $\mathbf{1}, U_j$, $j=1, \ldots, m-1$, and relations
\begin{center}
\begin{minipage}[b]{15cm}
 \begin{minipage}[b]{2cm}
  \begin{eqnarray*}
  \textrm{H1:}\\
  \textrm{H2:}\\
  \textrm{H3:}
  \end{eqnarray*}
 \end{minipage}
 \hspace{2cm}
 \begin{minipage}[b]{7cm}
  \begin{eqnarray*}
  U_i^2 & = & \delta U_i,\\
  U_i U_j & = & U_j U_i, \quad |i-j|>1,\\
  U_i U_{i+1} U_i - U_i & = & U_{i+1} U_i U_{i+1} - U_{i+1},
  \end{eqnarray*}
 \end{minipage}
\end{minipage}
\end{center}
where $\delta = 2$, and the analogue of (\ref{SU(N)condition}).

There is a $q$-version of this algebra, which is a representation of a Hecke algebra. This is the centralizer of a representation of the quantum group $SU(n)_q$ (or the universal enveloping algebra), with a deformation of (\ref{eqn:Hecke-1}) to
\begin{equation} \label{eqn:Hecke-1'}
(q^{-1} - g_j)(q + g_j) = 0.
\end{equation}
The invertible generators $g_j$, $j=1,\ldots,m-1$, satisfy the relations (\ref{eqn:Hecke-2}), (\ref{eqn:braiding_relation}), (\ref{eqn:Hecke-1'}) and the vanishing of the $q$-antisymmetrizer \cite{di_francesco/zuber:1990}
\begin{equation}\label{SU(N)q condition}
\sum_{\sigma \in S_m} (-q)^{|I_{\sigma}|} g_{\sigma} = 0,
\end{equation}
where $g_{\sigma} = \prod_{i \in I_{\sigma}} g_i$ if $\sigma = \prod_{i \in I_{\sigma}} (i,i+1)$.
Then writing $g_j = q^{-1} - U_j$, we are interested in the generalized Temperley-Lieb algebra generated by self-adjoint operators $\mathbf{1}, U_j$, $j=1, \ldots, m-1$, satisfying H1-H3 and the analogue of (\ref{SU(N)q condition}), where now $\delta = q+q^{-1}$.
In the cases $n=2,3$, which we are interested in, (\ref{SU(N)q condition}) reduces for $SU(2)$ to the Temperley-Lieb condition
\begin{equation}
U_i U_{i \pm 1} U_i - U_i = 0, \label{eqn:SU(2)q_condition}
\end{equation}
and for $SU(3)$ it is
\begin{equation} \label{eqn:SU(3)q_condition}
\left( U_i - U_{i+2} U_{i+1} U_i + U_{i+1} \right) \left( U_{i+1} U_{i+2} U_{i+1} - U_{i+1} \right) = 0.
\end{equation}

There are minor errors in a parallel discussion in Section 2 of the published version of \cite{evans/pugh:2009iii} which have been corrected in the arXiv version.

\subsection{Braided subfactors and modular invariants} \label{sect:VerlindeAlgebra}

Let $A$ and $B$ be type $\mathrm{III}$ von Neumann factors. A unital $\ast$-homomorphism $\rho:A\rightarrow B$ is called a $B$-$A$ morphism.
The positive number $d_{\rho}=[B:\rho(A)]^{1/2}$ is called the statistical dimension of $\rho$; here $[B:\rho(A)]$ is the Jones-Kosaki index \cite{jones:1983, kosaki:1986} of the subfactor $\rho(A)\subset B$.
Some $B$-$A$ morphism $\rho'$ is called equivalent to $\rho$ if $\rho'=\mathrm{Ad}(u)\circ\rho$ for some unitary $u\in B$.
The equivalence class $[\rho]$ of $\rho$ is called the $B$-$A$ sector of $\rho$. If $\rho$ and $\sigma$ are $B$-$A$ morphisms with finite statistical dimensions, then the vector space of intertwiners
$$\mathrm{Hom}(\rho,\sigma)=\{ t\in B: t\rho(a)=\sigma(a)t \,, \,\, a\in A \}$$
is finite-dimensional, and we denote its dimension by $\langle \rho, \sigma \rangle$.
A $B$-$A$ morphism is called irreducible if $\langle \rho,\rho \rangle=1$, i.e. if $\mathrm{Hom}(\rho,\rho) = \mathbb{C} \mathbf{1}_B$.
Then, if $\langle \rho, \tau \rangle \neq 0$ for some (possibly reducible) $B$-$A$ morphism $\tau$, then $[\rho]$ is called an irreducible subsector of $[\tau]$ with multiplicity $\langle \rho, \tau \rangle$.
An irreducible $A$-$B$ morphism $\overline{\rho}$ is a conjugate morphism of the irreducible $\rho$ if and only if $[\overline{\rho}\rho]$ contains the trivial sector $[\mathrm{id}_A]$ as a subsector, and then $\langle \rho\overline{\rho}, \mathrm{id}_B \rangle = 1 = \langle \overline{\rho}\rho, \mathrm{id}_A \rangle$ automatically \cite{izumi:1991}.

The Verlinde algebra is realised in the subfactor models by systems of endomorphisms ${}_N \mathcal{X}_N$ of the hyperfinite type $\mathrm{III}_1$ factor $N$.
That is, ${}_N \mathcal{X}_N$ denotes a finite system of finite index irreducible endomorphisms of a factor $N$ in the sense that different elements of ${}_N \mathcal{X}_N$ are not unitary equivalent, for any $\lambda \in {}_N \mathcal{X}_N$ there is a representative $\overline{\lambda} \in {}_N \mathcal{X}_N$ of the conjugate sector $[\overline{\lambda}]$, and ${}_N \mathcal{X}_N$ is closed under composition and subsequent irreducible decomposition.
In the case of WZW models associated to $SU(n)$ at level $k$, the Verlinde algebra is a non-degenerately braided system of endomorphisms ${}_N \mathcal{X}_N$, labelled by the positive energy representations of the loop group of $SU(n)_k$ on a type $\mathrm{III}_1$ factor $N$, with fusion rules $\lambda \mu = \bigoplus_{\nu} N_{\lambda \nu}^{\mu} \nu$ which exactly match those of the positive energy representations \cite{wassermann:1998}. The fusion matrices $N_{\lambda} = [N_{\rho \lambda}^{\sigma}]_{\rho,\sigma}$ are a family of commuting normal matrices which give a representation themselves of the fusion rules of the positive energy representations of the loop group of $SU(n)_k$, $N_{\lambda} N_{\mu} = \sum_{\nu} N_{\lambda \nu}^{\mu} N_{\nu}$.
This family $\{ N_{\lambda} \}$ of fusion matrices can be simultaneously diagonalised:
\begin{equation} \label{eqn:verlinde_formula}
N_{\lambda} = \sum_{\sigma} \frac{S_{\sigma, \lambda}}{S_{\sigma,1}} S_{\sigma} S_{\sigma}^{\ast},
\end{equation}
where $1$ is the trivial representation, and the eigenvalues $S_{\sigma, \lambda}/S_{\sigma,1}$ and eigenvectors $S_{\sigma} = [S_{\sigma, \mu}]_{\mu}$ are described by the statistics $S$ matrix.
Moreover, there is equality between the statistics $S$- and $T$- matrices and the Kac-Peterson modular $S$- and $T$- matrices which perform the conformal character transformations \cite{kac:1990}, thanks to \cite{frohlich/gabbiani:1990, fredenhagen/rehren/schroer:1992, wassermann:1998}.

The key structure in the conformal field theory is the modular invariant partition function $Z$. In the subfactor setting this is realised by
a braided subfactor $N \subset M$ where trivial (or permutation) invariants in the ambient factor $M$ when restricted to $N$ yield $Z$. This would mean that the dual canonical endomorphism is in $\Sigma({}_N \mathcal{X}_N)$, i.e. decomposes as a finite linear combination of endomorphisms in ${}_N \mathcal{X}_N$.
Indeed if this is the case for the inclusion $N \subset M$, then the process of $\alpha$-induction allows us to analyse the modular invariant,
providing two extensions of $\lambda$ on $N$ to endomorphisms $\alpha^{\pm}_{\lambda}$ of $M$, such that the matrix $Z_{\lambda,\mu} = \langle \alpha_{\lambda}^+, \alpha_{\mu}^- \rangle$ is a modular invariant \cite{bockenhauer/evans/kawahigashi:1999, bockenhauer/evans:2000, evans:2003}.

Let ${}_N \mathcal{X}_M$, ${}_M \mathcal{X}_M$ denote a system of endomorphisms consisting of a choice of representative endomorphism of each irreducible subsector of sectors of the form $[\lambda \overline{\iota}]$, $[\iota \lambda \overline{\iota}]$ respectively, for each $\lambda \in {}_N \mathcal{X}_N$, where $\iota: N \hookrightarrow M$ is the inclusion map which we may consider as an $M$-$N$ morphism, and $\overline{\iota}$ is a representative of its conjugate $N$-$M$ sector.
The action of the system  ${}_N \mathcal{X}_N$ on the $N$-$M$ sectors ${}_N \mathcal{X}_M$ produces a nimrep (non-negative matrix integer representation of the fusion rules) $G_{\lambda} G_{\mu} = \sum_{\nu} N_{\lambda \nu}^{\mu} G_{\nu}$,
whose spectrum reproduces exactly the diagonal part of the modular invariant, i.e.
\begin{equation} \label{eqn:verlinde_formulaG}
G_{\lambda} = \sum_{\sigma} \frac{S_{\sigma,\lambda}}{S_{\sigma,1}} \psi_{\sigma} \psi_{\sigma}^{\ast},
\end{equation}
with the spectrum of $G_{\lambda} = \{ S_{\mu, \lambda}/S_{\mu,1}$ with multiplicity $Z_{\mu,\mu} \}$ \cite[Theorem 4.16]{bockenhauer/evans/kawahigashi:2000}. The labels $\mu$ of the non-zero diagonal elements are called the exponents of $Z$, counting multiplicity. A modular invariant for which there exists a nimrep whose spectrum is described by the diagonal part of the invariant is said to be nimble.

The systems ${}_N \mathcal{X}_N$, ${}_N \mathcal{X}_M$, ${}_M \mathcal{X}_M$ are (the irreducible objects of) tensor categories of endomorphisms with the Hom-spaces as their morphisms. Thus ${}_N \mathcal{X}_N$ gives a braided modular tensor category, and ${}_N \mathcal{X}_M$ a module category.
The structure of the module category ${}_N \mathcal{X}_M$ is the same as a tensor functor $F$ from ${}_N \mathcal{X}_N$ to the category $\mathrm{Fun}({}_N \mathcal{X}_M,{}_N \mathcal{X}_M)$ of additive functors from ${}_N \mathcal{X}_M$ to itself, see \cite{ostrik:2003}. That is, $F$ is essentially the module category ${}_N \mathcal{X}_M$.

The classification of $SU(2)$ modular invariants is due to Cappelli, Itzykson and Zuber \cite{cappelli/itzykson/zuber:1987ii}.
They label the modular invariant with an $ADE$ graph $\mathcal{G}$ such that the diagonal part $Z_{\mu,\mu}$ of the invariant is exactly the multiplicity of the eigenvalue $S_{\mu,\rho}/S_{\mu,1}$ of $\mathcal{G}$, where 1, $\rho$ denote the trivial, fundamental representations respectively.
Since these $ADE$ graphs can be matched to the affine Dynkin diagrams -- the McKay graphs of the representation theory of the  finite subgroups of $SU(2)$ -- di Francesco and Zuber \cite{di_francesco/zuber:1990} were guided to find candidates for classifying graphs for $SU(3)$ modular invariants by first considering the McKay graphs of the finite subgroups of $SU(3)$ to produce a candidate list of $\mathcal{ADE}$ graphs whose spectra described the diagonal part of the modular invariant. They proposed candidates for most of the modular invariants, except for the conjugate invariants $\mathcal{A}^{\ast}$ as they restricted themselves to only look for graphs which are three-colourable.
In the subfactor theory, this is understood in the following way. Suppose $N \subset M$ is a braided subfactor which realises the modular invariant $Z_{\mathcal{G}}$. Evaluating the nimrep $G$ at the fundamental representation $\rho$, we obtain for the inclusion $N \subset M$ a matrix $G_{\rho}$, which is the adjacency matrix for the $ADE$ graph $\mathcal{G}$ which labels the modular invariant.
Every $SU(2)$ modular invariant is realised, and all nimreps are realised by subfactors \cite{ocneanu:2000ii, ocneanu:2002, xu:1998, bockenhauer/evans:1999i, bockenhauer/evans:1999ii, bockenhauer/evans/kawahigashi:1999, bockenhauer/evans/kawahigashi:2000}, apart from the tadpole nimreps of the orbifolds of the even $A$'s (see e.g. \cite{bockenhauer/evans:2001} for an explanation of the failure of the tadpole nimreps).
Behrend, Pearce, Petkova and Zuber \cite{behrend/pearce/petkova/zuber:2000} (see also \cite{zuber:2002}) systematically proposed nimreps as a framework for boundary conformal field theory. The $N$-$M$ system ${}_N \mathcal{X}_M$ corresponds to boundary fields in their language, and the $M$-$M$ system ${}_M \mathcal{X}_M$ to defect lines.
B\"{o}ckenhauer and Evans \cite{bockenhauer:1999} understood that nimrep graphs for the $SU(3)$ conjugate invariants were not three-colourable. This was also realised simultaneously by Behrend, Pearce, Petkova and Zuber \cite{behrend/pearce/petkova/zuber:2000} and Ocneanu \cite{ocneanu:2002}. The figures for the complete list of the $\mathcal{ADE}$ graphs are given in \cite{behrend/pearce/petkova/zuber:2000, ocneanu:2002, evans/pugh:2009i}.
The classification of $SU(3)$ modular invariants was shown to be complete by Gannon \cite{gannon:1994}, and the complete list is given in \cite{evans/pugh:2009ii}. Ocneanu claimed \cite{ocneanu:2000ii, ocneanu:2002} that all $SU(3)$ modular invariants were realised by subfactors and this was shown in \cite{xu:1998, bockenhauer/evans:1999i, bockenhauer/evans:1999ii, bockenhauer/evans/kawahigashi:1999, bockenhauer/evans:2001, bockenhauer/evans:2002, evans/pugh:2009i, evans/pugh:2009ii}.
However, the classification of nimreps is incomplete if one relaxes the condition that the nimrep be compatible with a modular invariant \cite{gannon:2001, graves:2010}.
Ostrik \cite{kirillov/ostrik:2002, ostrik:2003} took up a categorical description of subfactor $\alpha$-induction, see \cite[Remark 14]{ostrik:2003}, and this was taken further by Fjelstad, Fr\"{o}hlich, Fuchs, Schweigert and Runkel as a categorical framework for conformal field theory. See \cite{fuchs/runkel/schweigert:2010} for a review.

\subsection{Subfactors} \label{sect:subfactors}

Suppose we have a system of endomorphisms ${}_N \mathcal{X}_N$ of a type $\mathrm{III}_1$ factor $N$ for $SU(n)_k$, $k \leq \infty$, where $\rho$ denotes the endomorphism in ${}_N \mathcal{X}_N$ corresponding to the fundamental generator.
We can form the tunnel
\begin{equation} \label{eqn:sector_tunnel}
\cdots \subset \rho\overline{\rho}\rho(N) \subset \rho\overline{\rho}(N) \subset \rho(N) \subset N.
\end{equation}
By decomposing the sectors of $1, \overline{\rho}, \overline{\rho}\rho, \overline{\rho}\rho\overline{\rho}, \ldots$ into irreducible sectors we can obtain the Bratteli diagram of the higher relative commutants of $\rho(N) \subset N$. If $[\lambda^{(0)}]$ is an irreducible at an even level of the Bratteli diagram and $[\lambda^{(0)}] [\rho]$ decomposes into irreducibles as $[\lambda^{(0)}] [\rho] = \bigoplus_{i=1}^s [\lambda^{(i)}]$, for irreducible sectors $[\lambda^{(i)}]$, $i=1,\ldots,s$, then there is an edge from the vertex $[\lambda^{(0)}]$ in the Bratteli diagram to the vertices $[\lambda^{(i)}]$, whilst if $[\lambda^{(0)}]$ is an irreducible at an odd level of the Bratteli diagram, we consider instead the decomposition of $[\lambda^{(0)}] [\overline{\rho}]$ into irreducibles. The Bratteli diagram obtained in this way is identical to that obtained for the Jones-Wenzl type $\mathrm{II}_1$ $SU(n)$ subfactors \cite{wenzl:1988}.
The principal graph is the bipartite graph constructed by deleting at each level the vertices belonging to the old sectors (that is, any vertex at a given level which appeared at a previous level of the Bratteli diagram) and the edges emanating from them \cite[Definition 4.6.5]{goodman/de_la_harpe/jones:1989}.
The decomposition of the sectors of $1, \rho, \rho\overline{\rho}, \rho\overline{\rho}\rho, \ldots$ into irreducibles yields the dual principal graph in a similar way.

The decomposition of the sectors of the form $(\overline{\rho}\rho)^m$ and $(\overline{\rho}\rho)^m\overline{\rho}$ will not usually produce all the irreducible sectors in ${}_N \mathcal{X}_N$, so to obtain all the irreducible sectors we also consider the decomposition of more general products $\rho^m \overline{\rho}^l$ into irreducibles, $m,l \geq 0$.
In this way we recover the graph $\mathcal{A}$ when $n=2,3$, with vertices labelling the irreducible sectors, and the edges representing multiplication by the fundamental generator $\rho$. This corresponds to idempotent completion in the categorical language of Section \ref{sect:SU(2)categorical}. The principal graph, respectively dual principal graph is however only the 0-1, 0-$(n-1)$ part of the full graph $\mathcal{A}$, where the edges are now undirected, and where either would be the entire graph only in the case when $n = 2$ \cite{evans/kawahigashi:1994}.

There is an identification between intertwiners and explicit paths on the intertwining graph $\mathcal{A}$ \cite{izumi:1998}, \cite[Section 3.5]{evans/kawahigashi:1998}. For a sector $[\lambda_i]$ at an even level in the Bratteli diagram, the intertwiners $\mathrm{Hom}(\rho\lambda_i, \lambda_j)$ are identified with the edges from $[\lambda_i]$ to $[\lambda_j]$ on $\mathcal{A}$, whilst for a sector $[\lambda_i']$ at an odd level in the Bratteli diagram, the intertwiners $\mathrm{Hom}(\overline{\rho}\lambda_i', \lambda_j)$ are identified with the edges from $[\lambda_i']$ to $[\lambda_j]$. Let $T(a_j)$ denote an intertwiner labeled by an edge $a_j$ of $\mathcal{A}$, and for a path $x = a_1 a_2 \cdots a_s$ on $\mathcal{A}$, define $T(x) := T(a_1) T(a_2) \cdots T(a_s)$. Then the $T(x)$ are an orthogonal basis of the intertwiners between some endomorphisms. The spaces of intertwiners are the Hilbert spaces on which the system ${}_N \mathcal{X}_N$ acts.
In this way the spaces $\mathrm{Hom}(\rho^{m_1} \overline{\rho}^{l_1}, \rho^{m_2} \overline{\rho}^{l_2})$ of morphisms are identified with the span of all pairs $(x_1,x_2)$ of paths $x_1$, $x_2$ on $\mathcal{A}$ and its opposite graph $\mathcal{A}^{\mathrm{op}}$ where all the edges of $\mathcal{A}$ are reversed, where $x_j$ has $m_j$, $l_j$ edges on $\mathcal{A}$, $\mathcal{A}^{\mathrm{op}}$ respectively, $j=1,2$.
In particular, the algebras $\mathrm{Hom}(\rho^m \overline{\rho}^l, \rho^m \overline{\rho}^l)$ are identified with the path algebra
(in the usual operator algebraic sense \cite{evans/kawahigashi:1998})
on $\mathcal{A}$, $\mathcal{A}^{\mathrm{op}}$.
Jones projections $e_j$ \cite{jones:1983} for the tunnel (\ref{eqn:sector_tunnel}) are identified with those in the path algebra on $\mathcal{A}$, $\mathcal{A}^{\mathrm{op}}$, where edges are alternately on $\mathcal{A}$ and $\mathcal{A}^{\mathrm{op}}$ \cite{evans/kawahigashi:1994}.
The Jones projections in the path algebra are given by the product $c^{\ast}c$, where the annihilation and creation operators $c$, $c^{\ast}$ are defined in Section \ref{sect:nimreps}.

We now focus on the cases $n=2,3$, where for $k < \infty$, $q$ is a $(k+n)^{\mathrm{th}}$ root of unity.
For $n=2$, let $\rho = \overline{\rho}$ denote the endomorphism corresponding to the fundamental generator of $SU(2)_k$.
The tunnel (\ref{eqn:sector_tunnel}) defines Jones projections $e_j$ which generate the Temperley-Lieb algebra.
The intertwiner space is generated by the Jones projections $e_j$, so that $\mathrm{Hom} (\rho^m, \rho^m) \cong \mathcal{TL}_m := \mathrm{alg}(1, e_1, e_2, \ldots e_{m-1})$.
Jones-Wenzl projections $f_m = 1 - e_1 \vee \cdots \vee e_{m-1}$ are given by \cite{wenzl:1987}:
\begin{equation} \label{eqn:JW-projection}
f_{m+1} = f_m - \frac{[2]_q[m]_q}{[m+1]_q} f_m e_m f_m,
\end{equation}
where the quantum integer $[m]_q$ is defined by $[m]_q = (q^m - q^{-m})/(q - q^{-1})$. Here, $a \vee b$ denotes the projection such that $\mathrm{Ran}(a \vee b) = \mathrm{Ran}(a) + \mathrm{Ran}(b)$ where $a$, $b$ are projections on a Hilbert space.
For $m < k-1 \leq \infty$, the Jones-Wenzl projection $f_m$ is the minimal central projection corresponding to the new sector that appears at level $m$ in the Bratteli diagram.
For the fixed point algebra $(\otimes^m M_2)^{SU(2)}$ which is equal to the Temperley-Lieb algebra with $q=1$, the Jones-Wenzl projection $f_m$ is the projection on the $(m+1)$-dimensional representation indexed by $m$ in the intertwiner space $\mathrm{Hom}(\rho^m,\rho^m)$, where $\rho$ is the fundamental representation of $SU(2)$ on $\mathbb{C}^2$.

For $n=3$, we take the fundamental generator $\rho$ of $SU(3)_k$.
The Jones projections $e_j$, $j=1, \ldots, 2m-1$, for the tunnel (\ref{eqn:sector_tunnel}) are identified with projections in the algebras $\mathrm{Hom}(\rho^m \overline{\rho}^m, \rho^m \overline{\rho}^m)$, whilst the algebras $\mathrm{Hom} (\rho^m, \rho^m)$ are generated by the $A_2$-Temperley-Lieb operators $U_j$, $j=1, \ldots, m-1$.
For the fixed point algebra $(\otimes^m M_3)^{SU(3)}$, there is a generalized Jones-Wenzl projection $f_{(m,0)}$ which is the projection on the representation indexed by $\lambda_{(m,0)}$ in the intertwiner space $\mathrm{Hom}(\rho^m,\rho^m)$, where $\rho$ is the fundamental representation of $SU(3)$ on $\mathbb{C}^3$.
More generally, if we have an $A_2$-Temperley-Lieb algebra generated by self-adjoint operators $U_i$ with parameter $\delta = q + q^{-1}$, generalized Jones-Wenzl projections $f_{(m,0)}$ (called projectors in \cite{suciu:1997}, also called clasps \cite{kuperberg:1996}, magic elements \cite{ohtsuki/yamada:1997}) are defined
by \cite{suciu:1997}:
\begin{eqnarray*}
f_{(0,0)} & = & 1, \qquad \qquad f_{(1,0)} \;\; = \;\; \mathbf{1}, \\
f_{(m+1,0)} & = & f_{(m,0)} - \frac{[m]_q}{[m+1]_q} f_{(m,0)} U_m f_{(m,0)}.
\end{eqnarray*}

Note that the nimrep graph given by the module category ${}_N \mathcal{X}_M$ is not the principal graph for the braided inclusion $N \subset M$. The principal graph of a subfactor of index $<4$ can only be $A$, $D_{\mathrm{even}}$, $E_6$ or $E_8$, whilst the nimrep graph of a subfactor is any Coxeter-Dynkin diagram including $D_{\mathrm{odd}}$ and $E_7$. Indeed, the principal graphs of the braided inclusions $N \subset M$ usually have index which exceeds 4, and are those of the Goodman-de la Harpe-Jones construction and their generalizations discussed in \cite{goodman/de_la_harpe/jones:1989, evans/kawahigashi:1998, evans/pugh:2009ii}.
The even, odd vertices of the $ADE$ graphs are the $A$-$A$, $B$-$B$ systems respectively for a subfactor $A \subset B$. For the braided $SU(2)$ subfactors, all $A_n$ vertices (both even and odd) are represented as $N$-$N$ sectors, and all the vertices of the classifying graph $\mathcal{G}$ appear as $N$-$M$ sectors.

When we consider modular invariants, their module categories or nimreps, other graphs $\mathcal{G}$ will appear. The above intertwining discussion already leads us to the path Hilbert space $\mathbb{C}\mathcal{A}$, which is the vector space of paths on $\mathcal{A}$, identified with the intertwiners $T(x)$ where $x$ is a path on $\mathcal{A}$.
The module category ${}_N \mathcal{X}_M$ from a braided inclusion $N \subset M$ yields a nimrep $G$ and we obtain the path Hilbert space $\mathbb{C}\mathcal{G}$, the vector space of paths on $\mathcal{G} = G_{\rho}$, identified with the intertwiners $T(x)$ where $x$ is now a path on $\mathcal{G}$.
Denote by $(\mathbb{C}\mathcal{G})_j$ the space of paths of length $j$ on $\mathcal{G}$.
The path Hilbert space $\mathbb{C}\mathcal{G}$ is a graded algebra where multiplication $(\mathbb{C}\mathcal{G})_i \times (\mathbb{C}\mathcal{G})_j \rightarrow (\mathbb{C}\mathcal{G})_{i+j}$ of two paths $x \in (\mathbb{C}\mathcal{G})_i$, $y \in (\mathbb{C}\mathcal{G})_j$ is given by concatenation of paths $xy$, and is defined to be zero if $r(x) \neq s(y)$, where $s(x)$, $r(x)$ denote the source, range vertices of the path $x$ respectively.
The endomorphisms $\mathrm{End} \left( (\mathbb{C}\mathcal{G})_j \right)$ on $(\mathbb{C}\mathcal{G})_j$ are the $|(\mathbb{C}\mathcal{G})_j| \times |(\mathbb{C}\mathcal{G})_j|$-matrices, with rows and columns labeled by the paths of length $j$ on $\mathcal{G}$.
The $\mathrm{End} \left( (\mathbb{C}\mathcal{G})_j \right)$ have an algebra structure given by matrix multiplication.
Thus there are two different notions of path algebra of $\mathcal{G}$. In the theory of operator algebras, the path algebra of $\mathcal{G}$ is usually $\bigoplus_{j \geq 0} \mathrm{End} \left( (\mathbb{C}\mathcal{G})_j \right)$ \cite{evans/kawahigashi:1998}. In this paper however, we will work with the graded algebra $\mathbb{C}\mathcal{G}$, as in e.g. \cite{bocklandt:2008}.

\subsection{$SU(2)$ Categorical Approach} \label{sect:SU(2)categorical}

In this section we will describe the Verlinde algebra and fusion rules for $SU(2)$ in the diagrammatic and categorical language of the Temperley-Lieb algebra \cite{kauffman:1987, turaev:1994, yamagami:2003, cooper:2007}.

\begin{figure}[tb]
\begin{center}
  \includegraphics[width=45mm]{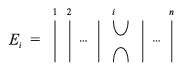}\\
 \caption{diagram $E_i \in \mathcal{T}_n$}\label{Fig:fig_Ei}
\end{center}
\end{figure}

Let $q$ be real or a root of unity, so that $\delta = [2]_q$ is real. Denote by $\mathcal{T}_{m,n}$ the set of all planar diagrams consisting of a rectangle with $m$, $n$ vertices along the top, bottom edge respectively, and with $(m+n)/2$ curves, called strings, inside the rectangle so that each vertex is the endpoint of exactly one string, and the strings do not cross each other.
Let $\mathcal{V}_{m,n}$ denote the free vector space over $\mathbb{C}$ with basis $\mathcal{T}_{m,n}$. Composition $RS$ of diagrams $R \in \mathcal{T}_{m,n}$, $S \in \mathcal{T}_{n,p}$ is given by gluing $S$ vertically below $R$ such that the vertices at the bottom of $R$ and the top of $S$ coincide, removing these vertices, and isotoping the glued strings if necessary to make them smooth. Any closed loops which may appear are removed, contributing a factor of $\delta$. The resulting diagram is in $\mathcal{T}_{m,p}$. This composition is clearly associative, and composition in $\mathcal{V} = \bigcup_{m,n \geq 0} \mathcal{V}_{m,n}$ is defined as its linear extension. The adjoint $R^{\ast} \in \mathcal{T}_{n,m}$ of a diagram in $R \in \mathcal{T}_{m,n}$ is given by reflecting $R$ about a horizontal line halfway between the top and bottom vertices of the diagram. This action is extended conjugate linearly to $\mathcal{V}$.
Let $E_i$ denote the diagram in $\mathcal{T}_n := \mathcal{T}_{n,n}$ illustrated in Figure \ref{Fig:fig_Ei}. For $\delta \geq 2$ there is an isomorphism $\mathcal{V}_n \cong \mathcal{TL}_n$ given by $\delta^{-1} E_i \rightarrow e_i$.

We will now define the Temperley-Lieb category $TL$ as a matrix category $TL = \mathrm{Mat}(C)$. We begin by defining $C$ to be the tensor category whose objects are projections in $\mathcal{V}_n := \mathcal{V}_{n,n}$, and whose morphisms $\mathrm{Hom}(p_1,p_2)$ between projections $p_i \in \mathcal{V}_{n_i}$, $i=1,2$, are given by the space $p_2 \mathcal{V}_{n_2,n_1} p_1$. We will use fraktur script to denote morphisms.
The tensor product is defined on the objects and morphisms by horizontal juxtaposition. The trivial object $\mathrm{id}_0$ is the empty diagram which is a projection in $\mathcal{V}_0$. (The category $C$ is the idempotent completion, or Karoubi envelope, of the category whose objects are non-negative integers, and whose morphisms are given by $\mathcal{V}_{m,n}$.)

In order to be able to take direct sums, we define the matrix category $TL = \mathrm{Mat}(C)$ to be the category with objects given by formal direct sums of objects in $C$, and morphisms $\mathrm{Hom}(p_1 \oplus \cdots \oplus p_{n_1}, q_1 \oplus \cdots \oplus q_{n_2})$ given by $n_2 \times n_1$ matrices, where the $i,j$-th entry is in $\mathrm{Hom}(p_j,q_i)$. The tensor product on $TL$ is given on objects by $(p_1 \oplus \cdots \oplus p_{n_1}) \otimes (q_1 \oplus \cdots \oplus q_{n_2}) = (p_1 \otimes q_1) \oplus (p_1 \otimes q_2) \oplus \cdots \oplus (p_{n_1} \otimes q_{n_2})$, and on morphisms by the usual tensor product on matrices with the tensor product for $C$ on matrix entries.
A projection $p \in TL$ is called simple if $\langle p,p \rangle = 1$.

We write $TL_n := \mathcal{V}_n$, and $\rho$ for the identity object in $TL_1$ consisting of a single vertical string. Then the identity diagram in $TL_n$, given by $n$ vertical strings, is expressed by $\rho^n := \otimes^n \rho$. We have $\mathrm{dim}(TL_0) = \mathrm{dim}(TL_1) = 1$ and $TL_0$, $TL_1$ have simple projections $f_0$ (the empty diagram), $f_1 = \rho$ respectively. Moving to $TL_2$, the identity diagram $\rho^2$ is a projection but is not simple, since $\langle \rho^2,\rho^2 \rangle = 2$. One of these morphisms is the identity diagram, the other is $\mathfrak{E}_1 = \mathrm{id}_{E_1}$. Since $E_1^2 = \delta E_1$, $e_1 = \delta^{-1} E_1$ is a projection. In fact, $e_1$ is isomorphic to $f_0$, as can be seen from the following isomorphisms $\psi: f_0 \rightarrow e_1$, $\psi^{\ast}: e_1 \rightarrow f_0$, where $\psi = \sqrt{[2]_q}^{-1}$ \includegraphics[width=5mm]{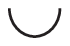}, and $\psi^{\ast}$ is defined by reflecting $\psi$ about its horizontal axis.
Then we have $\psi^{\ast}\psi = \mathfrak{f}_0 = \mathrm{id}_{f_0}$ and $\psi\psi^{\ast} = \mathfrak{e}_1 = \mathrm{id}_{e_1}$, where the morphism $\mathfrak{e}_1 = \delta^{-1} \mathfrak{E}_1$. Since $\langle f_0,\rho^2 \rangle = 1$, where the morphism is given by $\;$ \includegraphics[width=5mm]{fig_nakayama-cup} $\;$, and $\langle f_1,\rho^2 \rangle = 0$ (by parity), we have the decomposition $\rho^2 = f_0 \oplus f_2$, where $f_2$ is a simple projection in $TL_2$. In the same way, we obtain at each level $n$ that $\rho^n$ is a linear combination of $f_0, \ldots, f_{n-1}$ plus a new projection $f_n$, which turns out to be simple.
The morphisms $\mathfrak{f}_p = \mathrm{id}_{f_p}$ are Jones-Wenzl projections, and satisfy a similar recursion relation to (\ref{eqn:JW-projection}), with $e_m$ replaced by $\mathfrak{e}_m$.
These satisfy the properties \cite{wenzl:1988}:
\begin{center}
\includegraphics[width=80mm]{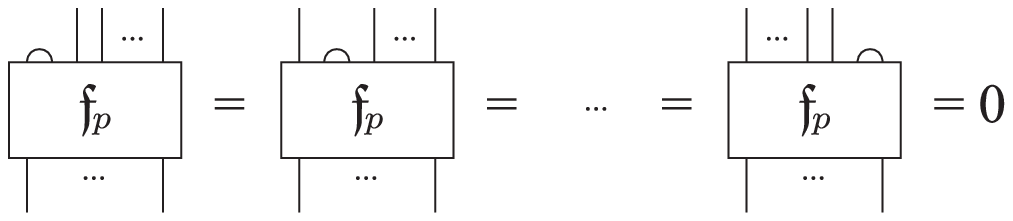} \\
\vspace{5mm}
\includegraphics[width=50mm]{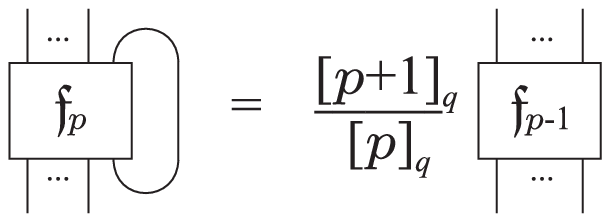}
$$\mathrm{tr}(\mathfrak{f}_p) = [p+1]_q,$$
\end{center}
where $\mathrm{tr}(\mathfrak{f}_p)$ is given by connecting the $i^{\mathrm{th}}$ string from the left along the top is connected to the $i^{\mathrm{th}}$ string from the left along the bottom for each $i = 1,\ldots,p$.
For $p \geq p'$ the Jones-Wenzl projections also satisfy the property $\mathfrak{f}_p (\mathrm{id}_{\rho^i} \otimes \mathfrak{f}_{p'} \otimes \mathrm{id}_{\rho^{p-p'-i}}) = \mathfrak{f}_p = (\mathrm{id}_{\rho^i} \otimes \mathfrak{f}_{p'} \otimes \mathrm{id}_{\rho^{p-p'-i}}) \mathfrak{f}_p$, for any $0 \leq i \leq p-p'$.
The morphisms $\mathfrak{f}_i$ and objects $f_i$ are denoted by $e_i$, $X_i$ respectively in \cite{cooper:2007}.
In \cite{morrison/peters/snyder:2008} there is some abuse of notation with both the objects and the morphisms given by the Jones-Wenzl projections denoted by $f^{(i)}$.

We have the relation
\begin{equation} \label{eqn:fusion_rule-fk}
f_p \otimes \rho \cong f_{p-1} \oplus f_{p+1}.
\end{equation}
This is seen from the isomorphisms $\psi: f_p \otimes \rho \rightarrow f_{p-1} \oplus f_{p+1}$, and $\psi^{-1} = \psi^{\ast}$, where
$$\psi = \left( \begin{array}{c} \includegraphics[width=22mm]{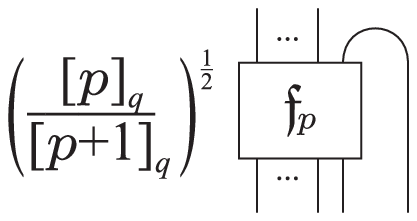} \\ \mathfrak{f}_{p+1} \end{array} \right),$$
and $\psi^{\ast}$ is defined as the transpose of $\psi$ where we replace each entry $a$ in $\psi$ by the reflection $a^{\ast}$ of $a$ about its horizontal axis.
Then it easy to check using the above properties and (\ref{eqn:JW-projection}) that $\psi\psi^{\ast} = \mathfrak{f}_{p-1} \oplus \mathfrak{f}_{p+1} = \mathrm{id}_{f_{p-1}} \oplus \mathrm{id}_{f_{p+1}}$ and $\psi^{\ast}\psi = \mathfrak{f}_p \otimes \mathrm{id}_{\rho} = \mathrm{id}_{f_p \otimes \rho}$.
Suppose that we have obtained simple projections $f_p$ for $0 \leq p \leq n$, and we obtain a new projection $f_{n+1}$ at level $n+1$ as above. From the relation (\ref{eqn:fusion_rule-fk}) we obtain the decomposition
\begin{equation} \label{eqn:decomp-Xk}
\rho^p = \bigoplus_{j=0}^{\lfloor p/2 \rfloor} \left[ \begin{array}{c} p \\ j \end{array} \right] f_{p-j},
\end{equation}
where $\left[ \begin{array}{c} p \\ j \end{array} \right] = C^p_j - C^p_{j-1}$ for binomial coefficients $C^p_j$, $0 \leq p \leq n$. Then from the relation (\ref{eqn:fusion_rule-fk}) with $p=n$ we obtain the decomposition (\ref{eqn:decomp-Xk}) with $p = n+1$. Since $f_p$ is simple for $0 \leq p \leq n$, $\langle \rho^{n+1},\rho^{n+1} \rangle = \langle f_{n+1},f_{n+1} \rangle + \sum_{j=0}^{\lfloor (n+1)/2 \rfloor} \left[ \begin{array}{c} n+1 \\ j \end{array} \right] - 1 = \langle f_{n+1},f_{n+1} \rangle + c_n - 1$, where $c_n = C^{2n}_n/(n+1)$ is the $n^{\mathrm{th}}$ Catalan number, which gives the dimension of $TL_n$ in the generic case. Thus we see that $\langle f_{n+1},f_{n+1} \rangle = 1$, so that $f_{n+1}$ is indeed simple.

In the generic case, $\delta \geq 2$, the Temperley-Lieb category $TL$ is semisimple, that is, every projection is a direct sum of simple projections, and for any pair of non-isomorphic simple projections $p_1$, $p_2$ we have $\langle p_1,p_2 \rangle = 0$.
We recover the infinite Dynkin diagram $A_{\infty}$, where vertices are labeled by the projections $f_i$ and edges represent tensoring by $\rho$.

In the non-generic case, $\delta = [2]_q < 2$, where $q$ is a $k+2^{\mathrm{th}}$ root of unity, we have $[k+2]_q = 0$. Then $\mathrm{tr}(\mathfrak{f}_{k+1}) = [k+2]_q = 0$. Thus the negligible morphisms are those in the unique proper tensor ideal in the Temperley-Lieb category generated by $\mathfrak{f}_{k+1}$ \cite{goodman/wenzl:2003}. The quotient $TL^{(k)} := TL/ \langle \mathfrak{f}_{k+1} \rangle$ is semisimple with simple objects $f_0 = \mathrm{id}_0, f_1 = \rho, f_2, \ldots, f_{k}$ which satisfy the fusion rules (\ref{eqn:fusion_rule-fk}) for $p < k$, and $f_{k} \otimes \rho \cong f_{k-1}$.
Thus we recover the Dynkin diagram $A_{k+1}$, where the vertices are labeled by the projections $f_i$ and edges represent tensoring by $\rho$.

\subsection{$SU(2)$ module categories} \label{sect:nimreps}

In this section we describe $SU(2)$ module categories in terms of preprojective algebras. Then in the subsequent Section \ref{sect:BGPA-GHJ} we relate this to  braided subfactors using the Goodman-de la Harpe-Jones construction \cite{goodman/de_la_harpe/jones:1989, evans/kawahigashi:1998, bockenhauer/evans/kawahigashi:2000} and its manifestation in the bipartite graph planar algebra construction \cite{jones:2000}.

As usual ${}_N \mathcal{X}_N$ will be a braided system of endomorphisms on a factor $N$, and $N \subset M$ will be a braided inclusion with classifying graph $\mathcal{G} = G_{\rho}$, of $ADE$ type, arising from the nimrep $G$ of ${}_N \mathcal{X}_N$ acting on ${}_N \mathcal{X}_M$.
Denote by $\mathcal{G}_0$, $\mathcal{G}_1$ the vertices, edges of $\mathcal{G}$ respectively.
The graph $\mathcal{G}$ is directed, and for every edge $a \in \mathcal{G}_1$ from vertex $i$ to $j$, there is a unique reverse edge $\widetilde{a} \in \mathcal{G}_1$ from $j$ to $i$.
As described in Section \ref{sect:subfactors}, the irreducible sectors in ${}_N \mathcal{X}_N$ label the vertices of the Dynkin diagram $A_{k+1}$, and the edges of $A_{k+1}$ represent multiplication by the fundamental generator, and the irreducible endomorphisms satisfy the same fusion rules as the projections $f_i$ in the category $TL^{(k)}$. We will use $\lambda_i$ to denote endomorphisms in ${}_N \mathcal{X}_N$, whilst $f_i$ will denote the abstract object in the category $TL^{(k)}$.

Semisimple module categories over $\mathcal{C}_q$ (where $\mathcal{C}_q = TL$ for $q = \pm 1$ or $q$ not a root of unity, and $\mathcal{C}_q = TL^{(k)}$ when $q$ is an $k+2^{\mathrm{th}}$ root of unity) where classified in \cite{etingof/ostrik:2004}: A semisimple $\mathcal{C}_q$-module category $\mathcal{D}$ is abelian, and is equivalent as an abelian category to the category $\mathcal{M}_I$ of $I$-graded vector spaces, where $I$ are the (isomorphism classes of) simple objects of $\mathcal{D}$. The structure of a $\mathcal{C}_q$-module category on $\mathcal{M}_I$ is the same as a tensor functor $F$ from $\mathcal{C}_q$ to $\mathrm{Fun}(\mathcal{M}_I,\mathcal{M}_I) \cong \mathcal{M}_{I \times I}$, the category of additive functors from $\mathcal{M}_I$ to itself.
Thus the module category $\mathcal{D} = {}_N \mathcal{X}_M$ gives rise to a monoidal functor $F$ from the Temperley-Lieb category $TL^{(k)}$ to $\mathrm{Fun}({}_N \mathcal{X}_M,{}_N \mathcal{X}_M)$, given by
\begin{equation} \label{eqn:functorF}
F(f) = \bigoplus_{i,j \in \mathcal{G}_0} G_{\lambda}(i,j) \, \mathbb{C}_{i,j},
\end{equation}
where $\lambda = \lambda_p$ is an irreducible endomorphism in ${}_N\mathcal{X}_N$ identified naturally with the Jones-Wenzl projection $f = f_p$.
The $\mathbb{C}_{i,j}$ are 1-dimensional $R$-$R$ bimodules, where $R = (\mathbb{C}\mathcal{G})_0$. The category of $R$-$R$ bimodules has a natural monoidal structure given by $\otimes_R$, or more explicitly, $E^{(1)} \otimes_R E^{(2)} = \bigoplus_{i,k \in \mathcal{G}_0} (E^{(1)} \otimes_R E^{(2)})_{i,k}$, where $(E^{(1)} \otimes_R E^{(2)})_{i,k} = \bigoplus_{j \in \mathcal{G}_0}(E^{(1)}_{i,j} \otimes E^{(2)}_{j,k})$ for all $R$-$R$ bimodules $E^{(r)} = \bigoplus_{i,j \in \mathcal{G}_0} E^{(r)}_{i,j}$, $r=1,2$.
Then we have $R$-$R$ bimodules $F(\rho) = \bigoplus_{i,j \in \mathcal{G}_0} \Delta_{\mathcal{G}}(i,j) \, \mathbb{C}_{i,j} = (\mathbb{C}\mathcal{G})_1$ and $F(\rho^m) = \otimes^m (\mathbb{C}\mathcal{G})_1 = (\mathbb{C}\mathcal{G})_m$.

The set of all edges $a$ form a basis for $(\mathbb{C}\mathcal{G})_1$.
The functor $F$ is defined on the morphisms of $TL$ by specifying annihilation and creation operators $c$, $c^{\ast}$ respectively \cite[Section 4]{ocneanu:2000i}:
\begin{eqnarray}
c(ab) & = & \delta_{b,\widetilde{a}} \frac{\sqrt{\mu_{r(a)}}}{\sqrt{\mu_{s(a)}}} \, s(a), \label{eqn:annihilation} \\
c^{\ast}(i) & = & \sum_{a \in \mathcal{G}_1:r(a)=i} \frac{\sqrt{\mu_{s(a)}}}{\sqrt{\mu_i}} \widetilde{a} a, \label{eqn:creation}
\end{eqnarray}
where $(\mu_j)_j$ is the Perron-Frobenius eigenvector for the Perron-Frobenius eigenvalue $\delta$ of $\mathcal{G}$. Then we set $F\left( \includegraphics[width=5mm]{fig_nakayama-cup} \right) = c^{\ast}$, $F\left( \includegraphics[width=5mm]{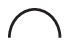} \right) = c$.

Let $\Sigma$ be the graded algebra $\Sigma = \bigoplus_{p=0}^{\infty} F(f_p)$, where the $p^{\mathrm{th}}$ graded part is $\Sigma_p = F(f_p) = F(\lambda_p)$. The multiplication $\mu$ is defined by $\mu_{p,l} = F(\mathfrak{f}_{p+l}): \Sigma_p \otimes_R \Sigma_l \rightarrow \Sigma_{p+l}$.

Preprojective algebras associated to graphs were introduced in \cite{gelfand/ponomarev:1979}, and it was shown that they are finite dimensional if and only if the graphs are of $ADE$ type. They have since found many other applications, including to Kleinian singularities \cite{crawley-boevey/holland:1998} and to Nakajima's quiver varieties \cite{malkin/ostrik/vybornov:2006}.
The preprojective algebra of $\mathcal{G}$ is the graded algebra defined by $\Pi = \mathbb{C}\mathcal{G}/\langle \mathrm{Im} \left( F \left( \includegraphics[width=5mm]{fig_nakayama-cup} \right) \right) \rangle$, where $\langle \mathrm{Im} \left( F \left( \includegraphics[width=5mm]{fig_nakayama-cup} \right) \right) \rangle \subset \mathbb{C}\mathcal{G}$ is the two-sided ideal generated by the image of the creation operators $\includegraphics[width=5mm]{fig_nakayama-cup}$ in $\mathbb{C}\mathcal{G}$. Its $p^{\mathrm{th}}$ graded part is $\Pi_p = (\mathbb{C}\mathcal{G})_p/\langle \mathrm{Im} \left( F \left( \includegraphics[width=5mm]{fig_nakayama-cup} \right) \right) \rangle_p$, where $\langle \mathrm{Im} \left( F \left( \includegraphics[width=5mm]{fig_nakayama-cup} \right) \right) \rangle_p$ is the restriction of $\langle \mathrm{Im} \left( F \left( \includegraphics[width=5mm]{fig_nakayama-cup} \right) \right) \rangle$ to $(\mathbb{C}\mathcal{G})_p$, which is equal to $\sum_{i=1}^{p-1}\mathrm{Im}(F(\mathfrak{e}_i))$, the linear span of the images in $\mathbb{C}\mathcal{G}_p$ of the morphisms $\mathfrak{e}_i = \mathrm{id}_{e_i}$.

Now $(\mathbb{C}\mathcal{G})_p/\sum_{i=1}^{p-1}\mathrm{Im}(F(\mathfrak{e}_i)) = (\mathbb{C}\mathcal{G})_p/\mathrm{Im}(F(\mathfrak{e}_1) \vee \cdots \vee F(\mathfrak{e}_{p-1}))$ is isomorphic to $\mathrm{ker}(F(\mathfrak{e}_1) \vee \cdots \vee F(\mathfrak{e}_{p-1})) = \mathrm{Im}(1-F(\mathfrak{e}_1) \vee \cdots \vee F(\mathfrak{e}_{p-1})) = \mathrm{Im}(F(\mathfrak{f}_p))$.
These are the essential paths $\mathrm{EssPath}_p = \mathrm{ker}(F(\mathfrak{e}_1) \vee \cdots \vee F(\mathfrak{e}_{p-1}))$ of Ocneanu \cite{ocneanu:2000i}, and we have
$$\Pi_p = (\mathbb{C}\mathcal{G})_p/\sum_{i=1}^{p-1}\mathrm{Im}(F(\mathfrak{e}_i)) \cong \mathrm{Im}(F(\mathfrak{f}_p)) = F(f_p) = \Sigma_p.$$
The isomorphism $\varphi: \Sigma_p \rightarrow \Pi_p$ is given by the natural inclusion of $\Sigma_p$ in $(\mathbb{C}\mathcal{G})_p$, then passing to the quotient $(\mathbb{C}\mathcal{G})_p/\langle \mathrm{Im} \left( F \left( \includegraphics[width=5mm]{fig_nakayama-cup} \right) \right) \rangle_p = \Pi_p$.
That $\varphi$ is an isomorphism as algebras is seen as follows, see \cite[Proposition 5.5.6]{cooper:2007}. Since the quotient map $\pi: \mathbb{C}\mathcal{G} \rightarrow \Pi$ is an algebra homomorphism, the multiplication of the images of $\Sigma_r$ and $\Sigma_s$ in $\Pi$ is equal to the multiplication of $F(f_r)$ and $F(f_s)$ in $\mathbb{C}\mathcal{G}$ and then taking the quotient. Now the image of $\mu_{r,s}(\Sigma_r \otimes \Sigma_s)$ in $\mathbb{C}\mathcal{G}$ is $F(f_{r+s})$, which is equal to the image of the multiple of $F(f_r)$ and $F(f_s)$ in $\mathbb{C}\mathcal{G}$ under $F(\mathfrak{f}_{r+s})$, since $\mathfrak{f}_p (\rho^i \otimes f_{p'} \otimes \rho^{p-p'-i}) = f_p$ for any $0 \leq i \leq p-p'$, $p' \leq p$. Since $\mathfrak{f}_{r+s} = \mathrm{id}_{{\rho}^{r+s}} + \phi$, where $\phi$ is a linear combination of $\mathfrak{e}_i$, we see that $\mathrm{Im}(F(\phi)) \subset \langle \mathrm{Im} \left( F \left( \includegraphics[width=5mm]{fig_nakayama-cup} \right) \right) \rangle$ so that $\pi \circ F(\mathfrak{f}_{r+s}) = \pi$. Thus $\varphi$ is an algebra homomorphism.

Applying the functor $F$ to the construction in Section \ref{sect:SU(2)categorical} we obtain the identification $F(f_p) \otimes_R F(\rho) \cong F(f_{p-1}) \oplus F(f_{p+1})$, which yields
\begin{equation} \label{eqn:SU(2)exact_seq-pre}
\Sigma_p \otimes_R (\mathbb{C}\mathcal{G})_1 \cong \Sigma_{p-1} \oplus \Sigma_{p+1}.
\end{equation}
Let $\Lambda$ be the graded coalgebra $\Lambda = F(f_0) \oplus F(f_1) \oplus F(f_0) = (\mathbb{C}\mathcal{G})_0 \oplus (\mathbb{C}\mathcal{G})_1 \oplus (\mathbb{C}\mathcal{G})_0$ with comultiplication $\Delta$, where $\Delta_{1,1}: \Lambda_2 \rightarrow \Lambda_1 \otimes_S \Lambda_1$ is given by $\Delta_{1,1} = F \left( \includegraphics[width=5mm]{fig_nakayama-cup} \right)$ and the other comultiplications are trivial \cite{cooper:2007}.
Let $\delta = [2]_q$ and suppose $[m]_q \neq 0$ for all $m \leq n$ for some $n \in \mathbb{N}$.
Then for all $p \leq n-2$ we obtain the following exact sequence:
\begin{equation} \label{eqn:SU(2)exact_seq}
0 \longrightarrow \Sigma_{p-1} \otimes_R \Lambda_2 \longrightarrow \Sigma_p \otimes_R \Lambda_1 \longrightarrow \Sigma_{p+1} \otimes_R \Lambda_0 \longrightarrow 0,
\end{equation}
where the connecting maps are given by the Koszul differential, the composite map $d_{p,l} = (\mu_{p,1} \otimes 1) \circ (1 \otimes \Delta_{1,l-1}) : \Sigma_p \otimes_R \Lambda_l \rightarrow \Sigma_p \otimes_R \Lambda_1 \otimes_R \Lambda_{l-1} = \Sigma_p \otimes_R \Sigma_1 \otimes_R \Lambda_{l-1} \rightarrow \Sigma_{p+1} \otimes_R \Lambda_{l-1}$.

Thus for $\delta = [2]_q < 2$, where $q$ is a $k+2^{\mathrm{th}}$ root of unity, we have $\Sigma = \bigoplus_{p=0}^{k} F(f_p)$, since $f_p = 0$ in $TL^{(k)}$ for $p \geq k+1$. This means that
$\mathrm{Im}(F(\mathfrak{e}_1) \vee \cdots \vee F(\mathfrak{e}_{k})) = (\mathbb{C}\mathcal{G})_{k+1}$.
The short exact sequence (\ref{eqn:SU(2)exact_seq}) degenerates for $p=k$ to give $0 \longrightarrow \Sigma_{k-1} \otimes_R \Lambda_2 \longrightarrow \Sigma_{k} \otimes_R \Lambda_1 \longrightarrow 0$, and we see that the pair $(\Pi,\Lambda)$ is almost Koszul, in the sense of \cite{brenner/butler/king:2002}, where the preprojective algebra $\Pi$ is a $(k,2)$-Koszul algebra \cite[Corollary 4.3]{brenner/butler/king:2002} \cite[Corollary 5.6.16]{cooper:2007}.

In the generic case, $\delta \geq 2$, there is an analogous pair $(\Pi,\Lambda)$ which is Koszul \cite{brenner/butler/king:2002}, where
Koszul duality is a generalization of the duality between symmetric and antisymmetric algebras.

\subsection{Bipartite graph planar algebras and the GHJ construction} \label{sect:BGPA-GHJ}

Jones \cite{jones:2000} introduced the graph planar algebra construction for a bipartite graph $\mathcal{G}$. We will show that the functor $F$ defined in (\ref{eqn:functorF}) recovers this bipartite graph planar algebra construction.

The planar algebra $P^{\mathcal{G}}$ of a finite bipartite graph $\mathcal{G}$, introduced in \cite{jones:2000}, is the path algebra on $\mathcal{G}$ where paths may start at any of the even vertices of $\mathcal{G}$, and where the $m^{\mathrm{th}}$ graded part $P^{\mathcal{G}}_m$ is given by all pairs of paths of length $m$ on $\mathcal{G}$ which start at the same even vertex and have the same end vertex.
Let $\mathcal{P}$ be the set of tangles in a disc with an even number of vertices on its outer disc, and a finite number of internal discs, each with an even number of vertices, such that each vertex is an endpoint of a string. Internal discs with $2m$ vertices on their boundary are labeled by  elements of $P^{\mathcal{G}}_m$. The presenting map $Z:\mathcal{P} \rightarrow P^{\mathcal{G}}$ is defined uniquely \cite[Theorem 3.1]{jones:2000}, up to isotopy, by first isotoping the strings of a tangle $T$ with internal discs in such a way that $T$ may be divided into horizontal strips where in each strip only cups, caps, internal discs or through strings appear. Then each cup, cap is given by the local operators (\ref{eqn:creation}), (\ref{eqn:annihilation}) respectively, which operate on the elements of $P^{\mathcal{G}}$ inserted in the internal discs, and the outer boundary of the tangle $T$ yields an element of $P^{\mathcal{G}}$.
The planar algebra $P^{\mathcal{G}}$ is a planar $\ast$-algebra, with $\ast$-operation defined on matrix units by $(x_1,x_2)^{\ast} = (x_2,x_1)$. The $\ast$-structure on a tangle $T$ is given by reflecting about a horizontal line which bisects $T$, and replacing every label of $T$ by its adjoint. The tangles $E_i$ in Figure \ref{Fig:fig_Ei} are thus self-adjoint.

We have a tower of algebras $P^{\mathcal{G}}_0 \subset P^{\mathcal{G}}_1 \subset P^{\mathcal{G}}_2 \subset \cdots \;$, where the inclusion $P^{\mathcal{G}}_m \subset P^{\mathcal{G}}_{m+1}$ is given by the graph $\mathcal{G}$.
There is a positive definite inner product defined from the trace on $P^{\mathcal{G}}$.
We have the inclusion $P^{\emptyset}_m := Z(\mathcal{V}_m) \subset P^{\mathcal{G}}_m$ for each $m$, and a double sequence
$$\begin{array}{ccccccc}
P^{\emptyset}_0 & \subset & P^{\emptyset}_1 & \subset & P^{\emptyset}_2 & \subset & \cdots \\
\cap & & \cap & & \cap & & \\
P^{\mathcal{G}}_0 & \subset & P^{\mathcal{G}}_1 & \subset & P^{\mathcal{G}}_2 & \subset & \cdots
\end{array}$$
Then $P^{\emptyset} := Z(\mathcal{V})$ is the embedding of the Temperley-Lieb algebra into the path algebra of $\mathcal{G}$, which is used to construct the Goodman-de la Harpe-Jones (GHJ) subfactors \cite{goodman/de_la_harpe/jones:1989}.
Let $\overline{P^{\mathcal{G}}}$ denote the von Neumann algebra GNS-completion of $P^{\mathcal{G}}$ with respect to the trace.
Then for $q = Z(\ast_{\mathcal{G}})$ the minimal projection in $P^{\mathcal{G}}_0$ corresponding to the distinguished vertex $\ast_{\mathcal{G}}$ of $\mathcal{G}$ with lowest Perron-Frobenius weight, we have an inclusion $q \overline{P^{\emptyset}} \subset q \overline{P^{\mathcal{G}}} q$ which gives the Goodman-de la Harpe-Jones subfactor $N_A \subset N_{\mathcal{G}}$, where $N_A' \cap N_{\mathcal{G}} = q P^{\mathcal{G}}_0 q = \mathbb{C}$.

The limit of the sequence of inclusions
$$\begin{array}{ccccccc}
P^{\mathcal{G}}_0 & \subset & P^{\mathcal{G}}_1 & \subset & P^{\mathcal{G}}_2 & \subset & \cdots \\
\cap & & \cap & & \cap & & \\
P^{\mathcal{G}}_1 & \subset & P^{\mathcal{G}}_2 & \subset & P^{\mathcal{G}}_3 & \subset & \cdots
\end{array}$$
gives a subfactor $N_{\mathcal{G}} \subset M_{\mathcal{G}}$ in a similar way.

Thus we obtain a commuting square of inclusions \cite{goodman/de_la_harpe/jones:1989}:
\begin{equation} \label{eqn:NA-NG_commuting_square}
\begin{array}{ccc}
N_{A} & \subset & N_{\mathcal{G}} \\
\cap & & \cap \\
M_{A} & \subset & M_{\mathcal{G}}
\end{array}
\end{equation}

This commuting square allows us to compute the dual canonical endomorphism $\theta$ of the GHJ subfactor, from which we construct the nimrep graph $\mathcal{G} = G_{\rho}$ \cite{bockenhauer/evans/kawahigashi:2000}.

Since the functor $F$ in Section \ref{sect:nimreps} is defined by the annihilation and creation operators in (\ref{eqn:annihilation}), (\ref{eqn:creation}), we see that $F$ is equivalent to the presenting map $Z$ above. The embedding $P^A \subset P^{\mathcal{G}}$ is given by the image under $F$ of the morphisms in the Temperley-Lieb category.

\subsection{$SU(3)$ Categorical Approach} \label{sect:SU(3)categorical}

In this section we describe the Verlinde algebra and fusion rules for $SU(3)$ in the diagrammatic and categorical language of Kuperberg spiders \cite{kuperberg:1996, evans/pugh:2009iii, cooper:2007}, namely the $A_2$-Temperley-Lieb category.

Let $q$ be real or a root of unity, so that $\delta = [2]_q$ is real. The $A_2$-Temperley-Lieb algebra is the generalized Temperley-Lieb algebra generated by a family $\{ U_j \}$ of self-adjoint operators which satisfy the relations H1-H3 and the vanishing of the $q$-antisymmetrizer for $SU(3)$ which gives (\ref{eqn:SU(3)q_condition}).

\begin{figure}[bt]
\begin{center}
  \includegraphics[width=40mm]{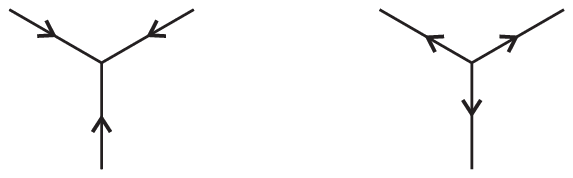}\\
 \caption{$A_2$ webs}\label{fig:A2-webs}
\end{center}
\end{figure}

We call a vertex a source vertex if the string attached to it has orientation away from the vertex. Similarly, a sink vertex will be a vertex where the string attached has orientation towards the vertex.
A string $s$ is a sequence of signs $+$, $-$.
For two (possibly empty) strings $s_1, s_2$, an $A_2$-$s_1,s_2$-tangle $T$ is a tangle on a rectangle with strings $s_1$, $s_2$ along the top, bottom edges respectively, generated by $A_2$ webs (see Figure \ref{fig:A2-webs}) such that every free end of $T$ is attached to a vertex along the top or bottom of the rectangle in a way that respects the orientation of the strings, every vertex has a string attached to it, and the tangle contains no closed loops or elliptic faces. Along the top edge the points $+$ are source vertices and $-$ are sink vertices, while along the bottom edge the roles are reversed.
We define the vector space $\mathcal{V}^{A_2}_{s_1,s_2}$ to be the free vector space over $\mathbb{C}$ with basis $\mathcal{T}^{A_2}_{s_1,s_2}$.

We define $V^{A_2}_{s_1,s_2}$ to be the quotient of $\mathcal{V}^{A_2}_{s_1,s_2}$ by the Kuperberg ideal generated by the Kuperberg relations K1-K3 below. That is, composition in $V^{A_2}_{s_1,s_2}$ is defined as follows. The composition $RS \in V^{A_2}_{s_1,s_3}$ of an $A_2$-$s_1,s_2$-tangle $R$ and an $A_2$-$s_2,s_3$-tangle $S$ is given by gluing $S$ vertically below $R$ such that the vertices at the bottom of $R$ and the top of $S$ coincide, removing these vertices, and isotoping the glued strings if necessary to make them smooth. Any closed loops which may appear are removed, contributing a factor of $\alpha = [3]_q$, as in relation K1 below. Any elliptic faces that appear are removed using relations K2, K3 below. The composition is associative and is extended linearly to elements in $V^{A_2}_{s_1,s_2}$.
\begin{center}
\begin{minipage}[b]{11.5cm}
 \begin{minipage}[t]{3cm}
  \parbox[t]{2cm}{\begin{eqnarray*}\textrm{K1:}\end{eqnarray*}}
 \end{minipage}
 \begin{minipage}[t]{5.5cm}
  \begin{center}
  \mbox{} \\
 \includegraphics[width=20mm]{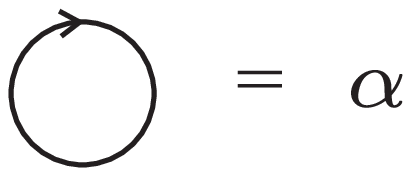}
  \end{center}
 \end{minipage}
 \begin{minipage}[t]{2cm}
  \mbox{} \\
  \parbox[t]{1cm}{}
 \end{minipage}
\end{minipage}
\begin{minipage}[b]{11.5cm}
 \begin{minipage}[t]{3cm}
  \parbox[t]{2cm}{\begin{eqnarray*}\textrm{K2:}\end{eqnarray*}}
 \end{minipage}
 \begin{minipage}[t]{5.5cm}
  \begin{center}
  \mbox{} \\
 \includegraphics[width=23mm]{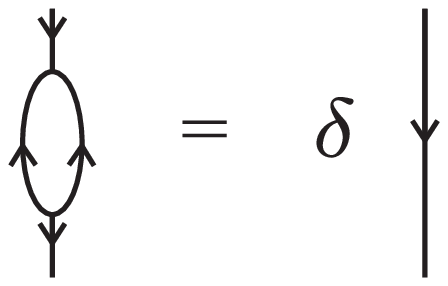}
  \end{center}
 \end{minipage}
 \begin{minipage}[t]{2cm}
  \mbox{} \\
  \parbox[t]{1cm}{}
 \end{minipage}
\end{minipage}
\begin{minipage}[b]{11.5cm}
 \begin{minipage}[t]{3cm}
  \parbox[t]{2cm}{\begin{eqnarray*}\textrm{K3:}\end{eqnarray*}}
 \end{minipage}
 \begin{minipage}[t]{5.5cm}
  \begin{center}
  \mbox{} \\
 \includegraphics[width=55mm]{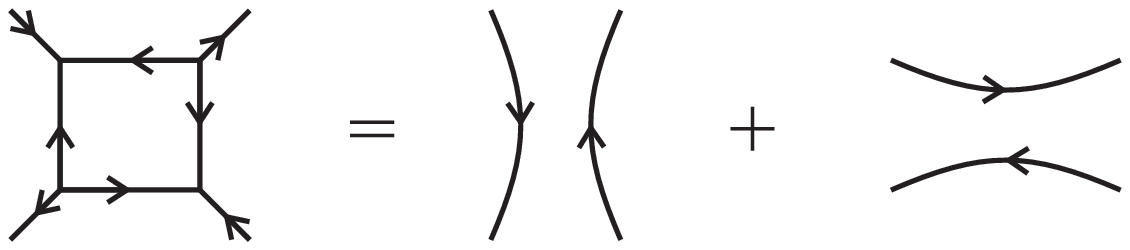}
  \end{center}
 \end{minipage}
 \begin{minipage}[t]{2cm}
  \mbox{} \\
  \parbox[t]{1cm}{}
 \end{minipage}
\end{minipage}
\end{center}

There is a braiding on $\mathcal{V}^{A_2}_{s_1,s_2}$, defined locally by the following linear combinations of local diagrams in $\mathcal{V}^{A_2}_{s_1,s_2}$ (see \cite{kuperberg:1996, suciu:1997}), for any $q \in \mathbb{C}$:
\begin{center}
\begin{minipage}[b]{16cm}
 \begin{minipage}[t]{4.5cm}
  \mbox{} \\
  \parbox[t]{1cm}{}
 \end{minipage}
 \begin{minipage}[t]{5.5cm}
  \mbox{} \\
   \includegraphics[width=55mm]{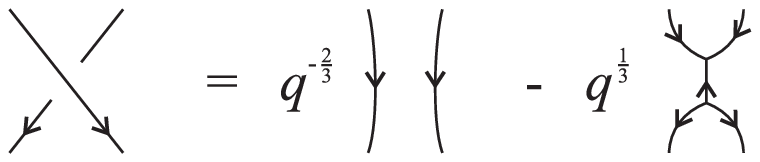}
 \end{minipage}
 \begin{minipage}[t]{5.5cm}
  \hfill
  \parbox[t]{2cm}{\begin{eqnarray}\label{braiding1}\end{eqnarray}}
 \end{minipage}
\end{minipage}
\end{center}

\noindent The braiding satisfies type II and type III Reidemeister moves, and a braiding fusion
relation \cite[Equations (8), (9)]{evans/pugh:2009iii}, provided $\delta = [2]_q$ and $\alpha = [3]_q$.

Thus it is sufficient to work over $V^{A_2}_{(m,n),(m',n')} := V^{A_2}_{+^m -^n, +^{m'} -^{n'}}$, where $+^k -^l$ is the string of $k$ signs $+$ followed by $l$ signs $-$, since, for any arbitrary string $s$ with $m$ signs $+$ and $n$ signs $-$ and string $s'$ with $m'$ signs $+$ and $n'$ signs $-$, there is an isomorphism $\iota$ between $V^{A_2}_{s,s'}$ and $V^{A_2}_{(m,n), (m',n')}$ given by using the braiding to permute the order of the signs in $s$ to $+^m -^n$, and the inverse braiding to permute the order of the signs in $s'$ to $+^{m'} -^{n'}$.

\begin{figure}[tb]
\begin{center}
  \includegraphics[width=40mm]{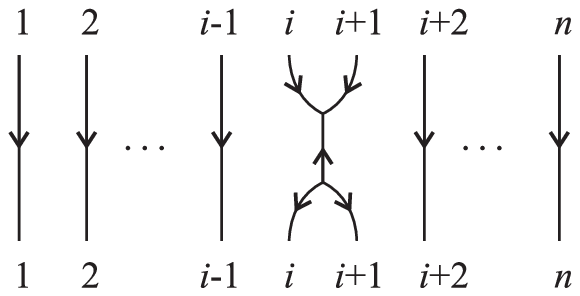}\\
 \caption{The tangle $W_i \in \mathcal{V}^{A_2}_{m}$.}\label{fig:W_i}
\end{center}
\end{figure}

A diagrammatic representation of the Hecke algebra for $SU(3)$ is as follows: Let $W_i \in \mathcal{V}^{A_2}_{m} := \mathcal{V}^{A_2}_{(m,0),(m,0)}$ be the tangle illustrated in Figure \ref{fig:W_i}.
A $\ast$-operation can be defined on $\mathcal{V}^{A_2}_{m}$, where for an $m$-tangle $T \in \mathcal{T}^{A_2}_{m}$, $T^{\ast}$ is the $m$-tangle obtained by reflecting $T$ about a horizontal line halfway between the top and bottom vertices of the tangle, and reversing the orientations on every string. Then $\ast$ on $\mathcal{V}^{A_2}_{m}$ is the conjugate linear extension of $\ast$ on $\mathcal{T}^{A_2}_{m}$. For $\delta \in \mathbb{R}$ (so $q \in \mathbb{R}$ or $q$ a root of unity), the $\ast$-operation leaves the Kuperberg ideal invariant due to the symmetry of the relations K1-K3.
For $m \in \mathbb{N} \cup \{ 0 \}$ we define the algebra $A_2\textrm{-}TL_m$ to be $\mathrm{alg}(\mathrm{id}_m, w_i|i = 1, \ldots, m-1)$, where $\mathrm{id}_m$ is the identity diagram given by $m$ vertical strings, and $w_i$ is the image of $W_i \in \mathcal{V}^{A_2}_{m}$ in the quotient space $V^{A_2}_m := V^{A_2}_{(m,0),(m,0)}$. The $w_i$'s in $A_2\textrm{-}TL_m$ are clearly self-adjoint, and satisfy the relations H1-H3 and (\ref{eqn:SU(3)q_condition}) \cite{evans/pugh:2009iii}.

Diagrammatically, the generalized Jones-Wenzl projections $f_{(m,0)}$ of Section \ref{sect:subfactors} are given as follows: $f_{(1,0)}$ is given by a single vertical string in $\mathcal{T}^{A_2}_{+, -}$, whilst $f_{(m,0)} \in V^{A_2}_{(m,0),(m,0)}$ is defined inductively by \cite[(2.1.0)-(2.1.2)]{suciu:1997}:
\begin{center}
\begin{minipage}[b]{16cm}
 \begin{minipage}[t]{1cm}
  \hfill
  \parbox[t]{0.5cm}{}
 \end{minipage}
 \begin{minipage}[t]{13cm}
  \begin{center}
  \mbox{} \\
   \includegraphics[width=80mm]{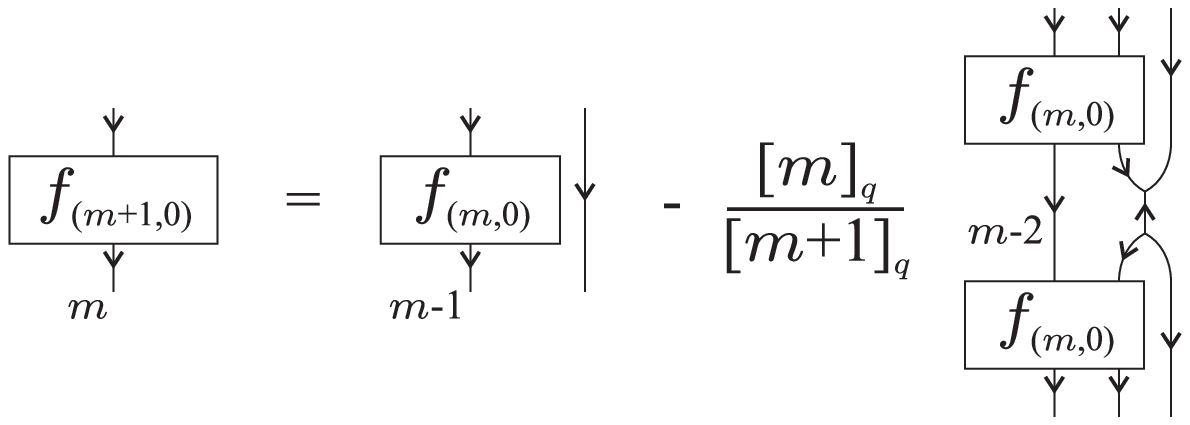}
  \end{center}
 \end{minipage}
 \begin{minipage}[t]{1.5cm}
  \hfill
  \vspace{2mm} \parbox[t]{1.5cm}{\begin{eqnarray}\label{eqn:f(k,0)}\end{eqnarray}}
 \end{minipage}
\end{minipage}
\end{center}

The generalized Jones-Wenzl projections $f_{(m,n)} \in V^{A_2}_{(m,n)} := V^{A_2}_{(m,n),(m,n)}$ are defined inductively by
\cite[(2.1.7)]{suciu:1997}:
\begin{center}
\begin{minipage}[b]{16cm}
 \begin{minipage}[t]{1cm}
  \hfill
  \parbox[t]{0.5cm}{}
 \end{minipage}
 \begin{minipage}[t]{13cm}
  \begin{center}
  \mbox{} \\
   \includegraphics[width=120mm]{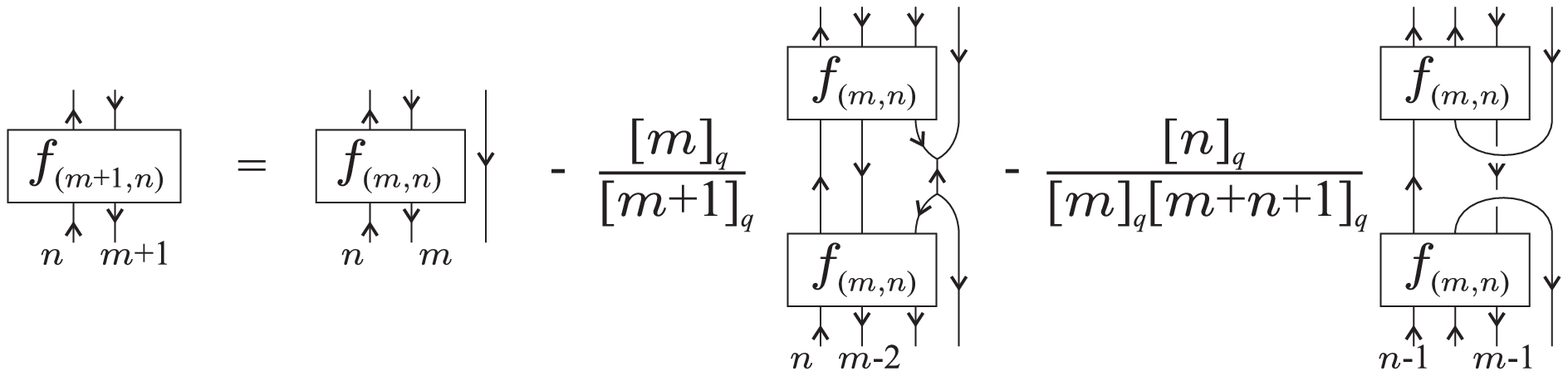}
  \end{center}
 \end{minipage}
 \begin{minipage}[t]{1.5cm}
  \hfill
  \vspace{2mm} \parbox[t]{1.5cm}{\begin{eqnarray}\label{eqn:f(m,n)}\end{eqnarray}}
 \end{minipage}
\end{minipage}
\end{center}

We will define the $A_2$-Temperley-Lieb category by $A_2\textrm{-}TL = \mathrm{Mat}(C^{A_2})$, where $C^{A_2}$ is the tensor category whose objects are projections in $V^{A_2}_{(m,n)}$, and whose morphisms are $\mathrm{Hom}(p_1,p_2) = p_2 V^{A_2}_{(m_2,n_2),(m_1,n_1)} p_1$, for projections $p_i \in V^{A_2}_{(m_i,n_i)}$, $i=1,2$. We write $A_2\textrm{-}TL_{(m,n)} = V^{A_2}_{(m,n)}$, and $\rho$, $\overline{\rho}$ for the identity projections in $A_2\textrm{-}TL_{(1,0)}$, $A_2\textrm{-}TL_{(0,1)}$ respectively consisting of a single string with orientation downwards, upwards respectively. Then the identity diagram in $A_2\textrm{-}TL_{(m,n)}$, given by $m+n$ vertical strings where the first $m$ strings have downwards orientation and the next $n$ have upwards orientation, is expressed as $\rho^m \overline{\rho}^n$.
We have $\mathrm{dim}(A_2\textrm{-}TL_{(0,0)}) = \mathrm{dim}(A_2\textrm{-}TL_{(1,0)}) = \mathrm{dim}(A_2\textrm{-}TL_{(0,1)}) = 1$ and $A_2\textrm{-}TL_{(0,0)}$, $A_2\textrm{-}TL_{(1,0)}$ and $A_2\textrm{-}TL_{(0,1)}$ have simple projections $f_{(0,0)}$ (the empty diagram), $f_{(1,0)} = \rho$ and $f_{(0,1)} = \overline{\rho}$ respectively. Moving to level 2, that is, $A_2\textrm{-}TL_{(k,l)}$ such that $k+l=2$, consider first $A_2\textrm{-}TL_{(2,0)}$. The identity diagram $\rho^2$ is a projection but is not simple, since $\langle \rho^2,\rho^2 \rangle = 2$. One of these morphisms is the identity diagram, the other is $\mathfrak{w}_1 = \mathrm{id}_{w_1}$. Now $u_1 = \delta^{-1} w_1$ is a projection, isomorphic to $f_{(0,1)}$, as can be seen from the isomorphisms $\psi: f_{(0,1)} \rightarrow u_1$ and $\psi^{\ast}: u_1 \rightarrow f_{(0,1)}$, where $\psi = \sqrt{\delta}^{-1}$ \includegraphics[width=5mm]{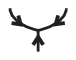}.
Then we have $\psi\psi^{\ast} = \mathfrak{u}_1 = \mathrm{id}_{u_1}$, $\psi^{\ast}\psi = \mathfrak{f}_{(0,1)} = \mathrm{id}_{f_{(0,1)}}$. Since $\langle f_{(0,1)},\rho^2 \rangle = 1$, where the morphism is given by $\,$ \includegraphics[width=5mm]{fig_nakayama-Yfork} $\,$, and $\langle f',\rho^2 \rangle = 0$ for $f' = f_{(0,0)}, f_{(1,0)}$ (by parity), we have the decomposition $\rho^2 = f_{(0,1)} \oplus f_{(2,0)}$, where $f_{(2,0)}$ is a simple projection in $A_2\textrm{-}TL_{(2,0)}$. Similarly, $\overline{\rho}^2 = f_{(1,0)} \oplus f_{(0,2)}$, where $f_{(0,2)}$ is a simple projection in $A_2\textrm{-}TL_{(0,2)}$. Finally at level 2 consider $A_2\textrm{-}TL_{(1,1)}$. The simple projection $e_1 := \alpha^{-1} E_1 \in \mathrm{Hom}(\rho \otimes \overline{\rho},\rho \otimes \overline{\rho})$ is isomorphic to $f_{(0,0)}$, as in the $SU(2)$ case (see Section \ref{sect:SU(2)categorical}). Since $\langle f_{(0,0)},\rho \otimes \overline{\rho} \rangle = 1$, where the morphism is given by $\,$ \includegraphics[width=5mm]{fig_nakayama-cup} $\,$, and $\langle f',\rho \otimes \overline{\rho} \rangle = 0$ for $f' = f_{(1,0)}, f_{(0,1)}$ (by parity), we have the decomposition $\rho \otimes \overline{\rho} = f_{(0,0)} \oplus f_{(1,1)}$, where $f_{(1,1)}$ is a simple projection in $A_2\textrm{-}TL_{(1,1)}$.

In the same way, we obtain we obtain at each level that $\rho^m \overline{\rho}^n$ is a linear combination of $f_{(i,j)}$ for $i,j \geq 0$, $0 \leq i+j < m+n$ such that $i-j \equiv m-n \textrm{ mod } 3$, plus a new projection $f_{(m,n)}$, which turns out to be simple.
The morphisms $\mathfrak{f}_{(p,l)} = \mathrm{id}_{f_{(p,l)}}$ are generalized Jones-Wenzl projections which satisfy the recursion relations (\ref{eqn:f(k,0)}) and (\ref{eqn:f(m,n)}).
These satisfy the properties \cite{suciu:1997}:
\begin{center}
\includegraphics[width=80mm]{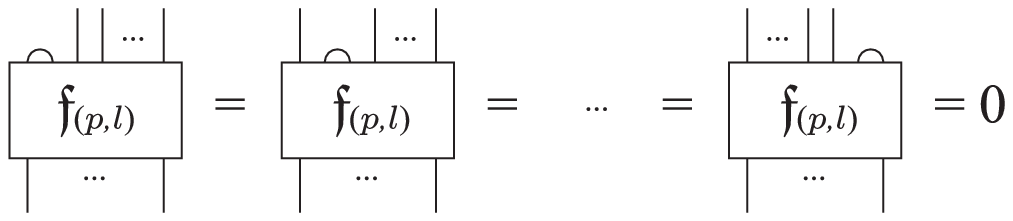} \\
\vspace{5mm}
\includegraphics[width=80mm]{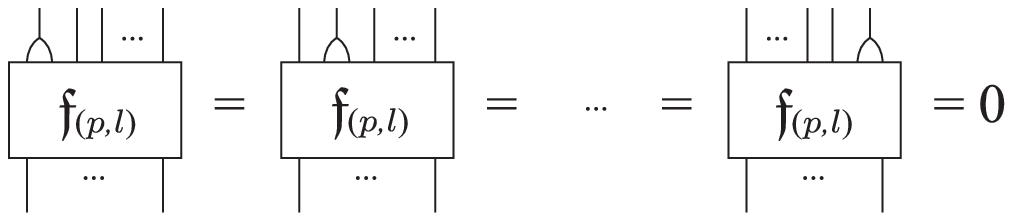} \\
\vspace{5mm}
\includegraphics[width=65mm]{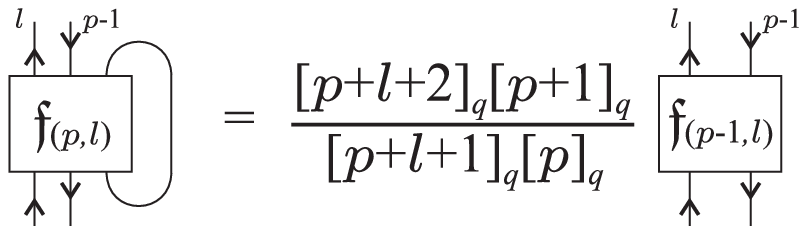}
$$\mathrm{tr}(\mathfrak{f}_{(p,l)}) = \frac{[p+1]_q[l+1]_q[p+l+2]_q}{[2]_q},$$
\end{center}
where $\mathrm{tr}$ is defined as in Section \ref{sect:SU(2)categorical}.
For $p \geq p'$ and $l \geq l'$, these generalized Jones-Wenzl projections also satisfy the property $\mathfrak{f}_{(p,l)} (\mathrm{id}_{\rho^i \overline{\rho}^j} \otimes \mathfrak{f}_{(p',l')} \otimes \mathrm{id}_{\rho^{p-p'-i} \overline{\rho}^{l-l'-j}}) = \mathfrak{f}_{(p,l)} = (\mathrm{id}_{\rho^i \overline{\rho}^j} \otimes \mathfrak{f}_{(p',l')} \otimes \mathrm{id}_{\rho^{p-p'-i} \overline{\rho}^{l-l'-j}}) \mathfrak{f}_{(p,l)}$, for any $0 \leq i \leq p-p'$, $0 \leq j \leq l-l'$. This property also holds if we conjugate either $\mathfrak{f}_{(p,l)}$ or $\mathrm{id}_{\rho^i \overline{\rho}^j} \otimes \mathfrak{f}_{(p',l')} \otimes \mathrm{id}_{\rho^{p-p'-i} \overline{\rho}^{l-l'-j}}$ by any braiding.

We have the relations
\begin{equation} \label{eqn:fusion_rule-f(k,l)}
f_{(p,l)} \otimes \rho \cong f_{(p,l-1)} \oplus f_{(p-1,l+1)} \oplus f_{(p+1,l)},
\end{equation}
given by the isomorphisms $\psi: f_{(p,l)} \otimes \rho \rightarrow f_{(p,l-1)} \oplus f_{(p-1,l+1)} \oplus f_{(p+1,l)}$, and $\psi^{-1} = \psi^{\ast}$, where
\begin{equation} \label{eqn:SU(3)-psi}
\psi = \left( \begin{array}{c} \includegraphics[width=50mm]{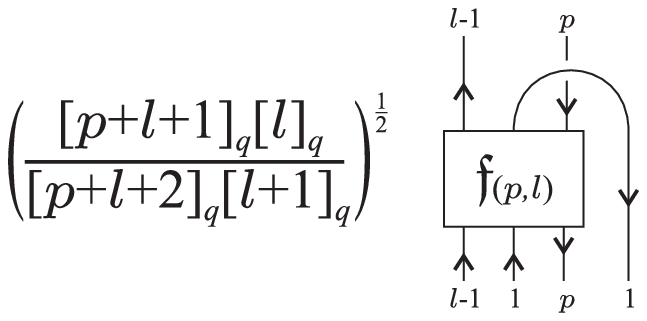} \\ \includegraphics[width=35mm]{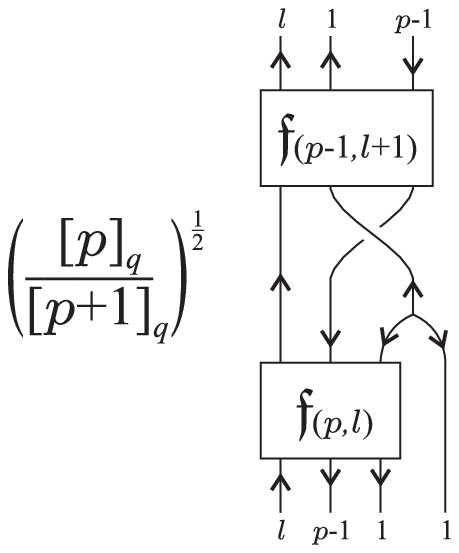} \\ \mathfrak{f}_{(p+1,l)} \end{array} \right).
\end{equation}
Then it easy to check using the above properties and (\ref{eqn:f(m,n)}) that $\psi\psi^{\ast} = \mathfrak{f}_{(p,l-1)} \oplus \mathfrak{f}_{(p-1,l+1)} \oplus \mathfrak{f}_{(p+1,l)}$ and $\psi^{\ast}\psi = \mathfrak{f}_{(p,l)} \otimes \mathrm{id}_{\rho}$.
Similarly we have
\begin{equation} \label{eqn:fusion_rule-f(k,l)conjugate}
f_{(p,l)} \otimes \overline{\rho} \cong f_{(p-1,l)} \oplus f_{(p+1,l-1)} \oplus f_{(p,l+1)},
\end{equation}
so the $f_{(p,l)}$ satisfy the fusion rules for $SU(3)$, given in (\ref{eqn:SU(3)fusion_rules}).
Suppose that we have obtained simple projections $f_{(p,l)}$ for $0 \leq p+l \leq n$, and we obtain a new projection $f_{(j,n-j+1)}$ at level $n+1$ as above, for some $0 \leq j \leq n+1$. From the relation (\ref{eqn:fusion_rule-f(k,l)}) with $p=j$, $l=n-j$ and from dimension considerations we see that $\langle f_{(j,n-j+1)},f_{(j,n-j+1)} \rangle = 1$, so that $f_{(j,n-j+1)}$ is indeed simple.

In the generic case, $\delta \geq 2$, the $A_2$-Temperley-Lieb category $A_2\textrm{-}TL$ is semisimple and for any pair of non-isomorphic simple projections $p_1$, $p_2$ we have $\langle p_1,p_2 \rangle = 0$.
We recover the infinite graph $\mathcal{A}^{(\infty)}$, where the vertices are labeled by the projections $f_{(p,l)}$ and the edges represent tensoring by $\rho$.

In the non-generic case where $q$ is a $k+3^{\mathrm{th}}$ root of unity, we have $\mathrm{tr}(\mathfrak{f}_{(p,l)}) = 0$ for $p+l=k+1$. By (\ref{eqn:f(m,n)}), $\mathfrak{f}_{(p',l')} = 0$ for all $p',l' \geq k+2$ if $\mathfrak{f}_{(p,l)} = 0$ for $p+l=k+1$. Thus the negligible morphisms are the ideal $\langle \mathfrak{f}_{(p,l)} | p+l=k+1 \rangle$ generated by $\mathfrak{f}_{(p,l)}$ such that $p+l=k+1$.
The quotient $A_2\textrm{-}TL^{(k)} := A_2\textrm{-}TL/ \langle \mathfrak{f}_{(p,l)} | p+l=k+1 \rangle$ is semisimple with simple objects $f_{(p,l)}$, $p,l \geq 0$ such that $p+l \leq k$ which satisfy the fusion rules (\ref{eqn:fusion_rule-f(k,l)}) and (\ref{eqn:fusion_rule-f(k,l)conjugate}), where $f_{(p',l')}$ is understood to be zero if $p'<0$, $l'<0$ or $p'+l' \geq k+1$.
Thus we recover the graph $\mathcal{A}^{(k+3)}$, where the vertices are labeled by the projections $f_{(p,l)}$ and the edges represent tensoring by $\rho$.

\subsection{$SU(3)$ module categories} \label{sect:nimrepsSU(3)}

In this section we describe $SU(3)$ module categories in terms of certain algebras of paths. Then in the subsequent Section \ref{sect:A2GPA-GHJ} we relate this to braided subfactors using the $SU(3)$ Goodman-de la Harpe-Jones construction \cite{evans/pugh:2009ii} and its manifestation in the $SU(3)$-graph planar algebra construction \cite{evans/pugh:2009iv}.

As usual ${}_N \mathcal{X}_{N} = \{ \lambda_{(p,l)} | \; 0 \leq p,l,p+l \leq k < \infty \}$ will be braided system of endomorphisms of $SU(3)_k$ on a factor $N$, and $N \subset M$ will be a braided inclusion with classifying graph $\mathcal{G} = G_{\rho}$, of $\mathcal{ADE}$ type, arising from the nimrep $G$ of ${}_N \mathcal{X}_N$ acting on ${}_N \mathcal{X}_M$. Then the module category gives rise to a monoidal functor $F$ from the $A_2$-Temperley-Lieb category $A_2\textrm{-}TL^{(k)}$ to $\mathrm{Fun}({}_N \mathcal{X}_M,{}_N \mathcal{X}_M)$, where $F$ is given by (\ref{eqn:functorF}), where now $\lambda = \lambda_{(p,l)}$ is an irreducible endomorphism in ${}_N \mathcal{X}_N$ identified with the generalized Jones-Wenzl projections $f = f_{(p,l)}$.
We denote by $\mathcal{G}^{\mathrm{op}}$ the opposite graph of $\mathcal{G}$ obtained by reversing the orientation of every edge of $\mathcal{G}$.
Then we have that $F(\rho^m \overline{\rho}^n)$ is the $R$-$R$ bimodule with basis given by all paths of length $m+n$ on
$\mathcal{G}$, $\mathcal{G}^{\mathrm{op}}$, where the first $m$ edges are on $\mathcal{G}$ and the last $n$ edges are on $\mathcal{G}^{\mathrm{op}}$, where $R=(\mathbb{C}\mathcal{G})_0$.
In particular $F(\rho^m) = (\mathbb{C}\mathcal{G})_m$.

\begin{figure}[tb]
\begin{center}
   \includegraphics[width=70mm]{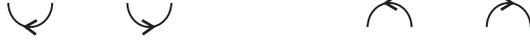}
\caption{left and right cups; left and right caps} \label{fig:cups&caps}
\end{center}
\end{figure}

\begin{figure}[tb]
\begin{center}
   \includegraphics[width=70mm]{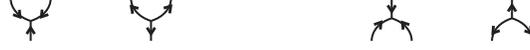}
\caption{incoming and outgoing Y-forks; incoming and outgoing inverted Y-forks} \label{fig:Y-forks}
\end{center}
\end{figure}

If $a \in \mathcal{G}_1$ is an edge on $\mathcal{G}$, we denote by $\widetilde{a} \in \mathcal{G}^{\mathrm{op}}_1$ the corresponding edge with opposite orientation on $\mathcal{G}^{\mathrm{op}}$.
We define annihilation operators $c_l$, $c_r$ by:
\begin{equation} \label{eqn:annihilation-lr}
c_l(a\widetilde{b}) = \delta_{s(a),s(b)} \frac{\sqrt{\mu_{r(a)}}}{\sqrt{\mu_{s(a)}}} \, s(a), \qquad c_r(\widetilde{b}a) = \delta_{r(a),r(b)} \frac{\sqrt{\mu_{s(a)}}}{\sqrt{\mu_{r(a)}}} \, r(a),
\end{equation}
and creation operators $c_l^{\ast}$, $c_r^{\ast}$ as their adjoints, where $a$ is an edge on $\mathcal{G}$ and $\widetilde{b}$ an edge on $\mathcal{G}^{\mathrm{op}}$, and $(\mu_j)_j$ is the Perron-Frobenius eigenvector for the Perron-Frobenius eigenvalue $\alpha$ of $\mathcal{G}$. Define the following fork operators $\curlyvee$, $\overline{\curlyvee}$ by:
\begin{eqnarray}
\curlyvee(\widetilde{a}) & = & \frac{1}{\sqrt{\mu_{s(a)} \mu_{r(a)}}} \sum_{b_1,b_2} W(\triangle_{s(a),r(a),r(b_1)}^{(a,b_1,b_2)}) b_1 b_2, \label{eqn:Yfork(in)} \\
\overline{\curlyvee}(a) & = & \frac{1}{\sqrt{\mu_{s(a)} \mu_{r(a)}}} \sum_{b_1,b_2} \overline{W(\triangle_{s(a),r(a),r(b_1)}^{(a,b_1,b_2)})} \widetilde{b_1} \widetilde{b_2}, \label{eqn:Yfork(out)}
\end{eqnarray}
where $\triangle_{s(a),r(a),r(b_1)}^{(a,b_1,b_2)}$ denotes the closed loop of length 3 on $\mathcal{G}$ along the edges $a$, $b_1$ and $b_2$, and $W(\triangle)$ are the Ocneanu cells on $\mathcal{G}$ constructed in \cite{evans/pugh:2009i}. We also define $\curlywedge = \overline{\curlyvee}^{\ast}$ and $\overline{\curlywedge} = \curlyvee^{\ast}$.
Then the functor $F$ is defined on the morphisms of $A_2\textrm{-}TL$ by assigning the following operators to the morphisms given in Figures  \ref{fig:cups&caps} and \ref{fig:Y-forks}: to the left, right caps the annihilation operators $c_l$, $c_r$ respectively, given by (\ref{eqn:annihilation-lr}), to the left, right cups the creation operators $c_l^{\ast}$, $c_r^{\ast}$ respectively, to the incoming, outgoing Y-forks the operators $\curlyvee$, $\overline{\curlyvee}$ respectively given in (\ref{eqn:Yfork(in)}), (\ref{eqn:Yfork(out)}), and to the incoming, outgoing inverted Y-forks the operators $\curlywedge$, $\overline{\curlywedge}$ respectively.

Let $\Sigma$ be the graded algebra $\Sigma = \bigoplus_{p=0}^{\infty} F(f_{(p,0)})$, where the $p^{\mathrm{th}}$ graded part is $\Sigma_p = F(f_{(p,0)}) = F(\lambda_{(p,0)})$. The multiplication $\mu$ is defined by $\mu_{p,l} = F(\mathfrak{f}_{(p+l,0)}): \Sigma_p \otimes_R \Sigma_l \rightarrow \Sigma_{p+l}$, where $\mathfrak{f}_{(p,l)} = \mathrm{id}_{f_{(p,l)}}$.

We define a graded algebra $\Pi$ by $\Pi = \mathbb{C}\mathcal{G}/\langle \mathrm{Im} \big( F \big( \includegraphics[width=5mm]{fig_nakayama-Yfork} \big) \big) \rangle$, where $\langle \mathrm{Im} \big( F \big( \includegraphics[width=5mm]{fig_nakayama-Yfork} \big) \big) \rangle \subset \mathbb{C}\mathcal{G}$ is the two-sided ideal generated by the image of the operators $\includegraphics[width=5mm]{fig_nakayama-Yfork}$ in $\mathbb{C}\mathcal{G}$. Its $p^{\mathrm{th}}$ graded part is $\Pi_p = (\mathbb{C}\mathcal{G})_p/\langle \mathrm{Im} \big( F \big( \includegraphics[width=5mm]{fig_nakayama-Yfork} \big) \big) \rangle_p$, where $\langle \mathrm{Im} \big( F \big( \includegraphics[width=5mm]{fig_nakayama-Yfork} \big) \big) \rangle_p$ is the restriction of $\langle \mathrm{Im} \big( F \big( \includegraphics[width=5mm]{fig_nakayama-Yfork} \big) \big) \rangle$ to $(\mathbb{C}\mathcal{G})_p$, which is equal to $\sum_{i=1}^{p-1}\mathrm{Im}(F(\mathfrak{U}_i))$, the union of the images on $\mathbb{C}\mathcal{G}_p$ of the morphisms $\mathfrak{U}_i = \mathrm{id}_{U_i}$.

The quotient $(\mathbb{C}\mathcal{G})_p/\sum_{i=1}^{p-1}\mathrm{Im}(F(\mathfrak{U}_i)) = (\mathbb{C}\mathcal{G})_p/\mathrm{Im}(F(\mathfrak{U}_1) \vee \cdots \vee F(\mathfrak{U}_{p-1}))$ is isomorphic to $\mathrm{ker}(F(\mathfrak{U}_1) \vee \cdots \vee F(\mathfrak{U}_{p-1}))$.
Clearly $\mathrm{ker}(F(\mathfrak{U}_1) \vee \cdots \vee F(\mathfrak{U}_{p-1})) \supset \mathrm{Im}(F(\mathfrak{f}_{(p,0)}))$ since $\mathfrak{U}_i \mathfrak{f}_{(p,0)} = 0$ for $i=1,\ldots,p-1$. For $a \in \mathrm{ker}(F(\mathfrak{U}_1) \vee \cdots \vee F(\mathfrak{U}_{p-1}))$, $a = a \cdot F(\mathfrak{f}_{(p,0)})$ since the only term in $\mathfrak{f}_{(p,0)}$ which does not contain a $\mathfrak{U}_i$ is the identity, which has coefficient 1. Thus $a \in \mathrm{Im}(F(\mathfrak{f}_{(p,0)}))$ so $\mathrm{ker}(F(\mathfrak{U}_1) \vee \cdots \vee F(\mathfrak{U}_{p-1})) = \mathrm{Im}(F(\mathfrak{f}_{(p,0)}))$.

Then we have
\begin{equation} \label{eqn:Pi=Sigma}
\Pi_p = (\mathbb{C}\mathcal{G})_p/\sum_{i=1}^{p-1}\mathrm{Im}(F(\mathfrak{U}_i)) \cong \mathrm{Im}(F(\mathfrak{f}_{(p,0)})) = F(f_{(p,0)}) = \Sigma_p.
\end{equation}
The isomorphism is given by the natural inclusion of $\Sigma_p$ in $(\mathbb{C}\mathcal{G})_p$, then passing to the quotient $(\mathbb{C}\mathcal{G})_p/\langle \mathrm{Im} \big( F \big( \includegraphics[width=5mm]{fig_nakayama-Yfork} \big) \big) \rangle_p = \Pi_p$.
That this map is an isomorphism as algebras follows by an analogous argument to that in the $SU(2)$ case \cite[Theorem 7.3.5]{cooper:2007}.

Applying the functor $F$ to the construction in Section \ref{sect:SU(3)categorical} we obtain the identifications $F(f_{(k,l)}) \otimes_R F(\rho) \cong F(f_{(p,l-1)}) \oplus F(f_{(p-1,l+1)}) \oplus F(f_{(p+1,l)})$ and $F(f_{(p-2,l)}) \otimes_R F(\overline{\rho}) \cong F(f_{(p-2,l)}) \oplus F(f_{(p,l-1)}) \oplus F(f_{(p-1,l+1)})$, which yield when $l=0$:
\begin{eqnarray}
\Sigma_p \otimes_R (\mathbb{C}\mathcal{G})_1 \cong F(f_{(p-1,1)}) \oplus \Sigma_{p+1}, \label{eqn:SU(3)exact_seq-pre1} \\
\Sigma_{p-1} \otimes_R (\mathbb{C}\mathcal{G}^{\mathrm{op}})_1 \cong \Sigma_{p-2} \oplus F(f_{(p-1,1)}). \label{eqn:SU(3)exact_seq-pre2}
\end{eqnarray}
Let $\Lambda$ be the graded coalgebra $\Lambda = F(f_{(0,0)}) \oplus F(f_{(1,0)}) \oplus F(f_{(0,1)}) \oplus F(f_{(0,0)}) = (\mathbb{C}\mathcal{G})_0 \oplus (\mathbb{C}\mathcal{G})_1 \oplus (\mathbb{C}\mathcal{G}^{\mathrm{op}})_1 \oplus (\mathbb{C}\mathcal{G})_0$. The comultiplication is $\Delta$ is given by $\Delta_{1,1} = F \big( \includegraphics[width=5mm]{fig_nakayama-Yfork} \big): \Lambda_2 \rightarrow \Lambda_1 \otimes_R \Lambda_1$, $\Delta_{1,2} = F \big( \includegraphics[width=4mm]{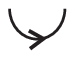} \big): \Lambda_3 \rightarrow \Lambda_1 \otimes_R \Lambda_2$, $\Delta_{2,1} = F \big( \includegraphics[width=4mm]{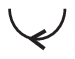} \big): \Lambda_3 \rightarrow \Lambda_2 \otimes_R \Lambda_1$, and the other comultiplications are trivial.
Let $\delta = [2]_q$ and suppose $[m]_q \neq 0$ for all $m \leq n$ for some $n \in \mathbb{N}$.
Then for all $p \leq n-3$ we obtain the following exact sequences:
\begin{eqnarray}
0 \longrightarrow F(f_{(p-1,1)}) \stackrel{F(\psi_2)}{\longrightarrow} \Sigma_p \otimes_R \Lambda_1 \stackrel{d}{\longrightarrow} \Sigma_{p+1} \otimes_R \Lambda_0 \longrightarrow 0, \label{eqn:SU(3)exact_seq-1} \\
0 \longrightarrow \Sigma_{p-2} \otimes_R \Lambda_3 \stackrel{d}{\longrightarrow} \Sigma_{p-1} \otimes_R \Lambda_2 \stackrel{F(\mathfrak{f}_{(p-1,1)})}{\longrightarrow} F(f_{(p-1,1)}) \longrightarrow 0, \label{eqn:SU(3)exact_seq-2}
\end{eqnarray}
where $\psi = (\psi_i)$ is the isomorphism given in (\ref{eqn:SU(3)-psi}) and $d$ is the Koszul differential defined in Section \ref{sect:nimreps}.
These sequences can be combined to give another exact sequence:
\begin{equation} \label{eqn:SU(3)exact_seq}
0 \longrightarrow \Sigma_{p-2} \otimes_R \Lambda_3 \longrightarrow \Sigma_{p-1} \otimes_R \Lambda_2 \longrightarrow \Sigma_p \otimes_R \Lambda_1 \longrightarrow \Sigma_{p+1} \otimes_R \Lambda_0 \longrightarrow 0,
\end{equation}
where now all the connecting maps are given by the Koszul differential $d$.

Thus when $q$ is a $k+3^{\mathrm{th}}$ root of unity, we have $\Sigma = \bigoplus_{p=0}^{k} F(f_{(p,0)})$, since $f_{(p,0)} = 0$ in $A_2\textrm{-}TL^{(k)}$ for $p \geq k+1$. This means that $\mathrm{Im}(F(\mathfrak{U}_1) \vee \cdots \vee F(\mathfrak{U}_{k})) = (\mathbb{C}\mathcal{G})_{k+1}$.
The exact sequence (\ref{eqn:SU(3)exact_seq}) degenerates for $p=k$ to give $0 \longrightarrow \Sigma_{k-2} \otimes_R \Lambda_3 \longrightarrow \Sigma_{k-1} \otimes_R \Lambda_2 \longrightarrow \Sigma_{k} \otimes_R \Lambda_1 \longrightarrow 0$, and for $p=k+1$ it degenerates to give $0 \longrightarrow \Sigma_{k-1} \otimes_R \Lambda_3 \longrightarrow \Sigma_{k} \otimes_R \Lambda_2 \longrightarrow 0$.
Then we see that the pair $(\Pi,\Lambda)$ is almost Koszul, where the algebra $\Pi$ is a $(k,3)$-Koszul algebra \cite[Corollary 7.4.19]{cooper:2007}.

In the generic case, $\delta \geq 2$, there is an analogous pair $(\Pi,\Lambda)$ which is Koszul \cite[Corollary 7.3.9]{cooper:2007}.

The isomorphism between the two algebras $\Pi$ and $\Sigma$ is a key ingredient in the determination of the Hilbert series in Section \ref{sect:Hilbert-SU(3)}.

\subsection{$SU(3)$-graph planar algebras and the $SU(3)$-GHJ construction} \label{sect:A2GPA-GHJ}

In \cite{evans/pugh:2009iv} we introduced the $A_2$-graph planar algebra construction for an $SU(3)$ $\mathcal{ADE}$ graph $\mathcal{G}$. The $A_2$-graph planar algebra $P^{\mathcal{G}}$ of an $SU(3)$ $\mathcal{ADE}$ graph $\mathcal{G}$ is the path algebra on $\mathcal{G}$ and $\mathcal{G}^{\mathrm{op}}$. We will show that the functor $F$ defined in Section \ref{sect:nimrepsSU(3)} recovers this $A_2$-graph planar algebra construction.

The presenting map $Z:\mathcal{P} \rightarrow P^{\mathcal{G}}$ is defined uniquely \cite[Theorem 5.1]{evans/pugh:2009iv}, up to isotopy, by first isotoping the strings of $T$ in such a way that the diagram $T$ may be divided into horizontal strips so that each horizontal strip only contains the following elements: a (left or right) cup, a (left or right) cap, an (incoming or outgoing) Y-fork, or an (incoming or outgoing) inverted Y-fork, see Figures \ref{fig:cups&caps} and \ref{fig:Y-forks}. Then $Z$ assigns to the left, right caps the annihilation operators $c_l$, $c_r$ respectively, given by (\ref{eqn:annihilation-lr}), to the left, right cups the creation operators $c_l^{\ast}$, $c_r^{\ast}$ respectively, to the incoming, outgoing Y-forks the operators $\curlyvee$, $\overline{\curlyvee}$ respectively given in (\ref{eqn:Yfork(in)}), (\ref{eqn:Yfork(out)}), and to the incoming, outgoing inverted Y-forks the operators $\curlywedge$, $\overline{\curlywedge}$ respectively.

We have a tower of algebras $P^{\mathcal{G}}_{0,0} \subset P^{\mathcal{G}}_{0,1} \subset P^{\mathcal{G}}_{0,2} \subset \cdots \;$, where the inclusion $P^{\mathcal{G}}_m \subset P^{\mathcal{G}}_{m+1}$ is given by the $m$-$(m+1)$ part of the graph $\mathcal{G}$.
There is a positive definite inner product defined from the trace on $P^{\mathcal{G}}$.
We have the inclusion $P^{\emptyset}_{0,m} := Z(\mathcal{V}^{A_2}_m) \subset P^{\mathcal{G}}_{0,m}$ for each $m$, and we have a double sequence
$$\begin{array}{ccccccc}
P^{\emptyset}_{0,0} & \subset & P^{\emptyset}_{0,1} & \subset & P^{\emptyset}_{0,2} & \subset & \cdots \\
\cap & & \cap & & \cap & & \\
P^{\mathcal{G}}_{0,0} & \subset & P^{\mathcal{G}}_{0,1} & \subset & P^{\mathcal{G}}_{0,2} & \subset & \cdots
\end{array}$$
Then $P^{\emptyset} := Z(\mathcal{V}^{A_2})$ is the embedding of the $A_2$-Temperley-Lieb algebra into the path algebra of $\mathcal{G}$, which is used to construct the $A_2$-Goodman-de la Harpe-Jones subfactors \cite{evans/pugh:2009ii}.
Let $\overline{P^{\mathcal{G}}}$ denote the GNS-completion of $P^{\mathcal{G}}$ with respect to the trace.
Then for $q = Z(\ast_{\mathcal{G}})$ the minimal projection in $P^{\mathcal{G}}_0$ corresponding to the vertex $\ast_{\mathcal{G}}$ of $\mathcal{G}$, we have an inclusion $q \overline{P^{\emptyset}} \subset q \overline{P^{\mathcal{G}}} q$ which gives the $A_2$-Goodman-de la Harpe-Jones subfactor $N_{\mathcal{A}} \subset N_{\mathcal{G}}$, where $N_{\mathcal{A}}' \cap N_{\mathcal{G}} = q P^{\mathcal{G}}_{0,0} q = \mathbb{C}$ and the sequence $\{ q P^{\emptyset}_{0,m} \subset P^{\mathcal{G}}_{0,m} \}_m$ is a periodic sequence of commuting squares of period 3, in the sense of Wenzl in \cite{wenzl:1988}.

Thus we obtain a commuting square of inclusions as in (\ref{eqn:NA-NG_commuting_square}), which allows us to compute the dual canonical endomorphism $\theta$ of the $A_2$-Goodman-de la Harpe-Jones subfactor, from which we construct the nimrep graph $\mathcal{G} = G_{\rho}$ \cite{evans/pugh:2009ii}.

Since the functor $F$ in Section \ref{sect:nimrepsSU(3)} is defined by the annihilation and creation operators given by (\ref{eqn:annihilation-lr}), and the incoming, outgoing (inverted) Y-fork operators given in (\ref{eqn:Yfork(in)}), (\ref{eqn:Yfork(out)}), we see that $F$ is equivalent to the presenting map $Z$ above. The embedding $P^A \subset P^{\mathcal{G}}$ is given by the image under $F$ of the morphisms in the $A_2$-Temperley-Lieb category.

\section{Hilbert series of the almost Calabi-Yau algebras} \label{sect:Hilbert-SU(3)}

In Section \ref{sect:almostCYalg} we introduce an algebra $A(\mathcal{G},W)$ associated to a finite graph $\mathcal{G}$ (an $SU(3)$ $\mathcal{ADE}$ graph or the McKay graph $\mathcal{G}_{\Gamma}$ of finite subgroup $\Gamma \subset SU(3)$) which carries a cell system $W$. In the case where $\mathcal{G}$ is an $SU(3)$ $\mathcal{ADE}$ graph, these algebras are called almost Calabi-Yau algebras and are shown to be isomorphic to the almost Koszul algebras $\Pi$ of Section \ref{sect:nimrepsSU(3)}. We will determine a formula for the Hilbert series which counts the dimensions of these algebras in Section \ref{sect:Hilbert_series}.

\subsection{The algebras $A(\mathcal{G},W)$} \label{sect:almostCYalg}

In this section we introduce the algebra $A(\mathcal{G},W)$ associated to a finite graph $\mathcal{G}$ which carries a cell system $W$.

For any finite directed graph $\mathcal{G}$, let $[\mathbb{C}\mathcal{G}, \mathbb{C}\mathcal{G}]$ denote the subspace of $\mathbb{C}\mathcal{G}$ spanned by all commutators of the form $xy - yx$, for $x,y \in \mathbb{C}\mathcal{G}$. If $x,y$ are paths in $\mathbb{C}\mathcal{G}$ such that $r(x) = s(y)$ but $r(y) \neq s(x)$, then $xy - yx = xy$, so in the quotient $\mathbb{C}\mathcal{G} / [\mathbb{C}\mathcal{G}, \mathbb{C}\mathcal{G}]$ the path $xy$ will be zero. Then any non-cyclic path, i.e. any path $x$ such that $r(x) \neq s(x)$, will be zero in $\mathbb{C}\mathcal{G} / [\mathbb{C}\mathcal{G}, \mathbb{C}\mathcal{G}]$. If $x = a_1 a_2 \cdots a_k$ is a cyclic path in $\mathbb{C}\mathcal{G}$, then $a_1 a_2 \cdots a_k - a_k a_1 \cdots a_{k-1} = 0$ in $\mathbb{C}\mathcal{G} / [\mathbb{C}\mathcal{G}, \mathbb{C}\mathcal{G}]$, so $a_1 a_2 \cdots a_k$ is identified with $a_k a_1 \cdots a_{k-1}$. Similarly, $x = a_1 a_2 \cdots a_k$ is identified with every cyclic permutation of the edges $a_j$, $j=1,\ldots,k$. So the commutator quotient $\mathbb{C}\mathcal{G} / [\mathbb{C}\mathcal{G}, \mathbb{C}\mathcal{G}]$ may be identified, up to cyclic permutation of the arrows, with the vector space spanned by cyclic paths in $\mathcal{G}$.
One defines a derivation $\partial_a : \mathbb{C}\mathcal{G} / [\mathbb{C}\mathcal{G}, \mathbb{C}\mathcal{G}] \rightarrow \mathbb{C}\mathcal{G}$ by
$\partial_{a} (a_1 \cdots a_n) = \sum_{j} a_{j+1} \cdots a_n a_1 \cdots a_{j-1}$,
where the summation is over all indices $j$ such that $a_j = a$.
Then for a potential $\Phi \in \mathbb{C}\mathcal{G} / [\mathbb{C}\mathcal{G}, \mathbb{C}\mathcal{G}]$, which is some linear combination of cyclic paths in $\mathcal{G}$, we define the algebra
$$A(\mathbb{C}\mathcal{G}, \Phi) = \mathbb{C}\mathcal{G} / \langle \rho_a \rangle,$$
which is the quotient of the path algebra by the two-sided ideal generated by the relations $\rho_a = \partial_a \Phi \in \mathbb{C}\mathcal{G}$, for all edges $a$ of $\mathcal{G}$.
If $\Phi$ is homogeneous, we define the Hilbert series $H_A$ for $A(\mathbb{C}\mathcal{G}, \Phi)$ as $H_A(t) = \sum_{p=0}^{\infty} H_{ji}^p t^p$, where the $H_{ji}^p$ are matrices which count the dimension of the subspace $\{ i x j | \; x \in A(\mathbb{C}\mathcal{G}, \Phi)_p \}$, where $A(\mathbb{C}\mathcal{G}, \Phi)_p$ is the subspace of $A(\mathbb{C}\mathcal{G}, \Phi)$ of all paths of length $p$, and $i,j \in A(\mathbb{C}\mathcal{G}, \Phi)_0$.

Suppose $A(\mathbb{C} \mathcal{G}, \Phi)$ is a Calabi-Yau algebra of dimension $d = 3$ and that $\textrm{deg} \; \Phi = 3$, that is, $\Phi$ is a linear combination of cyclic paths of length 3 on $\mathcal{G}$. Then $H_A(t)$ is given by (\ref{eqn:H(t)-CYd}) \cite[Theorem 4.6]{bocklandt:2008}.

For $\mathcal{G}$ an $SU(3)$ $\mathcal{ADE}$ graph or the McKay graph $\mathcal{G}_{\Gamma}$ of a finite subgroup $\Gamma \subset SU(3)$, we define a homogeneous potential $\Phi$ by \cite[Remark 4.5.7]{ginzburg:2006}:
\begin{equation} \label{eqn:potential-Phi}
\Phi = \sum_{a,b,c \in \mathcal{G}_1} W(\triangle^{(a,b,c)}_{s(a),s(b),s(c)}) \cdot \triangle^{(a,b,c)}_{s(a),s(b),s(c)} \quad \in \mathbb{C} \mathcal{G} / [\mathbb{C} \mathcal{G}, \mathbb{C} \mathcal{G}],
\end{equation}
for a cell system $W$ \cite{ocneanu:2000ii}, and we will denote by $A(\mathcal{G},W)$ the algebra $A(\mathbb{C} \mathcal{G}, \Phi)$.

Now let $\mathcal{G}$ be an $SU(3)$ $\mathcal{ADE}$ graph and let $(x_1,x_2) \in \mathrm{End} \left( {}_{(\mathbb{C}\mathcal{G})_0} (\mathbb{C}\mathcal{G})_p {}_{(\mathbb{C}\mathcal{G})_0} \right)$ be matrix units, where $x_1,x_2 \in (\mathbb{C}\mathcal{G})_p$ for $p = 0,1,\ldots \;$, as in Section \ref{sect:Hecke}, which act on $\mathbb{C}\mathcal{G}$ by $(x_1,x_2)y = \delta_{y,x_2} x_1$.
Applying the functor $F$ to the morphisms $\mathfrak{U}_p$ in $A_2\mathrm{-}TL$ we obtain a representation of the Hecke algebra on $\mathbb{C}\mathcal{G}$ given by (c.f. (\ref{eqn:Yfork(in)}) and Figure \ref{fig:W_i}):
\begin{equation} \label{Def:U_k}
F(\mathfrak{U}_p) = \sum_{x, a_i, b} \phi_{s(a_1)}^{-1} \phi_{r(a_2)}^{-1} W(\triangle_{s(b),r(b),r(a_3)}^{(b,a_3,a_4)}) \overline{W(\triangle_{s(b),r(b),r(a_1)}^{(b,a_1,a_2)})} \, (x a_1 a_2, x a_3 a_4),
\end{equation}
where the summation is over all paths $x$ of length $p-1$ and edges $a_i, b$ of $\mathcal{G}$ such that the paths $x a_1 a_2$, $x a_3 a_4$ make sense. These operators were shown to satisfy the relations H1-H3 in \cite{evans/pugh:2009ii}.
Let $\langle \rho_a \rangle_p$ denote the restriction of the ideal $\langle \rho_a \rangle$ in $(\mathbb{C}\mathcal{G})_p$, which is isomorphic to $\sum_{i=1}^{p-1}\mathrm{Im}(F(\mathfrak{U}_i))$. Then $A(\mathcal{G},W)_p \cong (\mathbb{C}\mathcal{G})_p / \sum_{i=1}^{p-1}\mathrm{Im}(F(\mathfrak{U}_i)) = \Pi_p \cong \Sigma_p$, where $\Pi$, $\Sigma$ are the graded algebras defined in Section \ref{sect:nimrepsSU(3)}.

\subsection{Hilbert Series of $A(\mathcal{G},W)$} \label{sect:Hilbert_series}

In this section we give the Hilbert series of the algebra $A(\mathcal{G},W)$ where $\mathcal{G}$ is an $SU(3)$ $\mathcal{ADE}$ graph $\mathcal{G}$ with cell system $W$.
In this case we will call $A(\mathcal{G},W)$ an almost Calabi-Yau algebra.

When $\mathcal{G} = \mathcal{G}_{\Gamma}$ is the McKay graph of a finite subgroup $\Gamma \subset SU(3)$, $A(\mathcal{G},W)$ is a Calabi-Yau algebra of dimension 3 \cite[Theorem 4.4.6]{ginzburg:2006} and its Hilbert series is thus given by (\ref{eqn:H(t)-CYd}).

Let ${}_N \mathcal{X}_{N} = \{ \lambda_{(p,l)} | \; 0 \leq p,q,p+l \leq k < \infty \}$ denote a non-degenerately braided
system of endomorphisms on a type $\mathrm{III}_1$ factor $N$, which is generated by $\rho = \lambda_{(1,0)}$ and its conjugate $\overline{\rho} = \lambda_{(0,1)}$, where the irreducible endomorphisms $\lambda_{(p,l)}$ satisfy the fusion rules of $SU(3)_k$ given in (\ref{eqn:SU(3)fusion_rules}).

Let $N \subset M$ be a braided subfactor with nimrep $G$ which realises the modular invariant $Z_{\mathcal{G}}$ at level $k$, where $\mathcal{G} = G_{\rho}$, and let $I = {}_N \mathcal{X}_M$. The nimrep $G$ sends $\lambda \in {}_N \mathcal{X}_{N}$ to the graph $G_{\lambda}$.
Let $S = \bigoplus_{p=0}^k \lambda_{(p,0)}$, so that $F(S) = \bigoplus_{p=0}^k F(\lambda_{(p,0)}) = \Sigma$. Then since $A(\mathcal{G},W) = \Pi \cong \Sigma$ (see Section \ref{sect:almostCYalg}),
\begin{equation} \label{eqn:F(S)=A}
F(S) \cong A(\mathcal{G},W),
\end{equation}
that is, the algebra $A(\mathcal{G},W)$ is given by the module category associated to the inclusion $N \subset M$.

The algebra $S$ has another description, based on \cite{malkin/ostrik/vybornov:2006} in the context of $SU(2)$, and the proof of Theorem 4.4.6 in \cite{ginzburg:2006}.
The tensor algebra $T\rho = \bigoplus_{j=0}^{\infty} \rho^j$ is the free algebra in the category $A_2\textrm{-}TL$ at level $k$
generated by $\rho$. Under the functor $F$ defined by (\ref{eqn:functorF}), $T\rho$ maps to $F(T\rho) = \mathbb{C} \mathcal{G}$, c.f. Section \ref{sect:nimrepsSU(3)}.
Let $T'$ be the quotient of $T\rho$ by the two-sided ideal $\langle \overline{\rho} \rangle$ generated by $\overline{\rho} \subset \rho^2$. It can be shown inductively for each grade that $T_j' = S_j$. For $j=0,1$, the result is trivial. Since $\rho^2 = \overline{\rho} + \lambda_{(2,0)}$ by (\ref{eqn:SU(3)fusion_rules}), we have $T_2' = \lambda_{(2,0)} = S_2$. Now consider $j=3$. From (\ref{eqn:SU(3)fusion_rules}) we obtain $\rho^3 = \rho\overline{\rho} + \lambda_{(1,1)} + \lambda_{(3,0)}$. Since $\rho\overline{\rho} = \mathrm{id} + \lambda_{(1,1)}$, we see that the ideal $\langle \mathrm{id} + \lambda_{(1,1)} \rangle \subset \langle \overline{\rho} \rangle$. Then $\lambda_{(1,1)} \in \langle \overline{\rho} \rangle$ and $T_3' = \lambda_{(3,0)} = S_3$. The situation for general $j$ is similar. Thus we see that $S$ is the symmetric algebra given by the quotient of the tensor algebra $T\rho$ by $\langle \overline{\rho} \rangle$,
see also \cite[(4.5.1)]{ginzburg:2006}. Let $\pi$, $\gamma$ denote the composite morphisms $\pi: \lambda_{(0,0)} \hookrightarrow \rho^3 \hookrightarrow T\rho$ and $\gamma: \overline{\rho} \hookrightarrow \rho^2 \hookrightarrow T\rho$. Then $\pi$ maps under $F$ to $F(\pi): \mathbb{C}I \rightarrow \mathbb{C}\mathcal{G}$, which sends $\mathbf{1} \in \mathbb{C}I$ to the potential $\Phi$ of (\ref{eqn:potential-Phi}), and $F(\gamma): F(\overline{\rho}) \rightarrow \mathbb{C}\mathcal{G}$ sends the reverse edge $\widetilde{a}$ of $a$ to the relation $\rho_a = \partial_a \Phi$.
Then under $F$, $S = T\rho / \langle \overline{\rho} \rangle$ is mapped to $F(T\rho) / \langle F(\overline{\rho}) \rangle \cong \mathbb{C}\mathcal{G} / \langle F(\gamma(\overline{\rho})) \rangle = \mathbb{C}\mathcal{G} / \langle \rho_a \rangle = A(\mathcal{G},W)$.
Reversing the argument, exactness of $F$ yields $F(S) = F(T\rho) / \langle F(\overline{\rho}) \rangle$ which by the above discussion is isomorphic to $\mathbb{C}\mathcal{G} / \langle F(\gamma(\overline{\rho})) \rangle = \mathbb{C}\mathcal{G} / \langle \rho_a \rangle = A(\mathcal{G},W)$.

The fusion rules (\ref{eqn:SU(3)fusion_rules}) yield the recursion relation
$\lambda_{(j+1,0)} \oplus (\overline{\rho} \otimes \lambda_{(j-1,0)}) = (\rho \otimes \lambda_{(j,0)}) \oplus \lambda_{(j-2,0)}$, $j=2,3,\ldots,k-1$, and so each $\lambda_{(j,0)}$ can be written recursively in terms of the three irreducible endomorphisms $\rho$, $\overline{\rho}$ and $\lambda_{(0,0)} = \mathrm{id}$.
Summing over all $j$, and using $t$ to keep track of the grading, we find that $S = \bigoplus_{p=0}^k \lambda_{(p,0)} t^p$ satisfies
\begin{equation} \label{eqn:pre-Hilbert_eqn}
\left( S \oplus (t^2 \overline{\rho} \otimes S) \oplus t^{k+3} \lambda_{(k,0)} \right) = \lambda_{(0,0)} \oplus (t \rho \otimes S) \oplus t^3 S.
\end{equation}
Then applying $F$ to (\ref{eqn:pre-Hilbert_eqn}) we obtain
the following equation for the Hilbert series of $A$:
$$H_A(t) - t \Delta_{\mathcal{G}} H_A(t) + t^2 \Delta_{\mathcal{G}}^T H_A(t) - t^3 H_A(t) = \mathbf{1} - t^{k+3} P,$$
where $P = G_{\lambda_{(k,0)}}$. From $P^3 = G_{\lambda_{(k,0)}} \otimes G_{\lambda_{(k,0)}} \otimes G_{\lambda_{(k,0)}} = G_{\lambda_{(0,0)}} = \textbf{1}$, we see that the matrix $P$ is an automorphism of the graph $\mathcal{G}$ of order 3.
Since each $G_{\lambda_{(j,0)}}$ can be written recursively in terms of $G_{\rho} = \Delta_{\mathcal{G}}$, $G_{\overline{\rho}} = \Delta_{\mathcal{G}}^T$ and $G_{\lambda_{(0,0)}} = \mathbf{1}$, the permutation $P$ can be determined for each graph $\mathcal{G}$ separately using standard Mathematica computations, and we obtain the following result:

\begin{Thm} \label{thm:SU(3)Hilbert}
Let $H_A(t)$ denote the Hilbert series of $A(\mathcal{G},W)$, for an $SU(3)$ $\mathcal{ADE}$ graph $\mathcal{G}$ with adjacency matrix $\Delta_{\mathcal{G}}$, Coxeter number $h=k+3$ and cell system $W$.
Then
\begin{equation} \label{eqn:Hilbert_Series-SU(3)ADE}
H_A (t) = \frac{1 - P t^h}{1 - \Delta_{\mathcal{G}} t + \Delta_{\mathcal{G}}^T t^2 - t^3},
\end{equation}
where $P$ is the permutation matrix corresponding to a $\mathbb{Z}_3$ symmetry of the graph.
It is the identity for $\mathcal{D}^{(n)}$, $\mathcal{A}^{(n)\ast}$, $n \geq 5$, $\mathcal{E}^{(8)\ast}$, $\mathcal{E}_l^{(12)}$, $l=1,2,4,5$, and $\mathcal{E}^{(24)}$. For the remaining graphs $\mathcal{A}^{(n)}$, $\mathcal{D}^{(n) \ast}$ and $\mathcal{E}^{(8)}$, let $V$ be the permutation matrix corresponding to the clockwise rotation of the graph by $2 \pi /3$. Then
$$ P = \left\{
\begin{array}{cl} V^2 & \mbox{ for } \quad \mathcal{A}^{(n)}, n \geq 4, \\
                  V & \mbox{ for } \quad \mathcal{E}^{(8)}, \\
                  V^{2n} & \mbox{ for } \quad \mathcal{D}^{(n) \ast}, n \geq 5.
\end{array} \right.$$
\end{Thm}

The numerator and denominator in (\ref{eqn:Hilbert_Series-SU(3)ADE}) commute, since any permutation matrix which corresponds to a symmetry of the graph $\mathcal{G}$ commutes with $\Delta_{\mathcal{G}}$ and $\Delta_{\mathcal{G}}^T$.

We warn that we have not yet realised the graph $\mathcal{E}_4^{(12)}$ as the nimrep produced by a subfactor, as we have not been able to construct a cell system on $\mathcal{E}_4^{(12)}$. However, Ocneanu has claimed that the graph $\mathcal{E}_4^{(12)}$ does have a cell system built on it \cite{ocneanu:2000ii}, and hence the above proof would hold in this case also.

In \cite[Proposition 3.14]{brenner/butler/king:2002} the Hilbert series was given for a $(p,q)$-Koszul algebra (or almost Koszul algebra), where the permutation matrix $P$ is equal to the product of the permutation matrices given by the Nakayama permutations for $A$ and its Koszul dual $\Lambda$. It was shown in Section \ref{sect:nimrepsSU(3)} that $A = A(\mathcal{G},W)$ is a $(h-3,3)$-Koszul algebra. Then the Nakayama permutation for the coalgebra $\Lambda$ being trivial is equivalent to Nakayama permutation of $A$ being given by the permutation matrix $P$.

\section{Nakayama automorphism for $SU(3)$ $\mathcal{ADE}$ graphs} \label{sect:nakayama}

When $\mathcal{G}$ is an $SU(3)$ $\mathcal{ADE}$ graph, the dual $A^{\ast} = \mathrm{Hom}(A, \mathbb{C})$ of the algebra $A=A(\mathcal{G},W)$ is identified as an $A$-$A$ bimodule with ${}_1 A_{\beta}$, with standard left action and the right action twisted by an automorphism $\beta$, the Nakayama automorphism. In this section we determine the Nakayama automorphism $\beta$ in Theorem \ref{Thm:A^=1_A_beta}.

The explicit cell systems $W$ computed in \cite{evans/pugh:2009i} and knowledge of the Hilbert series (\ref{eqn:Hilbert_Series-SU(3)ADE}) for $A$ are key ingredients in both these results.
Thus these results are not proven for the $SU(3)$ $\mathcal{ADE}$ graph $\mathcal{E}_4^{(12)}$, since we were not able to compute an explicit cell system $W$ for this graph in \cite{evans/pugh:2009i}. In the remainder of the paper, any reference to an $SU(3)$ $\mathcal{ADE}$ graph will not include the graph $\mathcal{E}_4^{(12)}$.

We begin with some preliminary results, including Proposition \ref{Prop:C=1} whose lengthy proof will be the content of Section \ref{sect:PropC=1}.
This section is based closely on \cite[Section 4.2]{brenner/butler/king:2002}, which is in the context of preprojective algebras of the $ADE$ graphs in $SU(2)$.

Let $A = A(\mathcal{G},W)$ for an $SU(3)$ $\mathcal{ADE}$ graph $\mathcal{G}$ and let $A_p$ denote its $p^\mathrm{th}$ graded part. The edges $a \in \mathcal{G}_1$ are a basis for $A_1$.
With the potential $\Phi$ defined in (\ref{eqn:potential-Phi}), $A$ has a relation $\rho_a$ for each edge $a \in \mathcal{G}_1$ given by
\begin{equation} \label{eqn:relations-A}
\rho_{a} = \sum_{b,b' \in \mathcal{G}_1} W_{s(a),r(a),r(b)}^{(a,b,b')} b b',
\end{equation}
where $W_{s(a),r(a),r(b)}^{(a,b,b')} := W(\triangle_{s(a),r(a),r(b)}^{(a,b,b')})$.

Let $h=k+3$ denote the Coxeter number of $\mathcal{G}$. The image of the endomorphism $\lambda_{(k,0)}$ under the functor $F$ defines a unique permutation $\nu$ of the graph $\mathcal{G}$, which is described as follows.
If the permutation matrix $P$ in Theorem \ref{thm:SU(3)Hilbert} is the identity matrix, then the permutation $\nu$ on the graph $\mathcal{G}$ is just the identity. For the other graphs, the permutation $\nu$ is given on the vertices of $\mathcal{G}$ by the permutation matrix $P$ and on $\mathcal{G}_1$ by the unique permutation on the edges of $\mathcal{G}$ such that $s(\nu(a)) = \nu(s(a))$ and $r(\nu(a)) = \nu(r(a))$ (note that there are no double edges on the graphs $\mathcal{G}$ for which $P$ is non-trivial).

For any vertex $i$ of $\mathcal{G}$, it can be seen from the Hilbert series of $A(\mathcal{G},W)$ that the space $i \cdot A_{h-3} \cdot \nu(i)$ is one-dimensional. For $k$ a vertex for which there is an edge from $k$ to $i$, it can be seen from the Hilbert series of $A(\mathcal{G},W)$ that the space $i \cdot A_{h-4} \cdot \nu(k)$ is one-dimensional, except where there is a double edge from $k$ to $i$ on $\mathcal{G}$, in which case $\mathrm{dim}(i \cdot A_{h-4} \cdot \nu(k)) = 2$. Let $u_{i\nu(i)}$ denote a generator of $i \cdot A_{h-3} \cdot \nu(i)$, and let $\{ v_{i\nu(k)}^{m}, m \in \{ 1,p_{ki} \} \}$ denote a basis for $i \cdot A_{h-4} \cdot \nu(k)$, where $p_{ki}$ denotes the number of edges from $k$ to $i$. For each edge $i \stackrel{a}{\longrightarrow} j$ and each edge $k \stackrel{b}{\longrightarrow} i$ on $\mathcal{G}$ there are non-zero scalars $\lambda_{a}^{(m)}$, $\mu_{b}^{(m')}$ such that
\begin{equation} \label{eqn:I}
\lambda_{a}^{(m)} a v_{j\nu(i)}^{m} = u_{i\nu(i)} = \mu_{b}^{(m')} v_{i\nu(k)}^{m'} \nu(b),
\end{equation}
for any $m \in \{ 1,p_{ij} \}$, $m' \in \{ 1,p_{ki} \}$.

\begin{Def}
We call a cell system $W$ on $\mathcal{G}$ $\nu$-invariant if $W_{i,j,k}^{(a,b,c)} = W_{\nu(i)\nu(j)\nu(k)}^{(\nu(a), \nu(b), \nu(c))}$ for all triangles $i \stackrel{a}{\longrightarrow} j \stackrel{b}{\longrightarrow} k \stackrel{c}{\longrightarrow} i$ on $\mathcal{G}$.
\end{Def}

\begin{Prop} \label{Prop:lam/mu=C}
Let $W$ be a $\nu$-invariant cell system on $\mathcal{G}$.
There is a constant $C$ such that
\begin{equation} \label{eqn:def-C}
\frac{\lambda_{a}^{(m)}\lambda_{b}^{(m')}\lambda_{c}^{(m'')}}{\mu_{a}^{(m)}\mu_{b}^{(m')}\mu_{c}^{(m'')}} = C
\end{equation}
for all triangles $i \stackrel{a}{\longrightarrow} j \stackrel{b}{\longrightarrow} k \stackrel{c}{\longrightarrow} i$ on $\mathcal{G}$, and all $m \in \{ 1,p_{ij} \}$, $m' \in \{ 1,p_{jk} \}$, $m'' \in \{ 1,p_{ki} \}$.
\end{Prop}

\noindent \emph{Proof}:
The dual $A^{\ast}$ of $A$ is an $A$-$A$ bimodule with the products $\varphi x$, $x \varphi$ defined by $(\varphi x)(y) = \varphi(xy)$, $(x \varphi)(y) = \varphi(yx)$, for $\varphi \in A^{\ast}$ and $x,y \in A$.
The element dual to $v_{j\nu(i)}^m$ is $(v_{\nu(i)j}^m)^{\ast} = \lambda_{a}^{(m)} u^{\ast}_{\nu(i)i} a$, where $s(a) = i$, $r(a) = j$, and $u^{\ast}_{\nu(i)i}$ is the element dual to $u_{i\nu(i)}$, since
$$(v_{\nu(i)j}^m)^{\ast}(v_{j\nu(i)}^m) = \lambda_{a}^{(m)} (u^{\ast}_{\nu(i)i} a)(v_{j\nu(i)}^m) = \lambda_{a}^{(m)} u^{\ast}_{\nu(i)i} (a v_{j\nu(i)}^m) = u^{\ast}_{\nu(i)i}(u_{i\nu(i)}) = 1.$$
Similarly, the element $\mu_{a'}^{(m)} \nu(a') u^{\ast}_{\nu(j)j}$, where $s(a') = i$, $r(a') = j$, is also dual to $v_{j\nu(i)}^m$. Then (\ref{eqn:I}) dualises to give
\begin{equation} \label{eqn:II}
\lambda_{a}^{(m)} u^{\ast}_{\nu(i)i} a = (v_{\nu(i)j}^m)^{\ast} = \mu_{a'}^{(m)} \nu(a') u^{\ast}_{\nu(j)j}.
\end{equation}
Let $a = a'$ and $m \in \{ 1,p_{ij} \}$. Then multiplying on the right by $(\lambda_{a}^{(m)})^{-1} W_{ijk}^{(a,b,c)} b$ in (\ref{eqn:II}) for a vertex $k$ such that there is a triangle $i \stackrel{a}{\longrightarrow} j \stackrel{b}{\longrightarrow} k \stackrel{}{\longrightarrow} i$ on $\mathcal{G}$, we have
$$W_{ijk}^{(a,b,c)} u^{\ast}_{\nu(i)i} ab = \frac{\mu_{a}^{(m)}}{\lambda_{a}^{(m)}} W_{ijk}^{(a,b,c)} \nu(a) u^{\ast}_{\nu(j)j} b = \frac{\mu_{a}^{(m)} \mu_{b}^{(m')}}{\lambda_{a}^{(m)} \lambda_{b}^{(m')}} W_{ijk}^{(a,b,c)} \nu(a) \nu(b) u^{\ast}_{\nu(k)k},$$
where the second equality follows from (\ref{eqn:II}) for a choice of $m' \in \{ 1,p_{jk} \}$. Summing over all vertices $j$ and edges $a$, $b$ such that there is a triangle $i \stackrel{a}{\longrightarrow} j \stackrel{b}{\longrightarrow} k \stackrel{c}{\longrightarrow} i$ on $\mathcal{G}$, and making a choice of $m = m_j$, $m' = m_j'$ for each $j$, the L.H.S. is zero by the relation $\rho_c$ in (\ref{eqn:relations-A}), and so we obtain $\sigma_{(\textbf{m})} u^{\ast}_{\nu(k)k} = 0$ where
$$\sigma_{(\textbf{m})} = \sum_{j,a,b} \frac{\mu_{a}^{(m_j)} \mu_{b}^{(m_j')}}{\lambda_{a}^{(m_j)} \lambda_{b}^{(m_j')}} W_{ijk}^{(a,b,c)} \nu(a) \nu(b),$$
and $\textbf{m} = (m_1,m_1',m_2,m_2',\ldots,m_r,m_r')$ where $r$ is the number of vertices $j$ in the summation.
Suppose $\sigma_{(\textbf{m})} \neq 0$. Then there exists a non-zero $v \in k \cdot A_{h-5} \cdot \nu(i)$ such that $v \sigma_{(\textbf{m})} = u_{k\nu(k)} \in k \cdot A_{h-3} \cdot \nu(k)$, and $\sigma_{(\textbf{m})} u^{\ast}_{\nu(k)k} = V^{\ast} \neq 0$ (since $v \neq 0$) which is a contradiction.
Thus $\sigma_{(\textbf{m})} = 0$, which implies
\begin{equation} \label{mu/lamW=W}
\frac{\mu_{a}^{(m)} \mu_{b}^{(m')}}{\lambda_{a}^{(m)} \lambda_{b}^{(m')}} W_{ijk}^{(a,b,c)} = \xi_{ik} W_{\nu(i)\nu(j)\nu(k)}^{(\nu(a), \nu(b), \nu(c))},
\end{equation}
where $\xi_{ik} \in \mathbb{C}$ does not depend on $j$, $a$ or $b$. For $j=j_1$, the left hand side of (\ref{mu/lamW=W}) only depends on $m_1$ and $m_1'$, whilst for $j=j_2$, the left hand side of (\ref{mu/lamW=W}) now only depends on $m_2$ and $m_2'$. Thus $\xi_{ik}$ does not depend on $\textbf{m}$, but only on $i$, $k$.
Define $\xi_{ik}' := \frac{\mu_{c}^{(m'')}}{\lambda_{c}^{(m'')}} \xi_{ik}$, which only depends on $i$, $k$ and $m''$.
Then since $W$ is $\nu$-invariant, we have
\begin{equation} \label{mu/lam=xi}
\frac{\mu_{a}^{(m)}\mu_{b}^{(m')}\mu_{c}^{(m'')}}{\lambda_{a}^{(m)}\lambda_{b}^{(m')}\lambda_{c}^{(m'')}} = \xi_{ik}',
\end{equation}
which only depends on $i$, $k$ and $m''$. However, by a similar argument, the left hand side is also equal to $\xi_{kj}'$ (which only depends on $k$, $j$ and $m'$) and $\xi_{ji}'$ (which only depends on $j$, $i$ and $m$). Hence the left hand side does not depend on $i$, $j$, $k$ or the choices of $m$, $m'$ and $m''$. So we have $\xi_{ik}' = \xi_{kj}' = \xi_{ji}' =: 1/C$ for all triangles $i \stackrel{a}{\longrightarrow} j \stackrel{b}{\longrightarrow} k \stackrel{c}{\longrightarrow} i$ on $\mathcal{G}$.
\hfill
$\Box$

Since $\xi_{ik}'$ in (\ref{mu/lam=xi}) does not depend on the choice of $m \in \{ 1,p_{ij} \}$, we obtain as an immediate corollary:

\begin{Cor} \label{Cor:lam/mu(m)}
For any $\nu$-invariant cell system $W$ on $\mathcal{G}$,
$$\frac{\lambda_{a}^{(m)}}{\mu_{a}^{(m)}} = \frac{\lambda_{a}^{(m')}}{\mu_{a}^{(m')}}$$
for any edge $i \stackrel{a}{\longrightarrow} j$ on $\mathcal{G}$, and $m,m' \in \{ 1,p_{ij} \}$.
\end{Cor}

We will now define an alternative basis $\{ v_{j\nu(i)}^a \}$ for $i \cdot A_{h-4} \cdot \nu(k)$ such that for any edges $i \stackrel{a,a'}{\longrightarrow} j$ and each edge $k \stackrel{b,b'}{\longrightarrow} i$ on $\mathcal{G}$ we have $\lambda_{a'} a' v_{j\nu(i)}^{a} = \delta_{a,a'} u_{i\nu(i)}$, $\mu_{b'} v_{i\nu(k)}^{b} \nu(b') = \delta_{b,b'} u_{i\nu(i)}$, where $\lambda_a$, $\mu_a$ are non-zero scalars. This basis will be used in Section \ref{sect:Omega^4(A)}.

For any double edge $(a,a')$ on $\mathcal{G}$ with $s(a)=i$, $r(a)=j$, if we set $v_{j\nu(i)}^a := \lambda_{a'}^{(1)} v_{j\nu(i)}^1 - \lambda_{a'}^{(2)} v_{j\nu(i)}^2$ and $v_{j\nu(i)}^{a'} := \lambda_{a}^{(1)} v_{j\nu(i)}^1 - \lambda_{a}^{(2)} v_{j\nu(i)}^2$, then we have $a'v_{j\nu(i)}^a = 0 = av_{j\nu(i)}^{a'}$ and $\lambda_b bv_{j\nu(i)}^b = u_{i\nu(i)}$, $b \in \{ a,a' \}$, where $\lambda_b = \varepsilon_b\lambda_b^{(1)}\lambda_b^{(2)}/(\lambda_{a'}^{(1)}\lambda_{a}^{(2)}-\lambda_{a}^{(1)}\lambda_{a'}^{(2)})$ and $\varepsilon_a = 1$, $\varepsilon_{a'} = -1$.
Then from (\ref{eqn:I}) and Corollary \ref{Cor:lam/mu(m)}, we see that $v_{j\nu(i)}^a \nu(a') = 0 = v_{j\nu(i)}^{a'} \nu(a)$ and $\mu_b v_{j\nu(i)}^b \nu(b) = u_{i\nu(i)}$, $b \in \{ a,a' \}$, where $\mu_a = \mu_a^{(1)}\mu_a^{(2)}/(\lambda_{a'}^{(1)}\mu_{a}^{(2)} - \mu_{a}^{(1)}\lambda_{a'}^{(2)})$ and $\mu_{a'} = \mu_{a'}^{(1)}\mu_{a'}^{(2)}/(\lambda_{a}^{(1)}\mu_{a'}^{(2)} - \mu_{a'}^{(1)}\lambda_{a}^{(2)})$.

Let $v_{j\nu(i)}^a = v_{j\nu(i)}^1$, $\lambda_a = \lambda_a^{(1)}$ when there is only one edge $a$ from $i$ to $j$, and if there is a double edge $(a,a')$ from $i$ to $j$ let $v_{j\nu(i)}^b$, $\lambda_b$ be as defined in the previous paragraph, where $b \in \{ a,a' \}$.
Thus we have defined an alternative basis $\{ v_{j\nu(i)}^a \}$ for $i \cdot A_{h-4} \cdot \nu(k)$ such that there for each edge $i \stackrel{a}{\longrightarrow} j$ and each edge $k \stackrel{b}{\longrightarrow} i$ on $\mathcal{G}$ we have
\begin{equation} \label{eqn:Ia}
\lambda_{a} a v_{j\nu(i)}^{a} = u_{i\nu(i)} = \mu_{b} v_{i\nu(k)}^{b} \nu(b),
\end{equation}
where $\lambda_a$, $\mu_a$ are non-zero scalars.
Since $\mathrm{dim}(i \cdot A_{h-3} \cdot l) = 0$ for all $l \neq \nu(i)$, we have
\begin{equation} \label{eqn:Ib}
b v_{j\nu(i)}^{a} = 0 = v_{j\nu(i)}^{a} \nu(b)
\end{equation}
when $a \neq b$.
We will usually write $v_{j\nu(i)}$ for $v_{j\nu(i)}^a$ where there is only one edge $a$ from $i$ to $j$. Dualising (\ref{eqn:Ia}) we thus get:
\begin{equation} \label{eqn:IIa}
\lambda_{a} u^{\ast}_{\nu(i)i} a = (v_{\nu(i)j}^a)^{\ast} = \mu_{a'} \nu(a') u^{\ast}_{\nu(j)j}.
\end{equation}

\subsection{Computation of the constant $C$} \label{sect:PropC=1}

In this section we compute the value of the constant $C$ in Proposition \ref{Prop:lam/mu=C}.

\begin{Prop} \label{Prop:C=1}
Let $\mathcal{G}$ be an $SU(3)$ $\mathcal{ADE}$ graph, which is not $\mathcal{E}_4^{(12)}$. Then for the permutation $\nu$ of $\mathcal{G}$ defined
in Section \ref{sect:nakayama},
the constant $C$ in Proposition \ref{Prop:lam/mu=C} is 1.
\end{Prop}

The proof of Proposition \ref{Prop:C=1} is done in a case-by-case method, where we will use the Ocneanu cells $W(\triangle)$ computed in \cite{evans/pugh:2009i}.
For the graphs $\mathcal{D}^{(n)}$, $\mathcal{D}^{(n)\ast}$ and $\mathcal{E}_1^{(12)}$, we did not claim to have computed all cell systems up to equivalence in \cite{evans/pugh:2009i}.

We will begin by describing the general strategy.
We choose a vertex $i$ of $\mathcal{G}$ which is the source of only one edge $a$, and similarly the range of only one edge $b$. We denote by $j$, $k$ the range, source vertices of the edges $a$, $b$ respectively. Note that there must necessarily be (at least one) edge $c$ from $j$ to $k$ on $\mathcal{G}$.
We choose a non-zero path $u_{i\nu(i)} = [i \; j \; l_1 \; l_2 \; \cdots \; l_{h-6} \; \nu(k) \; \nu(i)] \in i \cdot A_{h-3} \cdot \nu(i)$ of length $h-3$ from $i$ to $\nu(i)$, where $l_1, \ldots, l_{h-6}$ are vertices of $\mathcal{G}$. We let the elements $v_{j\nu(i)}$ and $v_{i\nu(k)}$ be the paths $v_{j\nu(i)} = [j \; l_1 \; l_2 \; \cdots \; l_{h-6} \; \nu(k) \; \nu(i)] \in j\cdot A_{h-4} \cdot \nu(i)$ and $v_{i\nu(k)} = [i \; j \; l_1 \; l_2 \; \cdots \; l_{h-6} \; \nu(k)] \in i \cdot A_{h-4} \cdot \nu(k)$. Then $u_{i\nu(i)} = a v_{j\nu(i)} = v_{i\nu(k)} \nu(b)$ in $A$, and we have $\lambda_a = \mu_b = 1$. Note that since $j\cdot A_{h-4} \cdot \nu(i)$ and $i \cdot A_{h-4} \cdot \nu(k)$ are now one-dimensional, we omit the notation $m=1$ from $\lambda_{a'}^{(m)}$, $\mu_{a'}^{(m)}$.
We now form the paths $v_{j\nu(i)} \nu(a) \in j \cdot A_{h-3} \cdot \nu(j)$, $b v_{i\nu(k)} \in k \cdot A_{h-3} \cdot \nu(k)$, and transform these using (\ref{eqn:I}) so that we have
\begin{eqnarray}
v_{j\nu(i)} \nu(a) & = & \frac{\lambda_c}{\mu_a} c v_{k\nu(j)}, \label{eqn:vy1} \\
b v_{i\nu(k)} & = & \frac{\mu_c}{\lambda_b} v_{k\nu(j)} \nu(c), \label{eqn:yv1}
\end{eqnarray}
where the same path $v_{k\nu(j)}$ appears in the right hand side of both equalities.
By Proposition \ref{Prop:lam/mu=C}, $\lambda_c/\mu_a = C \mu_c \mu_b/\lambda_a \lambda_b = C \mu_c/\lambda_b$, since $\lambda_a = \mu_b = 1$, so that (\ref{eqn:vy1}) becomes
\begin{equation} \label{eqn:vy2}
v_{j\nu(i)}  \nu(a) = \frac{\mu_c}{\lambda_b} C c v_{k\nu(j)}.
\end{equation}
Then we obtain the value of $\mu_c/\lambda_b$ from (\ref{eqn:yv1}), and substituting into (\ref{eqn:vy2}) we can determine the value of $C$.

For the graphs where $\nu$ is the identity, that the constant $C$ which appears in (\ref{eqn:vy2}) is 1 follows from the following considerations.
There is a conjugation on the $\mathcal{A}$ graphs given by the conjugation on the representations of $SU(3)$ given in Section \ref{sect:verlinde-sector}.
There is also a conjugation on the other $SU(3)$ $\mathcal{ADE}$ graphs: The conjugation $\tau: {}_N \mathcal{X}_N \rightarrow {}_N \mathcal{X}_N$ on the braided system of endomorphisms of $SU(3)_k$ on a factor $N$, given by the conjugation on the representations of $SU(3)$, induces a conjugation $\tau: {}_N \mathcal{X}_M \rightarrow {}_N \mathcal{X}_M$ such that $G_{\overline{\lambda}} = \tau G_{\lambda} \tau$, where $G_{\lambda} a = \lambda a$ for $\lambda \in {}_N \mathcal{X}_N$, $a \in {}_N \mathcal{X}_M$.
We call a path symmetric if it is invariant under reversal of the path and taking the conjugate of the graph, i.e. a path $x = [a_1 \; a_2 \; \cdots \; a_{l-1} \; a_l]$ is symmetric if $x^c := [\overline{a_l} \; \overline{a_{l-1}} \; \cdots \; \overline{a_2} \; \overline{a_1}] = x$, where $\overline{v}$ denotes the image of the vertex $v$ under conjugation of the graph. We note that for a path to be symmetric, it must start at a vertex of colour $p$ and end at a vertex of colour $-p \textrm{ mod } 3$.
We choose a non-zero path $u_{i\nu(i)}$ which is symmetric. We then form the path $v_{j\nu(i)} \nu(a)$ as above, and transform this to a scalar multiple $d_1$ of $c v_{k\nu(j)}$, where $v_{k\nu(j)}$ is a basis path in $k \cdot A_{h-4} \cdot \nu(j)$. We also find a non-zero path $c' v' \in j \cdot A_{h-3} \cdot \nu(j)$, where $v' \in k \cdot A_{h-4} \cdot \nu(j)$ is a symmetric path of length $h-4$, and transform this to a scalar multiple $d_2$ of $c v_{k\nu(j)}$. Then we obtain that $v_{j\nu(i)} \nu(a) = d_1 d_2^{-1} c' v'$. Now, since $u_{i\nu(i)}$ is symmetric, $(v_{j\nu(i)} \nu(a))^c = b v_{i\nu(k)}$. Then $b v_{i\nu(k)} = (v_{j\nu(i)} \nu(a))^c = (d_1 d_2^{-1} c' v')^c = d_1 d_2^{-1} v' \nu(c') = d_1 v_{k\nu(j)} \nu(c)$. So we see that the same coefficient $d_1$ appears in both (\ref{eqn:yv1}) and (\ref{eqn:vy2}), so that $C=1$. However, this method will not work for the graph $\mathcal{E}_5^{(12)}$ even though it has $\nu = \mathrm{id}$, since for all the possible choices for the vertices $j$, $k$, one of these must be of colour 0, hence the conjugation of the graph does not interchange $j \leftrightarrow k$.

For the graphs $\mathcal{E}_1^{(12)}$, $\mathcal{E}_2^{(12)}$, $\mathcal{E}_5^{(12)}$ and $\mathcal{E}^{(24)}$ it is not clear that chosen paths are non-zero in $l \cdot A_{h-3} \cdot \nu(l)$, for $l=i,j,k$, and we will need to construct a basis for the space of all paths which start at the vertices $i$, $j$ (and also $k$ in the case of $\mathcal{E}_5^{(12)}$).
The computations of these basis paths are lengthy and are contained in the Appendix to \cite{evans/pugh:2010ii-arxiv}. We will summarize the results for the computation of $C$ for these graphs here, with the detailed computations given in the Appendix to \cite{evans/pugh:2010ii-arxiv}.
In what follows we will use the notation $a_{rs}$ to denote the edge from $r$ to $s$ where this edge is unique. \\

\noindent \emph{The identity $\mathcal{A}$ graphs}:

We will first compute the value of the constant $C$ for the graphs $\mathcal{A}^{(n)}$.
The infinite graph $\mathcal{A}^{(\infty)}$, illustrated in \cite[Figure 4]{evans/pugh:2009i}, is the fusion graph for the irreducible representations of $SU(3)$, whilst for finite $n$, the graph $\mathcal{A}^{(n)}$ is the fusion graph for the irreducible representations of $SU(3)_n$, as described in Section \ref{sect:verlinde-sector}.
We will write $(p,q)$ for the irreducible representation $\lambda_{(p,q)}$.
For $\mathcal{A}^{(n)}$ the automorphism $\nu$ is the clockwise rotation of the graph by $2 \pi /3$. Choosing the vertices $i=(0,0)$, $j=(1,0)$ and $k=(0,1)$, we have $\nu(i)=(n-3,0)$, $\nu(j)=(n-4,1)$ and $\nu(k)=(n-4,0)$.
The unique cell system $W$ (up to equivalence) was computed in \cite[Theorem 5.1]{evans/pugh:2009i}, and we use the same notation for the cells here. The cell system $W$ is $\nu$-invariant.

\begin{figure}[tb]
\begin{center}
  \includegraphics[width=155mm]{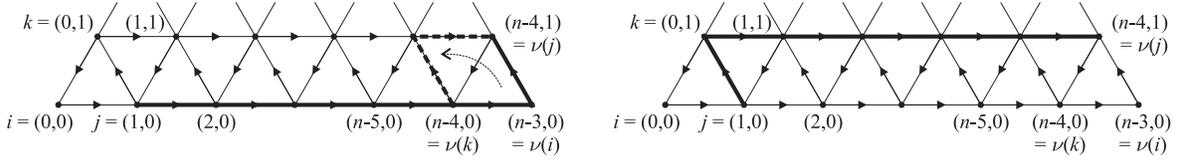}\\
  \caption{Two paths in $j \cdot A_{h-3} \cdot \nu(j)$ for $\mathcal{A}^{(n)}$: $v_{j\nu(i)} \nu(a_{ij})$ on the left and $a_{ij} v_{k\nu(j)}$ on the right.} \label{fig:pathsA-jv(j)}
\end{center}
\end{figure}

\begin{figure}[tb]
\begin{center}
  \includegraphics[width=155mm]{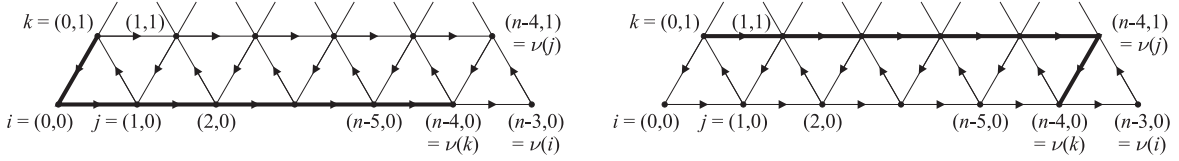}\\
  \caption{Two paths in $k \cdot A_{h-3} \cdot \nu(k)$ for $\mathcal{A}^{(n)}$: $a_{ki} v_{i\nu(k)}$ on the left and $v_{k\nu(j)} \nu(a_{ij})$ on the right.} \label{fig:pathsA-kv(k)}
\end{center}
\end{figure}

There is only one possible non-zero path of length $h-3$ from $i=(0,0)$ to $\nu(i)=(n-3,0)$, which is given by $u_{i\nu(i)} = [(0,0) \; (1,0) \; (2,0) \; \cdots \; (n-4,0) \; (n-5,0)]$. Any path of length $> h-3$ will be zero in $A$ since using the relations on $A$ it can be transformed to a path which begins $[i \; j \; k \; \cdots \; ] = 0$ since $[i \; j \; k] = \rho_{a_{ki}} = 0$ in $A$, where $a_{rs}$ is the edge from $r$ to $s$.
Then $v_{j\nu(i)} = [(1,0) \; (2,0) \; \cdots \; (n-4,0) \; (n-5,0)]$ and $v_{i\nu(k)} = [(0,0) \; (1,0) \; (2,0) \; \cdots \; (n-4,0)]$.
We form $v_{j\nu(i)} \nu(a_{ij})$ and using the relation $\rho_{\nu(a_{jk})}$ we obtain the path $[(1,0) \; (2,0) \; \cdots \; (n-4,0) \; (n-5,1) \; (n-4,1)]$, as shown on the left hand side of Figure \ref{fig:pathsA-jv(j)}, with coefficient $-W_{\nabla(n-5,0)}/W_{\triangle(n-4,0)}$, where we use the notation $W_{\triangle(k,m)} = W_{(k,m)(k+1,m)(k,m+1)}$ and $W_{\nabla(k,m)} = W_{(k+1,m)(k,m+1)(k+1,m+1)}$. Continuing in this way we obtain the path $a_{jk} v_{k\nu(j)} = [(1,0) \; (0,1) \; (1,1) \; (2,1) \; \cdots \; (n-4,1)]$, as shown on the right hand side of Figure \ref{fig:pathsA-jv(j)}, with coefficient
$$\xi = (-1)^{n-4} \frac{W_{\nabla(0,0)} W_{\nabla(1,0)} \ldots W_{\nabla(n-5,0)}}{W_{\triangle(1,0)} W_{\triangle(2,0)} \ldots W_{\triangle(n-4,0)}}.$$ Similarly, we form $a_{ki} v_{i\nu(k)}$ and transform using the relations to obtain
\begin{eqnarray*}
a_{ki} v_{i\nu(k)} & = & (-1)^{n-4} \frac{W_{\nabla(0,0)} W_{\nabla(1,0)} \ldots W_{\nabla(n-5,0)}}{W_{\triangle(0,0)} W_{\triangle(1,0)} \ldots W_{\triangle(n-5,0)}} a_{jk} v_{k\nu(j)} \;\; = \;\; \xi \; \frac{W_{\triangle(n-4,0)}}{W_{\triangle(0,0)}} a_{jk} v_{k\nu(j)} \\
& = & \xi \; a_{jk} v_{k\nu(j)},
\end{eqnarray*}
where the last equality follows from the $\nu$-invariance of the cell system $W$. Thus we see that $C=1$ for the graphs $\mathcal{A}^{(n)}$.
The only properties of the cells $W(\triangle)$ that are used here are their $\nu$-invariance and the fact that they are non-zero. \\

\begin{figure}[tb]
\begin{center}
  \includegraphics[width=55mm]{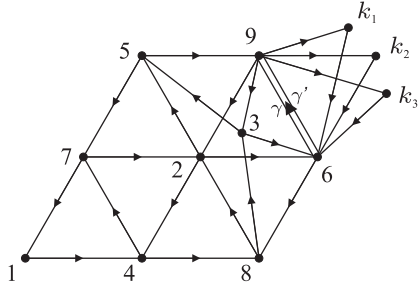}\\
  \caption{Graph $\mathcal{D}^{(9)}$} \label{fig:D(9)}
\end{center}
\end{figure}

\noindent \emph{The orbifold $\mathcal{D}$ graphs}:

We will now consider the graphs $\mathcal{D}^{(n)}$, which are $\mathbb{Z}_3$ orbifolds of the graphs $\mathcal{A}^{(n)}$. The graph $\mathcal{D}^{(9)}$ is illustrated in Figure \ref{fig:D(9)}. The weights $W(\triangle)$ for $\mathcal{A}^{(n)}$ are invariant under the $\mathbb{Z}_3$ symmetry of the graph given by rotation by $2\pi/3$. Thus there is an orbifold solution for the cell system $W$ on $\mathcal{D}^{(n)}$ where the weights $W(\triangle)$ are given by the corresponding weights for $\mathcal{A}^{(n)}$ \cite[Theorems 6.1 \& 6.2]{evans/pugh:2009i}. More precisely, excluding triangles $\triangle$ which contain one of the triplicated vertices $(k,k)_l$ in the case where $n = 3k+3$, the weight $W(\triangle_{i_1,i_2,i_3})$ for the triangle $\triangle_{i_1,i_2,i_3} = i_1 \rightarrow i_2 \rightarrow i_3 \rightarrow i_1$ on $\mathcal{D}^{(n)}$ is given by the weight $W(\triangle_{i_1^{(0)},i_2^{(1)},i_3^{(2)}}) = W(\triangle_{i_1^{(1)},i_2^{(2)},i_3^{(0)}}) = W(\triangle_{i_1^{(2)},i_2^{(0)},i_3^{(1)}})$ for $\mathcal{A}^{(n)}$, where $i_k^{(0)}$, $i_k^{(1)}$, $i_k^{(2)}$ are the three vertices of $\mathcal{A}^{(n)}$ which are identified under the $\mathbb{Z}_3$ action to give the vertex $i_k$ of $\mathcal{D}^{(n)}$, $k=1,2,3$.
If for a triangle $\triangle_{i_1,i_2,i_3}$ on $\mathcal{D}^{(n)}$ there is no choice of vertices $i_1^{(j_1)}$, $i_2^{(j_2)}$, $i_3^{(j_3)}$ on $\mathcal{A}^{(n)}$ which lie on a closed loop of length three $i_1^{(j_1)} \rightarrow i_2^{(j_2)} \rightarrow i_3^{(j_3)} \rightarrow i_1^{(j_1)}$, then we have $W(\triangle_{i_1,i_2,i_3}) = 0$.
When $n = 3k+3$, the weight $W(\triangle)$ for a triangle $\triangle$ which contain one of the triplicated vertices $(k,k)_l$ is just given by one third of the weight for the corresponding triangle on $\mathcal{A}^{(3k+3)}$.
Thus the relations (\ref{eqn:relations-A}) for $\mathcal{D}^{(n)}$ are given precisely by the relations for $\mathcal{A}^{(n)}$, except for the relations $\rho_{\gamma}$, $\rho_{\gamma'}$ in the case where $n = 3k+3$, which involve the triplicated vertices $(k,k)_l$. However, these last two relations are not used to show $C=1$, and thus the result for $\mathcal{D}^{(n)}$ follows from the result for $\mathcal{A}^{(n)}$ under the orbifold procedure.
The vertex $(n-3,0)$ of $\mathcal{A}^{(n)}$ is identified with the distinguished vertex $(0,0)$ of $\mathcal{A}^{(n)}$. Thus with $i$ the distinguished vertex of $\mathcal{D}^{(n)}$ with lowest Perron-Frobenius weight, the element $u_{i\nu(i)}$ is a closed loop of length $n-3$ starting and ending at $i$. We see that the permutation $\nu$ must be the identity for $\mathcal{D}^{(n)}$.

We will illustrate the general $\mathcal{D}$ case by giving the computations for the graph $\mathcal{D}^{(9)}$. We label the vertices as in Figure \ref{fig:D(9)}, and denote the two edges in the double edge by $\gamma$ and $\gamma'$. We write $W_{s(a)s(b)s(c)}$ for $W(\triangle^{(abc)}_{s(a)s(b)s(c)}$ where $a,b,c \not \in \{ \gamma, \gamma' \}$, and $W^{\eta}_{s(a)s(b)s(c)}$ when one of $a$, $b$ or $c$ is the edge $\eta \in \{ \gamma, \gamma' \}$.
We note that $W_{269}^{\gamma'} = W_{369}^{\gamma} = 0$ \cite[Theorem 6.2]{evans/pugh:2009i}. We choose the vertices $i=1$, $j=4$ and $k=7$. Then the path $u_{i\nu(i)}$ is obtained from the corresponding element for $\mathcal{A}^{(9)}$, so $u_{i\nu(i)} = u_{11} = [1483571]$ which is symmetric. We have $v_{j\nu(i)} = v_{41} = [483571]$, $v_{i\nu(k)} = v_{17} = [148357]$. Since $\nu = \mathrm{id}$, we only need to transform $v_{j\nu(i)} \nu(a) = [4835714]$, where $a$ is the edge from $i$ to $j$, to a scalar multiple of a path $\gamma' v'$ where $v'$ is some symmetric path. We will underline at each stage the subpath of length 2 which we transform using the relations. We have
\begin{eqnarray*}
[4835\underline{714}] & = & \frac{W_{247}}{W_{147}} [483\underline{572}4] \;\; = \;\; \frac{W_{247}W_{259}}{W_{147}W_{257}} [48\underline{359}24] \;\; = \;\; \frac{W_{247}W_{259}W_{369}^{\gamma'}}{W_{147}W_{257}W_{359}} [4\underline{836}'924] \\
& = & \frac{W_{247}W_{259}W_{369}^{\gamma'}W_{268}}{W_{147}W_{257}W_{359}W_{368}} [\underline{482}6'924] \;\; = \;\; \frac{W_{247}W_{259}W_{369}^{\gamma'}W_{268}W_{247}}{W_{147}W_{257}W_{359}W_{368}W_{248}} [4726'924],
\end{eqnarray*}
where $[69]$, $[6'9]$ denote the edge from 6 to 9 along $\gamma$, $\gamma'$ respectively.
The path $[4726'924]$ is non-zero, with $[726'924]$ symmetric.
Here it is important that the cells $W(\triangle)$ for some of the triangles which contain one of the double edges are zero, otherwise the path $v_{j\nu(i)} = v_{41} = [483571]$ might be zero in $A$.
The only other property of the cells $W(\triangle)$ that is used here is the fact that those that appear in the coefficient of $[4726'924]$ in the above equality are non-zero.
If there was another inequivalent cell system on $\mathcal{D}$ such that this coefficient is non-zero and the Hilbert series of $A$ was given by (\ref{eqn:Hilbert_Series-SU(3)ADE}), then the result $C=1$ would hold for this cell system also. \\

\noindent \emph{The conjugate $\mathcal{A}^{\ast}$ and conjugate orbifold $\mathcal{D}^{\ast}$ graphs}:

The proof for these graphs is slightly different to that for all the other graphs in that we do not choose $i$ to be a vertex which is the source of only one edge. Here the vertex $i$ is chosen to be the source of exactly two edges.
First consider the graphs $\mathcal{A}^{(n)\ast}$, where $\nu = \mathrm{id}$.
The unique cell system $W$ (up to equivalence) was computed in \cite[Theorems 7.1, 7.3 \& 7.4]{evans/pugh:2009i}, and we use the same notation for the cells here.
The $\mathcal{A}^{(n)\ast}$ graphs are illustrated in \cite[Figure 11]{evans/pugh:2009i}. We illustrate the cases $n=7,8$ in Figure \ref{fig:A(7,8)star}. The labelling we use here for the vertices of $\mathcal{A}^{(2m+1)\ast}$ is the reverse of the labelling used in \cite{evans/pugh:2009i}.
The relations in $A(\mathcal{A}^{(n)\ast},W)$ are
\begin{eqnarray}
& W_{112}[121] + W_{111}[111] = 0, & \nonumber \label{eqn:rel-Astar-i} \\
& W_{a-1,a,a}[a(a-1)a] + W_{a,a,a}[aaa] + W_{a,a,a+1}[a(a+1)a] = 0, & \label{eqn:rel-Astar-ii}  \\
& W_{a,a,a+1}[aa(a+1)] + W_{a,a+1,a+1}[a(a+1)(a+1)] = 0, & \label{eqn:rel-Astar-iii} \\
& W_{a,a,a+1}[(a+1)aa] + W_{a,a+1,a+1}[(a+1)(a+1)a] = 0, & \label{eqn:rel-Astar-iv}
\end{eqnarray}
where $a=2,\ldots,p-1$ in (\ref{eqn:rel-Astar-ii}), and $a=1,\ldots,a'$ in (\ref{eqn:rel-Astar-iii}), (\ref{eqn:rel-Astar-iv}), where $p = \lfloor (n-1)/2 \rfloor$, $a'=p-1$ for even $n$, and $a'=p-2$ for odd $n$.
For even $n$ we have the extra relation $W_{p-1,p,p}[p(p-1)p] + W_{p,p,p}[ppp] = 0$, and for odd $n$ we have the extra relation $[p(p-1)(p-1)] = [(p-1)(p-1)p] = 0$.

\begin{figure}[tb]
\begin{minipage}[t]{7.5cm}
\begin{center}
  \includegraphics[width=35mm]{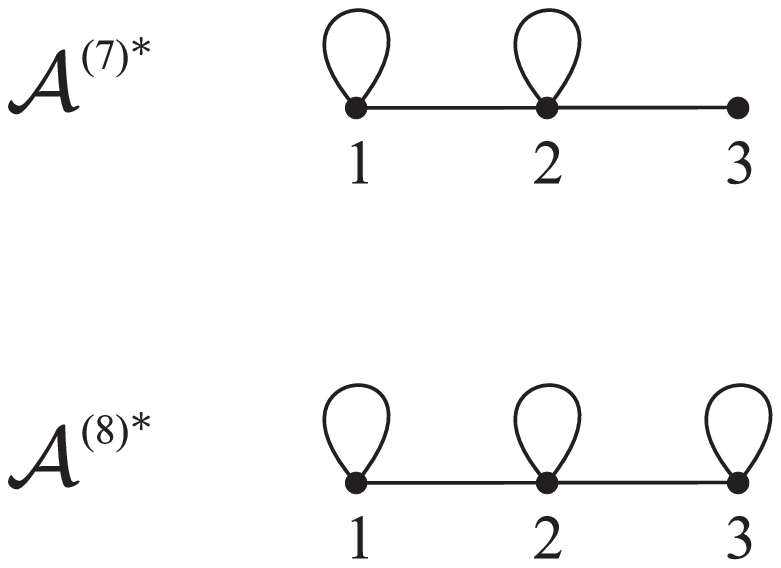}\\
  \caption{Graphs $\mathcal{A}^{(7)\ast}$, $\mathcal{A}^{(8)\ast}$} \label{fig:A(7,8)star}
\end{center}
\end{minipage}
\hfill
\begin{minipage}[t]{7.5cm}
\begin{center}
  \includegraphics[width=40mm]{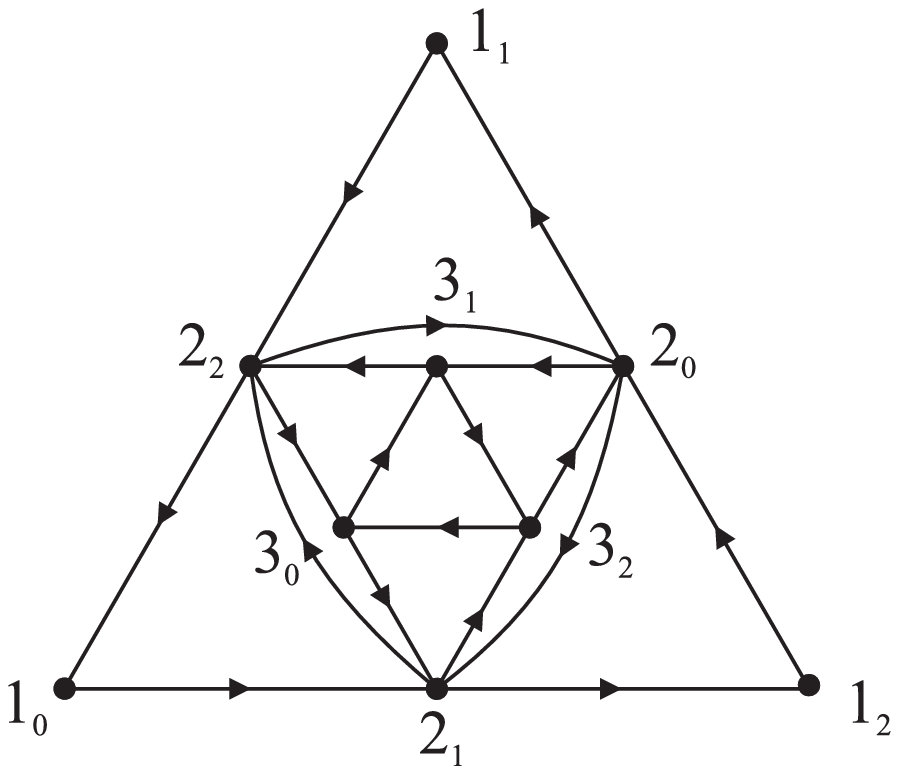}\\
  \caption{Graph $\mathcal{D}^{(7)\ast}$} \label{fig:D(7)star}
\end{center}
\end{minipage}
\end{figure}

We first consider the even case $n=2m+2$.
We choose the vertices $i=j=k=1$ as illustrated in Figure \ref{fig:A(7,8)star}. The element $u_{i\nu(i)} = u_{1\nu(1)}$ of length $2m+1$ is (up to some scalar) $[111\cdots1]$. We show this is non-zero by induction. When $m=1$, we see using the relations that $[1\underline{111}] = -(W_{112}/W_{111}) [\underline{112}1] = (W_{112}^2/W_{111}^2) [1\underline{221}] = -(W_{112}/W_{111}) [1211]$, so that all paths of length $3$ are equal, up to scalar multiple, in $A$.
We assume that all non-zero paths in $A(\mathcal{A}^{(2m+2)\ast},W)$ are proportional to the path $u_{1\nu(1)}^{(k)} = [111\cdots1]$ of length $2m-1$ for $m=k$.
For $m=k+1$, any path in $A(\mathcal{A}^{(2k+4)\ast},W)$ of length $2k+1$ must have one of the following forms, where $a\{j\} = [v_0 v_1 \cdots v_{2k-1}]$ is a path of length $2k-1$ with $v_{2k-1}=j$: \hspace{2mm} (I) $[a\{1\} 11] = -W_{112}/W_{111} [a\{1\} 21]$, \hspace{2mm} (II) $[a\{2\} 11] = -W_{122}/W_{112} [a\{2\} 21]$, \hspace{2mm} (III) $[a\{3\} 21]$.
Any path of form (I) is clearly proportional to $u_{1\nu(1)}^{(k+1)} = [111\cdots1]$ of length $2k-1$, since any path $a\{1\}$ in $A(\mathcal{A}^{(2k+4)\ast},W)$ of length $2k-3$ must be proportional to $u_{1\nu(1)}^{(k)}$ (we note that a non-zero coefficient when $m=k$ may become zero when $m=k+1$ since we replace the quantum integer $[s]_{q(k)}$ by $[s]_{q(k+1)}$ in the weights $W$, where $q(m) = e^{2 \pi i/m}$).
We now consider paths of form (II). There are three cases: (a) Suppose $v_{2k-4} = 1$ in $a\{2\}$. Then using the relations we see that $[a\{2\} 11] = -W_{111}/W_{112} [a\{1\} 11]$, where $a\{1\}$ is obtained from $a\{2\}$ by replacing its last edge by a closed loop from 1 to 1. This path is now of form (I). (b) If $v_{2k-4} = 2$, then we have $[a\{2\} 11] = -W_{112}/W_{122} [a\{1\} 11]$, where now $a\{1\}$ is obtained from $a\{2\}$ by replacing its last edge by an edge from 2 to 1, and this new path is of form (I). (c) If $v_{2k-4} = 3$, then we move to consider $v_{2k-5}$. Then we will have three cases similar to (a)-(c). If we are in case (c) we consider $v_{2k-6}$, and continuing in this way get $v_{2k-l-1} \leq v_{2k-l}$ for some $l$, and we are in case (a) or (b).
Finally, consider paths of form (III). Suppose $v_{2k-4} = 2$. Then using the relations we obtain
\begin{eqnarray*}
[1 \cdots v_{2k-5} \underline{232}1] & = & -\frac{W_{222}}{W_{223}} [1 \cdots v_{2k-5} 2\underline{221}] - \frac{W_{122}}{W_{223}} [1 \cdots v_{2k-5} 2\underline{121}] \\
& = & \frac{W_{112}W_{222}}{W_{122}W_{223}} [1 \cdots v_{2k-5} \underline{221}1] + \frac{W_{111}W_{122}}{W_{112}W_{223}} [1 \cdots v_{2k-5} 2111] \\
& = & \left( \frac{W_{111}W_{122}}{W_{112}W_{223}} - \frac{W_{112}^2W_{222}}{W_{122}^2W_{223}} \right) [1 \cdots v_{2k-5} 2111],
\end{eqnarray*}
and we are in case (I). If $v_{2k-4} = 3,4$, we proceed as in (II).
Thus any path in $A(\mathcal{A}^{(2k+4)\ast},W)$ of length $2k+1$ is equal to $\xi u_{1\nu(1)}^{(k+1)} = \xi [111\cdots1]$, for some $\xi \in \mathbb{R}$ (possibly zero).

Thus for $\mathcal{A}^{(2m+2)\ast}$ we choose $u_{i\nu(i)} = u_{1\nu(1)} = [111\cdots1]$, and we have $v_{j\nu(i)} = v_{i\nu(k)} = [111\cdots1]$ of length $2m-2$. Then $v_{j\nu(i)} \nu(a_{ij}) = [111\cdots1] = u_{1\nu(1)} = a_{jk} v_{k\nu(j)}$, where $v_{k\nu(j)} = v_{1\nu(1)} = [111\cdots1]$. Similarly $a_{ki} v_{i\nu(k)} = v_{k\nu(j)} \nu(a_{jk})$. Thus we obtain $C=1$.
The situation for $\mathcal{A}^{(2m+1)\ast}$ follows similarly.

The graphs $\mathcal{D}^{(n)\ast}$ are (three-colourable) unfolded versions of the graphs $\mathcal{A}^{(n)\ast}$, where we replace every vertex $v$ of $\mathcal{A}^{(n)\ast}$ by three vertices $v_0$, $v_1$, $v_2$, where $v_a$ is of colour $a$, such that there are edges $v_0 \rightarrow w_1$, $v_1 \rightarrow w_2$ and $v_2 \rightarrow w_0$ if and only if there is an edge from $v$ to $w$ on $\mathcal{A}^{(n)\ast}$. The graph $\mathcal{D}^{(7)\ast}$ is illustrated in Figure \ref{fig:D(7)star}. Then the proof for $\mathcal{D}^{(n)\ast}$ follows exactly as the proof for $\mathcal{A}^{(n)\ast}$ where now we write a suffix for each vertex indicating the colour of that vertex. Due to the three-coloured nature of the graphs $\mathcal{D}^{(n)\ast}$, we see that the maximum path $u_{i\nu(i)} = u_{1_0\nu(1_0)}$ of length $n-3$ must end at the vertex $\nu(1_0) = 1_r$, where $r \equiv n \textrm{ mod } 3$. Thus the permutation $\nu$ should be the identity, clockwise rotation of the graph by $2 \pi /3$, or anticlockwise rotation of the graph by $2 \pi /3$, for $n \equiv 0, 1, 2 \textrm{ mod } 3$ respectively. Since $\nu$ is not always the identity for the graphs $\mathcal{D}^{(n)\ast}$ we require $\nu$-invariance of the cells here. The orbifold cell systems computed in \cite[Theorems 8.1 \& 8.2]{evans/pugh:2009i} are $\nu$-invariant.
If another inequivalent $\nu$-invariant cell system could be found such that $u_{i\nu(i)} = u_{1_0\nu(1_0)}$ is non-zero and the Hilbert series of $A$ was given by (\ref{eqn:Hilbert_Series-SU(3)ADE}), then the result $C=1$ would hold for this cell system also. \\

\newpage
\noindent \emph{The graphs $\mathcal{E}^{(8)}$, $\mathcal{E}^{(8)\ast}$ for the conformal embeddings $SU(3)_5 \subset SU(6)_1$ and $SU(3)_5 \subset SU(6)_1 \rtimes \mathbb{Z}_3$} \cite{xu:1998, bockenhauer/evans:1999i, evans/pugh:2009ii}:

\begin{figure}[tb]
\begin{minipage}[t]{7.5cm}
\begin{center}
  \includegraphics[width=35mm]{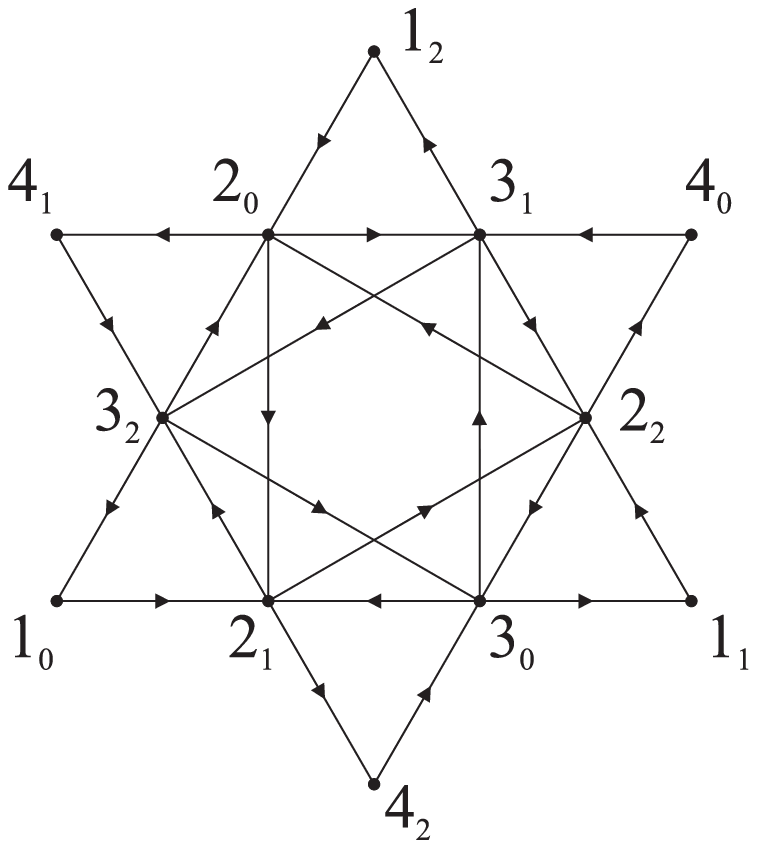}\\
  \caption{Graph $\mathcal{E}^{(8)}$} \label{fig:E(8)}
\end{center}
\end{minipage}
\hfill
\begin{minipage}[t]{7.5cm}
\begin{center}
  \includegraphics[width=15mm]{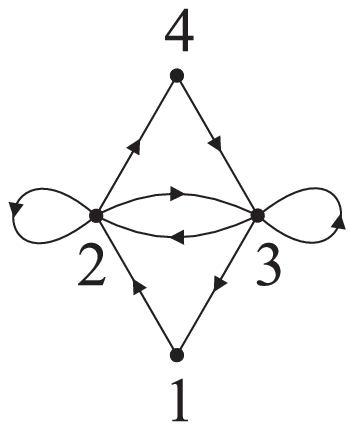}\\
  \caption{Graph $\mathcal{E}^{(8)\ast}$} \label{fig:E(8)star}
\end{center}
\end{minipage}
\end{figure}

We will first consider the graph $\mathcal{E}^{(8)\ast}$, illustrated in Figure \ref{fig:E(8)star}. The unique cell system $W$ (up to equivalence) was computed in \cite[Theorem 10.1]{evans/pugh:2009i}, and we use the same notation for the cells here. The quotient algebra $A$ has the relations
\begin{eqnarray}
& [123] = [231] = [324] = [432] = 0, & \label{eqn:rel-E8star-i} \\
& [222] = \frac{-W_{223}}{W_{222}} [232], \quad [333] = \frac{-W_{233}}{W_{333}} [323], & \label{eqn:rel-E8star-ii} \\
& \frac{-W_{123}}{W_{223}} [312] = [322] + [332], & \label{eqn:rel-E8star-iii} \\
& \frac{-W_{234}}{W_{223}} [243] = [223] + [233], & \label{eqn:rel-E8star-iv}
\end{eqnarray}
where $W_{222} = -W_{333}$ and $W_{223} = W_{233}$.

We choose the vertices $i=1$, $j=2$ and $k=3$. We choose the path $[122331]$ of length 5 to be the path $u_{i\nu(i)}$, which is nonzero since $[1\underline{223}31] = -(W_{243}/W_{223}) [124331]$ using (\ref{eqn:rel-E8star-iv}) and (\ref{eqn:rel-E8star-i}), and no relations can be used on $[124331]$ except for the one which transforms it back to $[122331]$. The only other relation which can be used on $[122331]$ is $[12\underline{233}1] = -(W_{243}/W_{233}) [122431]$. Then $v_{j\nu(i)} \nu(a_{ij}) = v_{2\nu(1)} \nu(a_{12}) = [223312]$, which we transform to the (non-zero) path $[233122]$:
\begin{eqnarray*}
[223\underline{312}] & = & -\frac{W_{223}}{W_{123}} [2\underline{233}22] - \frac{W_{233}}{W_{123}} [22\underline{333}2] = \frac{W_{223}}{W_{123}} [2\underline{232}22] + \frac{W_{233}^2}{W_{123}W_{333}} [2\underline{23232}] \\
& = & - 2\frac{W_{222}}{W_{123}} [222222],
\end{eqnarray*}
where in the penultimate equality we also used (\ref{eqn:rel-E8star-i}) and in the last equality we used (\ref{eqn:rel-E8star-ii}) once for $[2\underline{232}22]$ and twice for $[2\underline{23232}]$.
Similarly,
\begin{eqnarray*}
[23\underline{312}2] & = & -\frac{W_{223}}{W_{123}} [233222] - \frac{W_{233}}{W_{123}} [233322] = - 2\frac{W_{222}}{W_{123}} [222222],
\end{eqnarray*}
so that $v_{j\nu(i)} \nu(a_{ij}) = [223312] = [233122] = a_{jk} v'$, where $v' = [33122]$ is symmetric.
Here we need the fact that the coefficient of $[222222]$ in both the above equalities is not zero for the cell system $W$, as well as the fact that they are non-zero.

Since the graph $\mathcal{E}^{(8)}$, illustrated in Figure \ref{fig:E(8)}, is the (three-colourable) unfolded version of $\mathcal{E}^{(8)\ast}$, the result for $\mathcal{E}^{(8)}$ follows in the same way as the result for $\mathcal{D}^{(2n+1)\ast}$ follows from $\mathcal{A}^{(2n+1)\ast}$. The unique cell system $W$ (up to equivalence) was computed in \cite[Theorem 9.1]{evans/pugh:2009i}. Due to the three-coloured nature of $\mathcal{E}^{(8)}$, we see that since $u_{i\nu(i)} = u_{1_0\nu(1_0)}$ has length 5, $\nu$ must be the non-trivial permutation which sends $1_0 \rightarrow 1_2$. Since $\nu$ is not the identity for the graph $\mathcal{E}^{(8)}$ we require $\nu$-invariance of the cells here. The cell system $W$ for $\mathcal{E}^{(8)}$, computed in \cite[Theorem 9.1]{evans/pugh:2009i} is $\nu$-invariant. \\

\noindent \emph{The graph $\mathcal{E}_1^{(12)}$ for the conformal embedding $SU(3)_9 \subset (E_6)_1$} \cite{xu:1998, bockenhauer/evans:1999ii, evans/pugh:2009ii}:

\begin{figure}[tb]
\begin{center}
\includegraphics[width=100mm]{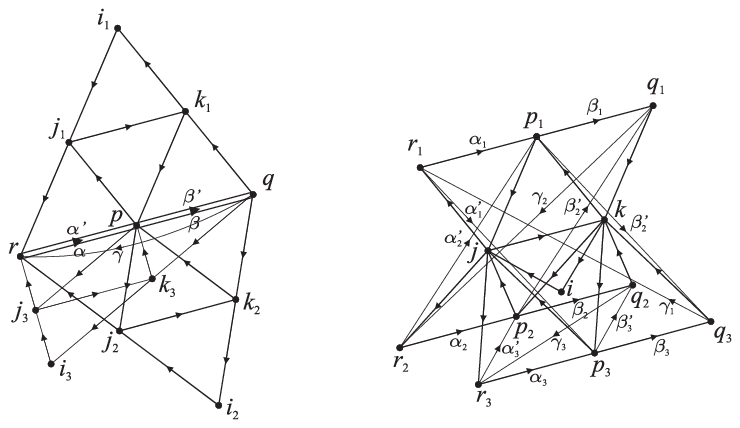}\\
 \caption{Graphs $\mathcal{E}_1^{(12)}$ and $\mathcal{E}_2^{(12)}$} \label{fig:E12(12)}
\end{center}
\end{figure}

Two inequivalent solutions for the cell system $W^{\pm}$ for the graph $\mathcal{E}_1^{(12)}$, illustrated in Figure \ref{fig:E12(12)}, were computed in \cite[Theorem 12.1]{evans/pugh:2009i}. We will use the solution $W^+$. The solution $W^-$ is obtained from $W^+$ by taking the conjugation of the graph $\mathcal{E}_1^{(12)}$, that is, $W^{-(abc)}_{ijk} = W^{+(\varphi(a)\varphi(b)\varphi(c))}_{\varphi(i)\varphi(j)\varphi(k)}$ where $\varphi$ is the map which reflects the graph $\mathcal{E}_1^{(12)}$ about the plane which passes through the vertices $i_1$, $i_2$, $i_3$ and $p$, and reverses the direction of each edge. For $\mathcal{E}_1^{(12)}$ the automorphism $\nu$ is the identity.
We choose the vertices $i=i_s$, $j=j_s$ and $k=k_s$, for some $s \in \{ 1,2,3 \}$.
We first computed a basis for the space of paths which start from the vertices $i_s$, $j_s$. The explicit details of these computations are given in the Appendix to \cite{evans/pugh:2010ii-arxiv}.
We will denote by $[rp]$, $[r'p]$ the path which goes along the edge $\alpha$, $\alpha'$ respectively, and similarly by $[pq]$, $[p'q]$ the path which goes along the edge $\beta$, $\beta'$ respectively.
Let $u_{i\nu(i)} = u_{i_s\nu(i_s)}$ be the (non-zero) path $u_{i\nu(i)} = [i_{s}j_{s}rpqrpqk_{s}i_{s}]$ of length 9.
Using the computations contained in the Appendix to \cite{evans/pugh:2010ii-arxiv}, we obtain
$v_{j\nu(i)} \nu(a_{ij}) = [j_{s}rpqrpqk_{s}i_{s}j_{s}] = d [j_{s}k_{s}pqk_{s}i_{s}j_{s}rpj_{s}] = d a_{jk} v'$,
where $v' = [k_{s}pqk_{s}i_{s}j_{s}rpj_{s}]$ is symmetric and
$d$ is a non-zero scalar. Hence $C=1$.
In the computations for the graph $\mathcal{E}_1^{(12)}$ we have used the orbifold cell system which was constructed explicitly in \cite{evans/pugh:2009i}. If another inequivalent cell system could be found which satisfied the non-vanishing of many coefficients which appear in the computations of a basis for $A$, and such that the Hilbert series of $A$ was given by (\ref{eqn:Hilbert_Series-SU(3)ADE}), then the result $C=1$ would hold for this cell system also. \\

\noindent \emph{The graph $\mathcal{E}_2^{(12)}$ for the conformal embedding $SU(3)_9 \subset (E_6)_1 \rtimes \mathbb{Z}_3$} \cite{bockenhauer/evans:2002, evans/pugh:2009ii}:

The graph $\mathcal{E}_2^{(12)}$, illustrated in Figure \ref{fig:E12(12)}, is a $\mathbb{Z}_3$ orbifold of $\mathcal{E}_1^{(12)}$.
Every cell system $W$ for $\mathcal{E}_2^{(12)}$ is equivalent to either a cell system $W^+$ or the inequivalent cell system $W^-$ \cite[Theorem 11.1]{evans/pugh:2009i}. We will use the solution $W^+$. The solution $W^-$ is obtained from $W^+$ by taking the conjugation of the graph $\mathcal{E}_1^{(12)}$, that is, $W^{-(abc)}_{ijk} = W^{+(\varphi(a)\varphi(b)\varphi(c))}_{\varphi(i)\varphi(j)\varphi(k)}$ where $\varphi$ is the map which reflects the graph $\mathcal{E}_2^{(12)}$ about the plane which passes through the vertices $i$, $p_1$, $p_2$ and $p_3$, and reverses the direction of each edge.
For the graph $\mathcal{E}_2^{(12)}$ the automorphism $\nu$ is the identity.
We choose $i$, $j$, $k$ as labelled on the graph $\mathcal{E}_2^{(12)}$ in Figure \ref{fig:E12(12)}.
We first computed a basis for the space of paths which start from the vertices $i$, $k$. The explicit details of these computations are given in the Appendix to \cite{evans/pugh:2010ii-arxiv}.
Here $q=e^{2 \pi i/12}$, and we will write $[m]$ for the quantum integer $[m]_q$. We let $a_{\pm}$, $b_{\pm}$ denote the scalars $a_{\pm} = \sqrt{[2][4]\pm\sqrt{[2][4]}}$, $b_{\pm} = \sqrt{[2]^2\pm\sqrt{[2][4]}}$.
Let $u_{i\nu(i)}$ be the path $[ijr_{1}p_{1}jkp_{1}q_{1}ki]$, which is symmetric. This path is non-zero in $A$ since
$$\left( \frac{\sqrt{[3][4]}}{\sqrt{[2]}} + a_- \frac{1}{\mu} \right) [\underline{ijr_{1}p_{1}jk}p_{1}q_{1}ki] = - a_+ [\underline{ijr_{1}p_{1}q_{1}kp_{1}}q_{1}ki] = - b_- [ijr_{1}p_{1}q_{1}r_{2}p_{1}q_{1}ki].$$
Then $v_{j\nu(i)} \nu(a_{ij})$ is given by $[kijr_{1}p_{1}jkp_{1}q_{1}k] = d [kp_{1}q_{1}kp_{1}jr_{1}p_{1}jk] = d a_{jk} v'$, where $v'$ is symmetric and $d$ is a non-zero scalar, and we obtain $C=1$. \\

\noindent \emph{The graph $\mathcal{E}_5^{(12)}$ for the twisted orbifold $\mathcal{D}^{(12)}$ invariant}:

\begin{figure}[tb]
\begin{center}
\includegraphics[width=70mm]{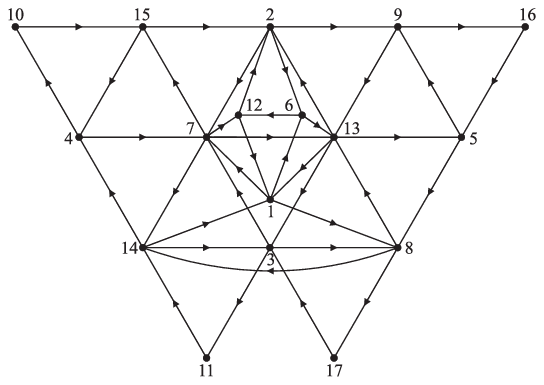}\\
 \caption{Labelled graph $\mathcal{E}_5^{(12)}$} \label{fig:E5(12)}
\end{center}
\end{figure}

The Moore-Seiberg invariant $Z_{\mathcal{E}_{MS}^{(12)}}$ is realised by a braided subfactor which produces the nimrep $\mathcal{E}_5^{(12)}$ \cite[Section 5.4]{evans/pugh:2009ii}, illustrated in Figure \ref{fig:E5(12)}. The unique cell system $W$ (up to equivalence) was computed in \cite[Theorem 13.1]{evans/pugh:2009i}.
For the graph $\mathcal{E}_5^{(12)}$ the automorphism $\nu$ is the identity. There are four possible choices for the vertices $j$, $k$: these are (3,8), (5,9), (14,3) and (15, 4). We see that it will not be possible to use the quicker method described above for the case where $\nu = \mathrm{id}$, since conjugating the graph will not interchange the vertices $j \leftrightarrow k$ for any of the choices of $j$, $k$.
We choose the vertices $i=10$, $j=15$ and $k=4$.
We first computed a basis for the space of paths which start from the vertices 10, 15, 4. The explicit details of these computations are given in the Appendix to \cite{evans/pugh:2010ii-arxiv}.
Here $q=e^{2 \pi i/12}$, and we will write $[m]$ for the quantum integer $[m]_q$.

We choose $u_{i\nu(i)} = u_{10\nu(10)} = [10,15,2,9,16,5,8,14,4,10]$. Then
\begin{eqnarray*}
v_{15\nu(10)} \nu(a_{10,15}) & = & [15,2,9,16,5,8,14,4,10,15] \;\; = \;\; c [15,4,7,12,1,7,13,1,7,15] \\
& = & c a_{15,4} v_{4\nu(15)},
\end{eqnarray*}
where $c = -[2]\sqrt{[3]^3}\left( \sqrt{[6]^3}+[3][4] \right)/[4]^2[6]$, and $C=1$ follows from
\begin{eqnarray*}
a_{4,10} v_{10\nu(4)} & = & [4,10,15,2,9,16,5,8,14,4] \;\; = \;\; c [4,7,12,1,7,13,1,7,15,4] \\
& = & c v_{4\nu(15)} \nu(a_{15,4}).
\end{eqnarray*}

\noindent \emph{The graph $\mathcal{E}^{(24)}$ for the conformal embedding $SU(3)_{21} \subset (E_7)_1$} \cite{xu:1998, bockenhauer/evans:1999ii, evans/pugh:2009ii}:

\begin{figure}[tb]
\begin{center}
\includegraphics[width=100mm]{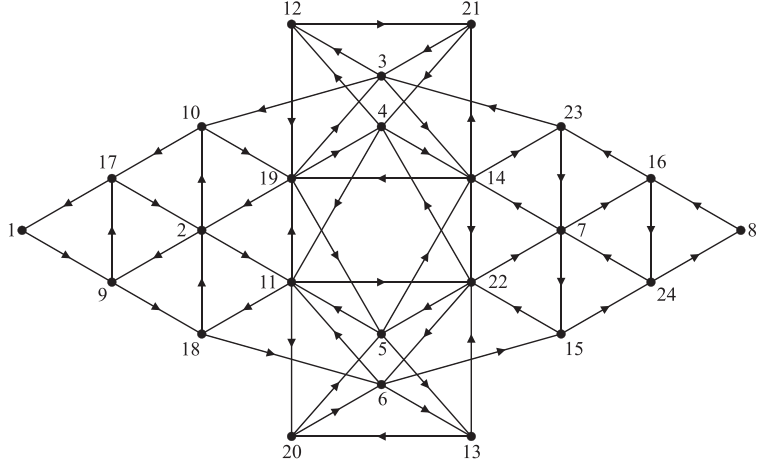}\\
 \caption{Labelled graph $\mathcal{E}^{(24)}$} \label{fig:E(24)}
\end{center}
\end{figure}

For the graph $\mathcal{E}^{(24)}$, illustrated in Figure \ref{fig:E(24)}, the automorphism $\nu$ is the identity.
The unique cell system $W$ (up to equivalence) was computed in \cite[Theorem 14.1]{evans/pugh:2009i}.
We choose the vertices $i=1$, $j=9$ and $k=17$.
We first computed a basis for the space of paths which start from the vertices 1, 9. The explicit details of these computations are given in the Appendix to \cite{evans/pugh:2010ii-arxiv}.
Here $q=e^{2 \pi i/24}$, and we will write $[m]$ for the quantum integer $[m]_q$.
Let $u$ be the path $[1,9,18,6,13,22,4,14,19,2,11,19,2,11,22,5,14,21,3,10,17,1]$, which is symmetric, and is non-zero as shown in the Appendix to \cite{evans/pugh:2010ii-arxiv}.
We choose $u_{i\nu(i)} = u_{1\nu(1)} := u$. Then $v_{j\nu(i)} \nu(a_{ij}) = v_{9\nu(1)} \nu(a_{1,9})$ is given by
\begin{eqnarray*}
\lefteqn{ [9,18,6,13,22,4,14,19,2,11,19,2,11,22,5,14,21,3,10,17,1,9] } \\
& = & -\frac{[4][9][10]\sqrt{[2][3]^5}}{[7]^2\sqrt{[5]^3}} [9,17,2,11,20,5,14,21,3,10,17,2,9,18,6,136,22,4,12,19,2,9] \\
& = & -\frac{[4][9][10]\sqrt{[2][3]^5}}{[7]^2\sqrt{[5]^3}} a_{jk} v',
\end{eqnarray*}
where $v'$ is symmetric. Then $C=1$.
\hfill
$\Box$

\subsection{Determining the Nakayama automorphism for $A(\mathcal{G},W)$}

We will now determine the Nakayama automorphism for the algebras $A=A(\mathcal{G},W)$, where $\mathcal{G}$ is an $SU(3)$ $\mathcal{ADE}$ graph which is not $\mathcal{E}_4^{(12)}$.

The algebra $A=A(\mathcal{G},W)$ is a Frobenius algebra, that is, there is a linear function $f:A \rightarrow \mathbb{C}$ such that $(x,y):=f(xy)$ is a non-degenerate bilinear form (this is equivalent to the statement that $A$ is isomorphic to its dual $A^{\ast} = \mathrm{Hom}(A,\mathbb{C})$ as left (or right) $A$-modules). There is an automorphism $\beta$ of $A$, called the Nakayama automorphism of $A$ (associated to $f$), such that $(x,y) = (y,\beta(x))$. Then there is an $A$-$A$ bimodule isomorphism $A^{\ast} \rightarrow {}_1 A_{\beta}$ \cite{yamagata:1996}.
Using the notation of Section \ref{sect:Hilbert-SU(3)}, we will define a non-degenerate form on $A$ by setting $f$ to be the function which is 0 on every element of $A$ of length $< h-3$, and 1 on $u_{i\nu(i)}$, for some $i \in \mathcal{G}_1$. Then using the relation $(x,y) = (y,\beta(x))$ this determines the value of $f$ on $u_{j\nu(j)}$, for all other $j \in \mathcal{G}_1$. We will normalize the $u_{j\nu(j)}$ such that $f(u_{j\nu(j)}) = 1$ for all $j \in \mathcal{G}_1$. From equation (\ref{eqn:Ia}) we see that $a^{\ast} = \lambda_a v_{j\nu(i)}^a$, $(v_{j\nu(i)}^a)^{\ast} = \mu_a \nu(a)$, where $a$ is an edge from $i$ to $j$ on $\mathcal{G}$.

\begin{Def} \label{def:NA}
Let $\beta$ denote the automorphism of $A$ defined on $\mathcal{G}$ by $\beta = \nu$, where $\nu$ is the permutation defined on $\mathcal{G}$ in Section \ref{sect:nakayama}.
\end{Def}

Let $W$ be a $\nu$-invariant cell system on $\mathcal{G}$. With $\beta$ defined above we have
$$\beta \left( \rho_{a} \right) \;\; = \;\; \sum_{b,b'} W_{i,r(b),k}^{(a,b,b')} \nu(b) \nu(b') \;\; = \;\; \sum_{b,b'} W_{\nu(i)\nu(r(b))\nu(k)}^{(\nu(a), \nu(b), \nu(b'))} \nu(b) \nu(b') \;\; = \;\; \rho_{\nu(a)},$$
for $a$ an edge from $k$ to $i$.
The following lemma shows that $\beta$ is the Nakayama automorphism for $A$:

\begin{Thm} \label{Thm:A^=1_A_beta}
The $A$-$A$ bimodule $A^{\ast}$ contains an element $u^{\ast}$ such that $A^{\ast} = A u^{\ast} = u^{\ast} A$, and $u^{\ast} a = \beta(a) u^{\ast}$ for all $a \in A$. Thus $A^{\ast}$ is isomorphic to ${}_1 A_{\beta}$ as an $A$-$A$ bimodule.
\end{Thm}

\noindent \emph{Proof}:
As in \cite[Corollary 4.7]{brenner/butler/king:2002}, any $u^{\ast} = \sum_i \zeta_i u^{\ast}_{\nu(i)i}$ with non-zero $\zeta_i \in \mathbb{C}$ will generate $A^{\ast}$ as both a left and right module. We will show that the $\zeta_i$ can be chosen such that $u^{\ast} a = \beta(a) u^{\ast}$ for each $a \in A_1$. By Definition \ref{def:NA} this becomes $\sum_l \zeta_l u^{\ast}_{\nu(l)l} a = \nu(a) \sum_l \zeta_l u^{\ast}_{\nu(l)l}$, so that $\zeta_i u^{\ast}_{\nu(i)i} a = \zeta_j \nu(a) u^{\ast}_{\nu(j)j} = \zeta_j (\lambda_{a}/\mu_{a}) u^{\ast}_{\nu(i)i} a$ by (\ref{eqn:IIa}). Thus we need to choose the $\zeta_i$ such that
\begin{equation} \label{eqn:III}
1 =  \frac{\zeta_j}{\zeta_i} \frac{\lambda_{a}}{\mu_{a}}.
\end{equation}
Similarly, for any triangle $i \stackrel{a}{\longrightarrow} j \stackrel{b}{\longrightarrow} k \stackrel{c}{\longrightarrow} i$ on $\mathcal{G}$, we need (\ref{eqn:III}) with $a$ replaced by $b$.
Then (\ref{eqn:III}) for $a$ and $b$ gives
$$1 = \frac{\zeta_i}{\zeta_k} \frac{\mu_{a}\mu_{b}}{\lambda_{a}\lambda_{b}} = \frac{\zeta_i}{\zeta_k} \frac{\lambda_{c}}{\mu_{c}},$$
by Propositions \ref{Prop:lam/mu=C} and \ref{Prop:C=1}.
Thus (\ref{eqn:III}) will also be satisfied with $i$, $j$ replaced by $k$, $i$ respectively. We may then consider a maximal connected subgraph $\mathcal{G}'$ of $\mathcal{G}$ such that for all vertices $i$, $j$ on $\mathcal{G}'$, with an edge $i \rightarrow j$, there is a choice of non-zero $\zeta_i$, $\zeta_j$ such that (\ref{eqn:III}) is satisfied. Suppose $\mathcal{G}' \neq \mathcal{G}$ so that there is a vertex $k$ of $\mathcal{G}'$ which is the source of an edge $\gamma$ to a vertex $l$ which is not on $\mathcal{G}'$. If we choose $\zeta_l = (\mu_{c}/\lambda_{c}) \zeta_k$ then (\ref{eqn:III}) is now satisfied for all vertices $i$, $j$ on the graph $\mathcal{G}''$, which is the graph obtained from $\mathcal{G}'$ by adjoining the vertex $l$ and the edge $c$. This contradicts the assumption that $\mathcal{G}'$ is the maximal subgraph such that (\ref{eqn:III}) is satisfied, hence $\mathcal{G}' = \mathcal{G}$.

Then there is an $A$-$A$ bimodule isomorphism $\phi: {}_1 A_{\beta} \rightarrow A^{\ast}$ given by $\phi(a) = a u^{\ast}$, such that $\phi(xa\beta(y)) = xa\beta(y)u^{\ast} = xau^{\ast}y = x\phi(a)y$, for all $a \in {}_1 A_{\beta}$, $x,y \in A$.
\hfill
$\Box$

\begin{Cor} \label{Cor:lambda=mu}
For every $a \in \mathcal{G}_1$ we have $\lambda_a = \mu_a$.
\end{Cor}

\noindent \emph{Proof}:
Using equation (\ref{eqn:Ia}), $1 = \lambda_a (a,v_{j\nu(i)}^a) = \lambda_a (v_{j\nu(i)}^a,\beta(a)) = \lambda_a/\mu_a$.
\hfill
$\Box$

\section{A finite resolution of $A$ as an $A$-$A$ bimodule} \label{sect:Omega^4(A)}

In this section we determine a finite resolution of $A = A(\mathcal{G},W)$ as an $A$-$A$ bimodule, where $\mathcal{G}$ is an $SU(3)$ $\mathcal{ADE}$ graph.

Let $A = A(\mathcal{G},W)$, where $\mathcal{G}$ is either an $SU(3)$ $\mathcal{ADE}$ graph or the McKay graph $\mathcal{G}_{\Gamma}$ of a finite subgroup $\Gamma \subset SU(3)$, and $W$ is a cell system on $\mathcal{G}$.
Consider the following complex:
\begin{equation} \label{seq:CY3_resolution}
A \otimes_R A[3] \stackrel{\mu_3}{\longrightarrow} A \otimes_R \widetilde{V} \otimes_R A[1] \stackrel{\mu_2}{\longrightarrow} A \otimes_R V \otimes_R A \stackrel{\mu_1}{\longrightarrow} A \otimes_R A \stackrel{\mu_0}{\longrightarrow} A \longrightarrow 0,
\end{equation}
where the $A$-$A$ bimodules $R$, $V$ and $\widetilde{V}$ are given by $R = (\mathbb{C}\mathcal{G})_0$, $V = (\mathbb{C}\mathcal{G})_1$ and $\widetilde{V} = \bigoplus_{a \in \mathcal{G}_1} \mathbb{C} \widetilde{a}$. We denote by $B[m]$ the space $B$ but with grading shifted by $m$.
The $A$-$A$ bimodule $\widetilde{V}$ is isomorphic to the dual space of $V$.
The connecting $A$-$A$ bimodule maps are given by
\begin{eqnarray*}
\mu_0(x \otimes y) & = & xy, \\
\mu_1(x \otimes a \otimes y) & = & x a \otimes y - x \otimes a y, \\
\mu_2(x \otimes \widetilde{a} \otimes y) & = & \sum_{b,b' \in \mathcal{G}_1} W_{abb'} (xb \otimes b' \otimes y + x \otimes b \otimes b'y), \\
\mu_3(x \otimes y) & = & \sum_{a \in \mathcal{G}_1} x a \otimes \widetilde{a} \otimes y - \sum_{a \in \mathcal{G}_1} x \otimes \widetilde{a} \otimes a y,
\end{eqnarray*}
where $x,y \in A$, $a \in \mathcal{G}_1$ and we denote by $W_{abb'}$ the cell $W(\triangle^{(a,b,b')}_{s(a),s(b),s(b')})$.

The sequence (\ref{seq:CY3_resolution}) is exact and the kernel of $\mu_3$ is zero if and only if $A$ is a Calabi-Yau algebra of dimension 3 \cite[Theorem 4.3]{bocklandt:2008}, e.g. if $\mathcal{G}$ is the McKay graph of a finite subgroup of $SU(3)$ -- see \cite{evans/pugh:2010i} for these subgroups and their graphs.
When $\mathcal{G}$ is an $SU(3)$ $\mathcal{ADE}$ graph, we will see that the kernel of $\mu_3$ is non-zero.

As in \cite{broomhead:2008}, applying the functor $\mathcal{F} = - \otimes_A R$ to the two-sided complex (\ref{seq:CY3_resolution}), we obtain the following one-sided complex of right $A$ modules:
\begin{equation} \label{seq:CY3_resolution-1sided}
A[3] \stackrel{\mathcal{F}(\mu_3)}{\longrightarrow} A \otimes_R \widetilde{V}[1] \stackrel{\mathcal{F}(\mu_2)}{\longrightarrow} A \otimes_R V \stackrel{\mathcal{F}(\mu_1)}{\longrightarrow} A \stackrel{\mathcal{F}(\mu_0)}{\longrightarrow} R \longrightarrow 0,
\end{equation}
where the connecting $A$ module maps are given by $\mathcal{F}(\mu_0)$ the projection of $A$ onto its zero-graded part $R$, $\mathcal{F}(\mu_1)(x \otimes a) = xa$, $\mathcal{F}(\mu_2)(x \otimes \widetilde{a}) = \sum_{b,b' \in \mathcal{G}_1} W_{abb'} xb \otimes b'$ and $\mathcal{F}(\mu_3)(x) = \sum_{a \in \mathcal{G}_1} xa \otimes \widetilde{a}$,
where $y \in A$, $a \in \mathcal{G}_1$.
The proof of \cite[Proposition 7.2.1]{broomhead:2008} carries over to our setting to show that the full complex (\ref{seq:CY3_resolution}) is exact if and only if the one-sided complex (\ref{seq:CY3_resolution-1sided}) is exact.

It was shown in Section \ref{sect:nimrepsSU(3)} that the one-sided sequence (\ref{eqn:SU(3)exact_seq}) is exact, for $k \leq h-2$.
Thus the bottom row of diagram (\ref{eqn:diag-d-F(mu)}) is shown to be exact everywhere except at the first term $A$ in degree $h-3$, where the Koszul differentials $d_i$
are given by $d_1(x \otimes a) = xa$, $d_2(x \otimes \widetilde{a}) = \sqrt{\mu_{s(a)}\mu_{r(a)}}^{-1} \sum_{b,b' \in \mathcal{G}_1} W_{abb'} xb \otimes b'$, and $d_3(x) = \sum_{a \in \mathcal{G}_1} (\sqrt{\mu_{r(a)}}/\sqrt{\mu_{s(a)}}) xa \otimes \widetilde{a}$.
There are isomorphisms $\psi_i$ such that the following diagram commutes:
\begin{equation} \label{eqn:diag-d-F(mu)}
\begin{array}{ccccccc}
A & \stackrel{\mathcal{F}(\mu_3)}{\longrightarrow} & A \otimes_R \widetilde{V} & \stackrel{\mathcal{F}(\mu_2)}{\longrightarrow} & A \otimes_R V & \stackrel{\mathcal{F}(\mu_1)}{\longrightarrow} & A \\
\downarrow \scriptstyle{\psi_3} & & \downarrow \scriptstyle{\psi_2} & & \downarrow \scriptstyle{\psi_1} & & \downarrow \scriptstyle{\psi_0} \\
A & \stackrel{d_3}{\longrightarrow} & A \otimes_R \widetilde{V} & \stackrel{d_2}{\longrightarrow} & A \otimes_R V & \stackrel{d_1}{\longrightarrow} & A \end{array}
\end{equation}
where $\psi_0 = \psi_1 = \mathrm{id}$, $\psi_2(x \otimes \widetilde{a}) = \sqrt{\mu_{s(a)}\mu_{r(a)}} \; x \otimes \widetilde{a}$ and $\psi_3(x) = \mu_{r(x)} x$.
Then the top row of diagram (\ref{eqn:diag-d-F(mu)}) is exact everywhere except at the first term $A$ in degree $h-3$, thus the complex (\ref{seq:CY3_resolution}) is shown to be exact everywhere except at the first term $A \otimes_R A[3]$. We will now compute the kernel of the map $\mu_3$, that is, the fourth syzygy $\Omega^4(A)$ of $A$.
The proof of the following theorem is based on the proof of Theorem 4.9 in \cite{brenner/butler/king:2002}.

\begin{Thm} \label{thm:Omega^4(A)}
Let $\mathcal{G}$ be an $SU(3)$ $\mathcal{ADE}$ graph which is not $\mathcal{E}_4^{(12)}$, and $W$ be a $\nu$-invariant cell system on $\mathcal{G}$. The fourth syzygy $\Omega^4(A) = \mathrm{Ker}(\mu_3)$ of $A=A(\mathcal{G},W)$ is isomorphic to ${}_1 A_{\beta^{-1}}$ as an $A$-$A$ bimodule, where $\beta$ is the Nakayama automorphism defined in Definition \ref{def:NA}. Thus
\begin{equation} \label{seq:q-CY3_resolution}
0 \rightarrow {}_1 A_{\beta^{-1}}[h] \stackrel{\mu_4}{\longrightarrow} A \otimes_R A[3] \stackrel{\mu_3}{\longrightarrow} A \otimes_R \widetilde{V} \otimes_R A[1] \stackrel{\mu_2}{\longrightarrow} A \otimes_R V \otimes_R A \stackrel{\mu_1}{\longrightarrow} A \otimes_R A \stackrel{\mu_0}{\longrightarrow} A \rightarrow 0
\end{equation}
is a finite resolution of $A$ as an $A$-$A$ bimodule, where $\mu_4(x) = x \sum_j x_j \otimes x_j^{\ast}$, the $x_j$ are a homogeneous basis for $A$ and $x_j^{\ast}$ is its dual basis under the non-degenerate form on $A$.
\end{Thm}

\noindent \emph{Proof}:
Since $A$ is a $(h-3,3)$-Koszul ring (see Section \ref{sect:nimrepsSU(3)}), by \cite[Theorem 3.15]{brenner/butler/king:2002} we see that $\Omega^4(A)$ is generated both as a left $A$-module and as a right $A$-module by its component $Z$ with total degree $h$, where
$$Z \subseteq (A \otimes_R A[3])_h = \sum_{r=0}^{h-3} A_r \otimes_R A_{h-3-r},$$
and that $Z \rightarrow A_0 \otimes_R A_{h-3}[3]$ is an $R$-$R$ bimodule isomorphism. Since $A_0 \cong R$ and $A_{h-3} \cong {}_1 R_{\nu}$ as $R$-$R$ bimodules, we see that $Z \cong {}_1 R_{\nu}$ as an $R$-$R$ bimodule, and has a basis given by elements $w_{i\nu(i)}$ which project onto $i \otimes u_{i\nu(i)}$ in the one-dimensional subspaces $i \cdot A_0 \otimes_R A_{h-3} \cdot \nu(i)$ of $A_0 \otimes_R A_{h-3}$. Since $\mu_3(w_{i\nu(i)}) = 0$ we must have
\begin{equation} \label{eqn:w(iv(i))}
w_{i\nu(i)} = i \sum_j x_j \otimes x_j^{\ast},
\end{equation}
where the $x_j$ are a homogeneous basis for $A$ and $x_j^{\ast}$ is its dual basis under the non-degenerate form on $A$, i.e. $x_j x_j^{\ast} = u_{s(x_j)\nu(s(x_j))}$.
To see this, consider the terms in $A_r \otimes_R A_{h-2-r}$ in $\mu_3(w_{i\nu(i)})$, for $r=1,\ldots,h-3$:
$$\sum_{\stackrel{a \in \mathcal{G}_1}{\scriptscriptstyle{y \in B_{r-1}^i}}} ya \otimes \widetilde{a} \otimes y^{\ast} - \sum_{\stackrel{a \in \mathcal{G}_1}{\scriptscriptstyle{y' \in B_{r}^i}}} y' \otimes \widetilde{a} \otimes a(y')^{\ast},$$
where $B_k^i$ denotes a basis of $i \cdot A_k$.
Now $y a = \sum_{x \in B_r^i} \lambda^{ya}_{x} x$ for any $y \in B_{r-1}^i$, $a \in \mathcal{G}_1$, where $\lambda^{ya}_{x} \in \mathbb{C}$.
Then for each $y' \in B_r^i$,
$$(\sum_{y \in B_{r-1}^i} y, a(y')^{\ast}) = (\sum_{y \in B_{r-1}^i} ya,(y')^{\ast}) = (\sum_{\stackrel{y \in B_{r-1}^i}{\scriptscriptstyle{x \in B_{r}^i}}} \lambda^{ya}_{x} x,(y')^{\ast}) = \sum_{\stackrel{y \in B_{r-1}^i}{\scriptscriptstyle{x \in B_{r}^i}}} \lambda^{ya}_{x} \delta_{y',x} = \sum_{y \in B_{r-1}^i} \lambda^{ya}_{y'},$$
so that $a(y')^{\ast} = \sum_{y \in B_{r-1}^i} \lambda^{ya}_{y'} y^{\ast}$. Then for fixed $a \in \mathcal{G}_1$ we have
$$\sum_{y \in B_{r-1}^i} y a \otimes \widetilde{a} \otimes y^{\ast} - \sum_{y' \in B_r^i} y' \otimes \widetilde{a} \otimes a(y')^{\ast} = \sum_{\stackrel{y \in B_{r-1}^i}{\scriptscriptstyle{y' \in B_{r}^i}}} \lambda^{ya}_{y'} y' \otimes \widetilde{a} \otimes y^{\ast} - \sum_{\stackrel{y \in B_{r-1}^i}{\scriptscriptstyle{y' \in B_{r}^i}}} y' \otimes \widetilde{a} \otimes \lambda^{ya}_{y'} y^{\ast} = 0,$$
as required.

Let $w = \sum_{i \in \mathcal{G}_0} w_{i\nu(i)}$. Then we have $Z = Rw = wR$ and by \cite[Theorem 3.15]{brenner/butler/king:2002} $\Omega^4(A) = Aw = wA$. Then there is an automorphism $\gamma$ of $A$ such that $xw = w\gamma(x)$ for all $x \in A$, which we will assume to have degree zero since $w$ is homogeneous. If we take $x=i \in A_0$ we see from (\ref{eqn:w(iv(i))}) that $\gamma(i) = \nu(i) = \beta(i)$, and so we must have $\gamma(a) = \gamma_a \nu(a) = \gamma_a \beta(a)$ for some $\gamma_a \in \mathbb{C} \setminus \{ 0 \}$.
If we now let $x = \nu(a) \in A_1$, we obtain $w_{i\nu(i)} \nu(a) = \gamma_{\nu(a)} a w_{j\nu(j)}$. Then by (\ref{eqn:w(iv(i))}) and comparing the terms in $A_1 \otimes_R A_{h-3}$, we have
$$\sum_{b \in \mathcal{G}_1} \lambda_{b} b \otimes v_{l\nu(i)}^b \nu(a) = \gamma_{\nu(a)} a \otimes u_{j\nu(j)},$$
using the identification $b^{\ast} = \lambda_b v_{l\nu(i)}^b$.
Using (\ref{eqn:Ia}), (\ref{eqn:Ib}) on the left hand side we obtain $(\lambda_a/\mu_a) a \otimes u_{j\nu(j)} = \gamma_{\nu(a)} a \otimes u_{j\nu(j)}$.
Thus $\gamma_{\nu(a)} = \lambda_a/\mu_a = 1$ for all $a \in \mathcal{G}_1$, using Corollary \ref{Cor:lambda=mu}, and we have $\gamma = \beta$ as required.
Then there is an $A$-$A$ bimodule isomorphism $\phi: {}_1 A_{\beta^{-1}} \rightarrow \Omega^4(A)$ given by $\phi(a) = aw$, such that $\phi(xa\beta^{-1}(y)) = xa\beta^{-1}yw = xawy = x\phi(a)y$, for all $a \in {}_1 A_{\beta^{-1}}$, $x,y \in A$.
The isomorphism $\phi$ thus defines an $A$-$A$ bimodule map $\mu_4: {}_1 A_{\beta^{-1}}[h] \longrightarrow A \otimes_R A[3]$ given by $\mu_4(x) = x \sum_j x_j \otimes x_j^{\ast}$.
\hfill
$\Box$

\paragraph{Acknowledgements}

This work was supported by the Marie Curie Research Training Network MRTN-CT-2006-031962 EU-NCG.
The authors would like to thanks Alastair King and the referees for their careful reading of earlier versions of this paper and their comments which have greatly improved the exposition.

\newpage
\begin{appendix}
\LARGE
\begin{center}
\textbf{Appendix to:} \\
\textbf{The Nakayama automorphism of the almost Calabi-Yau algebras associated to $SU(3)$ modular invariants}
\end{center}
\normalsize

\begin{center}
{David E. Evans and Mathew Pugh \\ $\,$ \\
        School of Mathematics, \\
        Cardiff University, \\
        Senghennydd Road, \\
        Cardiff, CF24 4AG, \\
        Wales, U.K.}
\end{center}

\normalsize

\section{Computation of certain basis paths for $A(\mathcal{G},W)$ for the graphs $\mathcal{E}_1^{(12)}$, $\mathcal{E}_2^{(12)}$, $\mathcal{E}_5^{(12)}$ and $\mathcal{E}^{(24)}$}

Here we provide the computations for a basis for (a subspace of) the quotient path space $A(\mathcal{G},W)$ for the exceptional $SU(3)$ $\mathcal{ADE}$ graphs $\mathcal{E}_1^{(12)}$, $\mathcal{E}_2^{(12)}$ and $\mathcal{E}^{(24)}$. These computations are needed for the proof of Proposition \ref{Prop:C=1} in \emph{The Nakayama automorphism of the almost Calabi-Yau algebras associated to $SU(3)$ modular invariants}. Due to the length of these computations, it was not practical to include them in the original paper, thus for completeness they are reproduced here instead. We do not write down a basis for the entire quotient space $A(\mathcal{G},W)$, but rather only for the subspaces $i \cdot A(\mathcal{G},W)$ and $j \cdot A(\mathcal{G},W)$ for two vertices $i$, $j$ of $\mathcal{G}$, where $i$ is a vertex which is the source of only one edge $a$ of $\mathcal{G}$, and $j$ is the range of this vertex.
For the graph $\mathcal{E}_5^{(12)}$ we also write a basis for the subspace $k \cdot A(\mathcal{G},W)$, where $k$ is the source of the unique edge $b$ whose range is the vertex $i$.
We also include details of the explicit computation of the constant $C=1$ for each of these graphs. These computations were only summarised in \cite{evans/pugh:2010ii}.

For each graph $\mathcal{G}$ we set $q=e^{2 \pi i/h}$, where $h$ is the Coxeter number of $\mathcal{G}$, and we will write $[m]$ for the quantum integer $[m]_q = (q^m - q^{-m})/(q-q^{-1})$.

\subsection{$\mathcal{E}_1^{(12)}$ graph for the conformal embedding $SU(3)_9 \subset (E_6)_1$}

Two inequivalent solutions for the cell system $W^{\pm}$ for the graph $\mathcal{E}_1^{(12)}$, illustrated in Figure \ref{fig:E12(12)}, were computed in \cite[Theorem 12.1]{evans/pugh:2009i}. We will use the solution $W^+$. The solution $W^-$ is obtained from $W^+$ by taking the conjugation of the graph $\mathcal{E}_1^{(12)}$, that is, $W^{-(abc)}_{ijk} = W^{+(\varphi(a)\varphi(b)\varphi(c))}_{\varphi(i)\varphi(j)\varphi(k)}$ where $\varphi$ is the map which reflects the graph $\mathcal{E}_1^{(12)}$ about the plane which passes through the vertices $i_1$, $i_2$, $i_3$ and $p$, and reverses the direction of each edge. For $\mathcal{E}_1^{(12)}$ the automorphism $\nu$ is the identity.
We choose $i=i_s$, $j=j_s$ and $k=k_s$, for some $s \in \{ 1,2,3 \}$.
We will write out a basis for the space of paths which start from the vertices $i_s$, $j_s$. We will denote by $l$ the length of the paths.
Suppose we have a basis for the paths of length $k$. We will consider every path of length $k+1$ obtainable by adding an extra edge to the end of each of the basis paths of length $k$. We will first list the basis paths of length $k+1$.
These will be followed by computations (contained within parentheses) showing how all the other paths obtained can be written in terms of these basis paths. We will mark basis paths with an asterisk, e.g. $[i_{s}j_{s}rpq]^*$.
We will underline the subpaths of length 2 on which we have used a relation at each stage, and if a subpath of length $k>2$ is underlined, this will indicate that we are using a relation on that path found when considering $l=k$.
Often when a relation is used on a path $b$ of length $k$, which gives $b$ as a linear combination of basis paths $b_i$ of length $k$, some of these paths $b_i$ are easily shown to be zero because the subpath given by the first $l<k$ edges of $b_i$ was shown to be zero when considering paths of length $l$. When this is the case, we will usually omit the paths $b_i$ which are known to be zero in this way.
We will denote by $[rp]$, $[r'p]$ the path which goes along the edge $\alpha$, $\alpha'$ respectively, and similarly by $[pq]$, $[p'q]$ the path which goes along the edge $\beta$, $\beta'$ respectively. Let $a_{\pm}$, $b_{\pm}$ denote the scalars $a_{\pm} = \sqrt{[2][4]\pm\sqrt{[2][4]}}$, $b_{\pm} = \sqrt{[2]^2\pm\sqrt{[2][4]}}$, $c_{\pm} = (a_{\pm}^2 b_{\pm})/(a_{\mp}^2 b_{\mp})$, and $\epsilon_s$ the third root of unity $\epsilon_s = e^{2\pi i (s-1)/3}$ for $s=1,2,3$. \\

\noindent
Paths starting at vertex $i_s$:
\scriptsize
\begin{eqnarray*}
l && \textrm{Paths} \\
\hline 2 && [i_{s}j_{s}r] \qquad ([i_{s}j_{s}k_{s}] = 0) \\
\hline 3 && [i_{s}\underline{j_{s}rp}] = -\epsilon_{s} \frac{a_-}{a_+} [i_{s}j_{s}r'p] -\overline{\epsilon_l} \frac{\sqrt{[3][4]}}{\sqrt{[2]}} \frac{1}{a_+} [\underline{i_{s}j_{s}k_{s}}p] = -\epsilon_{s} \frac{a_-}{a_+} [i_{s}j_{s}r'p] \\
\hline 4 && [i_{s}j_{s}\underline{rpj_{s}}] = -\epsilon_{s} \frac{a_-}{a_+} [i_{s}j_{s}r'pj_{s}], \qquad [i_{s}j_{s}rpq] \\
&& \left( [i_{s}j_{s}\underline{rpj_{s\pm1}}] = -\epsilon_{s\pm1} \frac{a_-}{a_+} [\underline{i_{s}j_{s}r'p}j_{s\pm1}] = \overline{\epsilon_{s}}\epsilon_{s\pm1} [i_{s}j_{s}rpj_{s\pm1}] \qquad \Rightarrow [i_{s}j_{s}rpj_{s\pm1}] = 0 \right) \\
&& \left( [i_{s}j_{s}\underline{rp'q}] = \frac{b_+}{b_-} [\underline{i_{s}j_{s}r'p}q] = -\overline{\epsilon_{s}} \frac{a_+ b_+}{a_- b_-} [i_{s}j_{s}rpq]^* \right) \\
\hline 5 && [i_{s}j_{s}rpqk_{s}], \qquad [i_{s}j_{s}r\underline{pj_{s}r}] = \epsilon_{s} \frac{[4]}{[2]} \frac{b_+}{a_-} [i_{s}j_{s}rpqr] \\
&& \left( [i_{s}j_{s}r\underline{pj_{s}k_{s}}] = -\epsilon_{s} \frac{\sqrt{[2]}}{\sqrt{[3][4]}} a_- [i_{s}j_{s}rpqk_{s}] - \overline{\epsilon_{s}} \frac{\sqrt{[2]}}{\sqrt{[3][4]}} a_+ [\underline{i_{s}j_{s}rp'q}k_{s}] = -\epsilon_{s} \frac{\sqrt{[2]}}{\sqrt{[3][4]}} a_- (1+c_{+}) [i_{s}j_{s}rpqk_{s}]^* \right) \\
&& \left( [i_{s}j_{s}r\underline{pqk_{s\pm1}}] = -\epsilon_{s} \frac{a_+}{a_-} [\underline{i_{s}j_{s}rp'q}k_{s\pm1}] = c_{+} [i_{s}j_{s}rpqk_{s\pm1}] \qquad \Rightarrow [i_{s}j_{s}rpqk_{s\pm1}] = 0 \right) \\
\hline 6 && [i_{s}j_{s}rp\underline{qrp}] = -\overline{\epsilon_{s}} \frac{[2]}{[4]} \frac{a_+}{b_-} [i_{s}j_{s}rpqk_{s}p] \\
&& \left( [i_{s}j_{s}rp\underline{qk_{s}i_{s}}] = 0 \right) \qquad \left( [i_{s}j_{s}rp\underline{qr'p}] = \epsilon_{s} \frac{[2]}{[4]} \frac{a_-}{b_+} [\underline{i_{s}j_{s}rpqk_{s}p}] = -\overline{\epsilon_{s}} \frac{a_- b_-}{a_+ b_+} [i_{s}j_{s}rpqrp]^* \right) \\
\hline 7 && [i_{s}j_{s}rpqrpq] \\
&& \left( [i_{s}j_{s}\underline{rp'q}] = \frac{b_+}{b_-} [i_{s}j_{s}rp\underline{qr'p}q] = \epsilon_{s} \frac{[2]}{[4]} \frac{a_-}{b_-} [\underline{i_{s}j_{s}rpqk_{s}p}q] = -\overline{\epsilon_{s}} \frac{a_-}{a_+} [i_{s}j_{s}rpqrpq]^* \right) \\
&& \left( [i_{s}j_{s}rpq\underline{rpj_{s\pm1}}] = -\epsilon_{s\pm1} \frac{a_-}{a_+} [i_{s}j_{s}rp\underline{qr'p}j_{s\pm1}] = - \epsilon_{s}\epsilon_{s\pm1} \frac{[2]}{[4]} \frac{a_-^2}{a_+ b_+} [\underline{i_{s}j_{s}rpk_{s}}pj_{s\pm1}] = -\overline{\epsilon_{s}}\epsilon_{s\pm1} c_{-} [i_{s}j_{s}rpqrpj_{s\pm1}] \right. \\
&& \qquad \Rightarrow [i_{s}j_{s}rpqrpj_{s\pm1}] = 0 \bigg) \qquad \left( [\underline{i_{s}j_{s}rpqrp}j_{s}] = -\overline{\epsilon_{s}} \frac{[2]}{[4]} \frac{a_+}{b_-} [i_{s}j_{s}rpq\underline{k_{s}pj_{s}}] = 0 \right) \\
\hline 8 && [i_{s}j_{s}rpqrpqk_{s}] \qquad \left( [i_{s}j_{s}rpqr\underline{pqr}] = 0 \right) \\
&& \left( [i_{s}j_{s}rpqr\underline{pqk_{s\pm1}}] = -\epsilon_{s\pm1} \frac{a_+}{a_-} [\underline{i_{s}j_{s}rpqrp'q}k_{s\pm1}] = \overline{\epsilon_{s}}\epsilon_{s\pm1} [i_{s}j_{s}rpqrpqk_{s\pm1}] \qquad \Rightarrow [i_{s}j_{s}rpqrpqk_{s\pm1}] = 0 \right) \\
\hline 9 && [i_{s}j_{s}rpqrpqk_{s}i_{s}] \qquad \left( [i_{s}j_{s}rpqrp\underline{qk_{s}p}] = 0 \right) \\
\hline 10 && \left( [i_{s}j_{s}rpqrpq\underline{k_{s}i_{s}j_{s}}] = 0 \right)
\end{eqnarray*}
\normalsize

\begin{figure}[tb]
\begin{center}
\includegraphics[width=100mm]{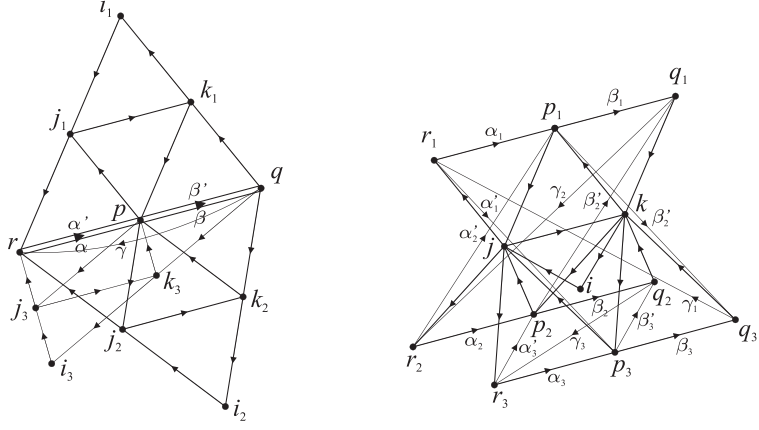}\\
 \caption{Graphs $\mathcal{E}_1^{(12)}$ and $\mathcal{E}_2^{(12)}$} \label{fig:E12(12)A}
\end{center}
\end{figure}

\noindent
Paths starting at vertex $j_s$:
\scriptsize

\normalsize

Let $u_{i\nu(i)} = u_{i_s\nu(i_s)}$ be the (non-zero) path $u_{i\nu(i)} = [i_{s}j_{s}rpqrpqk_{s}i_{s}]$ of length 9.
We have
\begin{eqnarray*}
[\underline{j_{s}k_{s}p}qk_{s}i_{s}j_{s}rpj_{s}] & = & -\epsilon_{s} \frac{\sqrt{[2]}}{\sqrt{[3][4]}} a_+ [j_{s}rpqk_{s}i_{s}j_{s}rpj_{s}] - \overline{\epsilon_{s}} \frac{\sqrt{[2]}}{\sqrt{[3][4]}} a_- [\underline{j_{s}r'pq}k_{s}i_{s}j_{s}rpj_{s}] \\
& = & -\epsilon_{s} \frac{\sqrt{[2]}}{\sqrt{[3][4]}} a_+ [j_{s}rpqk_{s}i_{s}j_{s}rpj_{s}] - \overline{\epsilon_{s}} \frac{\sqrt{[2]}}{\sqrt{[3][4]}} \frac{a_- b_-}{b_+} [\underline{j_{s}rp'qk_{s}i_{s}}j_{s}rpj_{s}] \\
& = & -\epsilon_{s} \frac{\sqrt{[2]}}{\sqrt{[3][4]}} a_+ (1-c_{-}) [\underline{j_{s}rpqk_{s}i_{s}j_{s}rp}j_{s}] \\
& = & -\epsilon_{s} \frac{\sqrt{[2]}}{\sqrt{[3][4]}} a_+ (1-c_{-}) \lambda^{(8)}_1 [j_{s}rpqk_{s}pq\underline{k_{s}pj_{s}}] \\
& = & \epsilon_{s} \frac{[2]}{[4]\sqrt{[3]}} a_+ (1-c_{-}) \lambda^{(8)}_1 [j_{s}rpqk_{s}pqk_{s}i_{s}j_{s}].
\end{eqnarray*}
Then
\begin{eqnarray*}
v_{j\nu(i)} \nu(a_{ij}) & = & [\underline{j_{s}rpqrp}qk_{s}i_{s}j_{s}] \\
& = & -\overline{\epsilon_{s}} \frac{[2]}{[4]} \frac{a_+}{b_-} [j_{s}rpqk_{s}pqk_{s}i_{s}j_{s}] - \overline{\epsilon_{s+1}} \frac{[2]}{[4]} \frac{a_+}{b_-} [\underline{j_{s}rpqk_{s+1}pqk_{s}}i_{s}j_{s}] \\
&& \qquad - \overline{\epsilon_{s-1}} \frac{[2]}{[4]} \frac{a_+}{b_-} [\underline{j_{s}rpqk_{s-1}pq}k_{s}i_{s}j_{s}] \\
& = & -\overline{\epsilon_{s}} \frac{[2]}{[4]} \frac{a_+}{b_-} [j_{s}rpqk_{s}pqk_{s}i_{s}j_{s}] - \overline{\epsilon_{s-1}} \frac{[2]}{[4]} \frac{a_+}{b_-} \frac{\lambda^{(5)}_1-\lambda^{(6)}_4}{\lambda^{(6)}_5} [\underline{j_{s}rpqk_{s+1}pqk_{s}}i_{s}j_{s}] \\
& = & -\overline{\epsilon_{s}} \frac{[2]}{[4]} \frac{a_+}{b_-} [j_{s}rpqk_{s}pqk_{s}i_{s}j_{s}] \\
& = & d [j_{s}k_{s}pqk_{s}i_{s}j_{s}rpj_{s}] \;\; = \;\; d a_{jk} v',
\end{eqnarray*}
where $v' = [k_{s}pqk_{s}i_{s}j_{s}rpj_{s}]$ is symmetric and $d = -\epsilon_{s} \sqrt{[3]} /b_- \lambda^{(8)}_1 (1-c_{-})$ is non-zero. Hence $C=1$. \\

\subsection{$\mathcal{E}_2^{(12)}$ graph for the conformal embedding $SU(3)_9 \subset (E_6)_1 \rtimes \mathbb{Z}_3$}

The graph $\mathcal{E}_2^{(12)}$, illustrated in Figure \ref{fig:E12(12)A}, is a $\mathbb{Z}_3$ orbifold of $\mathcal{E}_1^{(12)}$.
Every cell system $W$ for $\mathcal{E}_2^{(12)}$ is equivalent to either a cell system $W^+$ or the inequivalent cell system $W^-$ \cite[Theorem 11.1]{evans/pugh:2009i}. We will use the solution $W^+$. The solution $W^-$ is obtained from $W^+$ by taking the conjugation of the graph $\mathcal{E}_1^{(12)}$, that is, $W^{-(abc)}_{ijk} = W^{+(\varphi(a)\varphi(b)\varphi(c))}_{\varphi(i)\varphi(j)\varphi(k)}$ where $\varphi$ is the map which reflects the graph $\mathcal{E}_2^{(12)}$ about the plane which passes through the vertices $i$, $p_1$, $p_2$ and $p_3$, and reverses the direction of each edge.
For the graph $\mathcal{E}_2^{(12)}$ the automorphism $\nu$ is the identity.
We choose $i$, $j$, $k$ as labelled on the graph $\mathcal{E}_2^{(12)}$ in Figure \ref{fig:E12(12)A}.
We will write out a basis for the space of paths which start from the vertices $i$, $k$. Let $a_{\pm}$, $b_{\pm}$ denote the scalars $a_{\pm} = \sqrt{[2][4]\pm\sqrt{[2][4]}}$, $b_{\pm} = \sqrt{[2]^2\pm\sqrt{[2][4]}}$ and $c_{\pm} = (a_{\pm}^2 b_{\pm})/(a_{\mp}^2 b_{\mp})$. \\

\noindent
Paths starting at vertex $i$:
\scriptsize
\begin{eqnarray*}
l && \textrm{Paths} \\
\hline 2 && [ijr_{1}], \qquad [ijr_{2}], \qquad [ijr_{3}] \qquad \left( [ijk] = 0 \right) \\
\hline 3 && [i\underline{jr_{1}p_{1}}] = -\frac{a_+}{a_-} [ijr_{2}p_{1}], \qquad [i\underline{jr_{2}p_{2}}] = -\frac{a_+}{a_-} [ijr_{3}p_{2}], \qquad [i\underline{jr_{3}p_{3}}] = -\frac{a_+}{a_-} [ijr_{1}p_{3}] \\
\hline 4 && [\underline{ijr_{1}p_{1}}q_{1}] = -\frac{a_+}{a_-} [ij\underline{r_{2}p_{1}q_{1}}] = -\frac{a_+ b_+}{a_- b_-} [\underline{ijr_{2}p_{2}}q_{1}] = c_{+} [ijr_{3}p_{2}q_{1}], \\
&& [\underline{ijr_{2}p_{2}}q_{2}] = -\frac{a_+}{a_-} [ij\underline{r_{3}p_{2}q_{2}}] = -\frac{a_+ b_+}{a_- b_-} [\underline{ijr_{3}p_{3}}q_{2}] = c_{+} [ijr_{1}p_{3}q_{2}], \\
&& [\underline{ijr_{3}p_{3}}q_{3}] = -\frac{a_+}{a_-} [ij\underline{r_{1}p_{3}q_{3}}] = -\frac{a_+ b_+}{a_- b_-} [\underline{ijr_{1}p_{1}}q_{3}] = c_{+} [ijr_{2}p_{1}q_{3}], \\
&& [ij\underline{r_{1}p_{1}j}] = -\frac{a_+}{a_-} [\underline{ijr_{1}p_{3}}j] = [ij\underline{r_{3}p_{3}j}] = -\frac{a_+}{a_-} [\underline{ijr_{3}p_{2}}j] = [ij\underline{r_{2}p_{2}j}] = -\frac{a_+}{a_-} [ijr_{2}p_{1}j] \\
\hline 5 && [ijr_{1}\underline{p_{1}q_{1}r_{2}}] = \frac{a_+}{b_-} [ijr_{1}p_{1}jr_{2}], \qquad [ijr_{2}\underline{p_{2}q_{2}r_{3}}] = \frac{a_+}{b_-} [ijr_{2}p_{2}jr_{3}], \qquad [ijr_{3}\underline{p_{3}q_{3}r_{1}}] = \frac{a_+}{b_-} [ijr_{3}p_{3}jr_{1}], \qquad [ijr_{1}p_{1}q_{1}k] \\
&& \bigg( \textrm{Now } [ijr_{1}\underline{p_{1}q_{1}k}] = -\frac{\sqrt{[3][4]}}{\sqrt{[2]}} \frac{1}{a_+} [\underline{ijr_{1}p_{1}j}k] - \frac{a_-}{a_+} [ij\underline{r_{1}p_{1}q_{3}}k] = \frac{\sqrt{[3][4]}}{\sqrt{[2]}} \frac{1}{a_-} [ijr_{1}\underline{p_{3}jk}] - \frac{a_- b_-}{a_+ b_+} [ijr_{1}p_{3}q_{3}k] \\
&& \qquad = -[\underline{ijr_{1}p_{3}q_{2}}k] - \frac{a_+}{a_-} (1+c_{-}) [\underline{ijr_{1}p_{3}q_{3}}k] = -c_{-} [ijr_{2}\underline{p_{2}q_{2}k}] + (1+c_{+}) [i\underline{jr_{2}p_{1}}q_{3}k] \\
&& \qquad = \frac{a_-}{a_+} c_{-} [\underline{ijr_{2}p_{2}}q_{1}k] + \frac{\sqrt{[3][4]}}{\sqrt{[2]}} \frac{1}{a_+} c_{-} [\underline{ijr_{2}p_{2}}jk] - \frac{a_-}{a_+} (1+c_{+}) [ijr_{1}p_{1}q_{3}k] \\
&& \qquad = -c_{-}^2 [ijr_{1}\underline{p_{1}q_{1}k}] + \frac{\sqrt{[3][4]}}{\sqrt{[2]}} \frac{1}{a_+} c_{-} [ijr_{1}p_{1}jk] - \frac{a_-}{a_+} (1+c_{+}) [ijr_{1}p_{1}q_{3}k] \\
&& \qquad = \frac{a_-}{a_+} \left( c_{-}^2 - c_{+} - 1 \right) [ijr_{1}p_{1}q_{3}k] + \frac{\sqrt{[3][4]}}{\sqrt{[2]}} \frac{1}{a_+} \left( c_{-}^2 + c_{-} \right) [ijr_{1}p_{1}jk], \\
&& \textrm{ and comparing this with the first equality } [ijr_{1}p_{1}q_{1}k] = -\frac{\sqrt{[3][4]}}{\sqrt{[2]}} \frac{1}{a_+} [ijr_{1}p_{1}jk] - \frac{a_-}{a_+} [ijr_{1}p_{1}q_{3}k] \\
&& \textrm{ we obtain } [ijr_{1}p_{1}jk] = \mu [ijr_{1}p_{1}q_{3}k], \quad \textrm{ where } \mu = \frac{\sqrt{[2]}}{\sqrt{[3][4]}} a_- \frac{c_{+} - c_{-}^2}{1 + c_{-} + c_{-}^2}. \\
&& \textrm{Similarly } [ijr_{2}p_{2}jk] = \mu [ijr_{2}p_{2}q_{1}k] \textrm{ and } [ijr_{3}p_{3}jk] = \mu [ijr_{3}p_{3}q_{2}k]. \\
&& \textrm{Thus } [ijr_{1}\underline{p_{1}q_{1}k}] = -\frac{\sqrt{[3][4]}}{\sqrt{[2]}} \frac{1}{a_+} [ijr_{1}p_{1}jk] - \frac{a_-}{a_+} [\underline{ijr_{1}p_{1}q_{3}k}] = \left( -\frac{\sqrt{[3][4]}}{\sqrt{[2]}} \frac{1}{a_+} - \frac{a_-}{a_+} \frac{1}{\mu} \right) [\underline{ijr_{1}p_{1}j}k] \\
&& \qquad = \left( -\frac{\sqrt{[3][4]}}{\sqrt{[2]}} \frac{1}{a_+} - \frac{a_-}{a_+} \frac{1}{\mu} \right) [\underline{ijr_{2}p_{2}jk}] = [ijr_{2}p_{2}q_{2}k], \\
&& \textrm{ so } [ijr_{2}p_{2}q_{2}k] = [ijr_{1}p_{1}q_{1}k]^*, \qquad \textrm{ and similarly } [ijr_{3}p_{3}q_{3}k] = [ijr_{1}p_{1}q_{1}k]^* \bigg) \\
&& \bigg( \textrm{From the previous computation we have } [ijr_{1}p_{1}jk] = \left( -\frac{\sqrt{[3][4]}}{\sqrt{[2]}} \frac{1}{a_+} - \frac{a_-}{a_+} \frac{1}{\mu} \right)^{-1} [ijr_{1}p_{1}q_{1}k]^* \bigg) \\
&& \bigg( \textrm{Similarly } [ijr_{2}p_{2}jk] = \left( -\frac{\sqrt{[3][4]}}{\sqrt{[2]}} \frac{1}{a_+} - \frac{a_-}{a_+} \frac{1}{\mu} \right)^{-1} [ijr_{2}p_{2}q_{2}k] = \left( -\frac{\sqrt{[3][4]}}{\sqrt{[2]}} \frac{1}{a_+} - \frac{a_-}{a_+} \frac{1}{\mu} \right)^{-1} [ijr_{1}p_{1}q_{1}k]^*, \\
&& \textrm{ and } [ijr_{3}p_{3}jk] = \left( -\frac{\sqrt{[3][4]}}{\sqrt{[2]}} \frac{1}{a_+} - \frac{a_-}{a_+} \frac{1}{\mu} \right)^{-1} [ijr_{1}p_{1}q_{1}k]^* \bigg) \\
\hline 6 && [ijr_{1}p_{1}\underline{q_{1}r_{2}p_{1}}] = \frac{a_+}{b_-} [\underline{ijr_{1}p_{1}q_{1}k}p_{1}] = \frac{a_+}{b_-} [ijr_{3}p_{3}\underline{q_{3}kp_{1}}] = [ijr_{3}p_{3}q_{3}r_{1}p_{1}], \\
&& [ijr_{2}p_{2}\underline{q_{2}r_{3}p_{2}}] = \frac{a_+}{b_-} [\underline{ijr_{2}p_{2}q_{2}k}p_{2}] = \frac{a_+}{b_-} [ijr_{1}p_{1}\underline{q_{1}kp_{2}}] = [ijr_{1}p_{1}q_{1}r_{2}p_{2}], \\
&& [ijr_{3}p_{3}\underline{q_{3}r_{1}p_{3}}] = \frac{a_+}{b_-} [\underline{ijr_{3}p_{3}q_{3}k}p_{3}] = \frac{a_+}{b_-} [ijr_{2}p_{2}\underline{q_{2}kp_{3}}] = [ijr_{2}p_{2}q_{2}r_{3}p_{3}] \\
&& \left( [\underline{ijr_{1}p_{1}q_{1}k}i] = \left( -\frac{\sqrt{[3][4]}}{\sqrt{[2]}} \frac{1}{a_+} - \frac{a_-}{a_+} \frac{1}{\mu} \right) [ijr_{1}p_{1}\underline{jki}] = 0 \right) \\
\hline 7 && [ijr_{1}p_{1}q_{1}\underline{r_{2}p_{1}q_{1}}] = \frac{b_+}{b_-} [\underline{ijr_{1}p_{1}q_{1}r_{2}p_{2}}q_{1}] = \frac{b_+}{b_-} [ijr_{2}p_{2}q_{2}r_{3}p_{2}q_{1}], \\
&& [ijr_{2}p_{2}q_{2}\underline{r_{3}p_{2}q_{2}}] = \frac{b_+}{b_-} [\underline{ijr_{2}p_{2}q_{2}r_{3}p_{3}}q_{2}] = \frac{b_+}{b_-} [ijr_{3}p_{3}q_{3}r_{1}p_{3}q_{2}], \\
&& [ijr_{3}p_{3}q_{3}\underline{r_{1}p_{3}q_{3}}] = \frac{b_+}{b_-} [\underline{ijr_{3}p_{3}q_{3}r_{1}p_{1}}q_{3}] = \frac{b_+}{b_-} [ijr_{1}p_{1}q_{1}r_{2}p_{1}q_{3}], \\
&& \bigg( \textrm{Now } [\underline{ijr_{1}p_{1}q_{1}r_{2}p_{1}}j] = \frac{a_+}{b_-} [ijr_{1}p_{1}q_{1}kp_{1}j], \textrm{ but also } [ijr_{1}p_{1}q_{1}\underline{r_{2}p_{1}j}] = -\frac{a_-}{a_+} [\underline{ijr_{1}p_{1}q_{1}r_{2}p_{2}}j] = -\frac{a_-}{b_-} [ijr_{1}p_{1}q_{1}kp_{2}j], \\
&& \textrm{ and } [ijr_{1}p_{1}q_{1}\underline{r_{2}p_{1}j}] = -\frac{a_-}{a_+} [\underline{ijr_{1}p_{1}q_{1}r_{2}p_{2}}j] = -\frac{a_-}{a_+} [ijr_{2}p_{2}q_{2}\underline{r_{3}p_{2}j}] = \frac{a_-^2}{a_+^2} [\underline{ijr_{2}p_{2}q_{2}r_{3}p_{3}}j] = \frac{a_-^2}{a_+ b_-} [\underline{ijr_{2}p_{2}q_{2}k}p_{3}j] \\
&& \qquad = \frac{a_-^2}{a_+ b_-} [ijr_{1}p_{1}q_{1}kp_{3}j]. \textrm{ Now since } [ijr_{1}p_{1}q_{1}\underline{kij}] = 0, \\
&& -\frac{\sqrt{[4]}}{\sqrt{[2]}} [ijr_{1}p_{1}q_{1}\underline{kij}] = [\underline{ijr_{1}p_{1}q_{1}kp_{1}j}] + [\underline{ijr_{1}p_{1}q_{1}kp_{2}j}] + [\underline{ijr_{1}p_{1}q_{1}kp_{3}j}] = \frac{b_-}{a_-} \left( \frac{a_-}{a_+} + \frac{a_+ }{a_-} - 1 \right) [ijr_{1}p_{1}q_{1}r_{2}p_{1}j] = 0, \\
&& \textrm{ so } [ijr_{1}p_{1}q_{1}r_{2}p_{1}j] = [ijr_{2}p_{2}q_{2}r_{3}p_{2}j] = [ijr_{3}p_{3}q_{3}r_{1}p_{3}j] = 0. \bigg) \\
\hline 8 && [ijr_{1}p_{1}q_{1}r_{2}\underline{p_{1}q_{1}k}] = -\frac{a_-}{a_+} [\underline{ijr_{1}p_{1}q_{1}r_{2}p_{1}q_{3}}k] = -\frac{a_- b_-}{a_+ b_+} [ijr_{3}p_{3}q_{3}r_{1}\underline{p_{3}q_{3}k}] = c_{-} [\underline{ijr_{3}p_{3}q_{3}r_{1}p_{3}q_{2}}k] \\
&& \qquad = \frac{b_-}{b_+} c_{-} [ijr_{2}p_{2}q_{2}r_{3}p_{2}q_{2}k] \qquad \left( [ijr_{1}p_{1}q_{1}r_{2}\underline{p_{1}q_{1}r_{2}}] = \frac{a_+}{b_-} [\underline{ijr_{1}p_{1}q_{1}r_{2}p_{1}j}r_{2}] = 0 \right) \\
&& \left( [ijr_{2}p_{2}q_{2}r_{3}\underline{p_{2}q_{2}r_{3}}] = \frac{a_+}{b_-} [\underline{ijr_{2}p_{2}q_{2}r_{3}p_{2}j}r_{3}] = 0 \right) \qquad \left( [ijr_{3}p_{3}q_{3}r_{1}\underline{p_{3}q_{3}r_{1}}] = \frac{a_+}{b_-} [\underline{ijr_{3}p_{3}q_{3}r_{1}p_{3}j}r_{1}] = 0 \right) \\
\hline 9 && [ijr_{1}p_{1}q_{1}r_{2}p_{1}q_{1}ki] \qquad \left( [ijr_{1}p_{1}q_{1}r_{2}p_{1}\underline{q_{1}kp_{1}}] = \frac{b_-}{a_+} [\underline{ijr_{1}p_{1}q_{1}r_{2}p_{1}q_{1}r_{2}}p_{1}] = 0 \right) \\
&& \left( [\underline{ijr_{1}p_{1}q_{1}r_{2}p_{1}q_{1}k}p_{2}] = \frac{b_-}{b_+} c_{-} [ijr_{2}p_{2}q_{2}r_{3}p_{2}\underline{q_{2}kp_{2}}] = \frac{b_-^2}{a_+ b_+} c_{-} [\underline{ijr_{2}p_{2}q_{2}r_{3}p_{2}q_{2}r_{3}}p_{2}] = 0 \right) \\
&& \left( [\underline{ijr_{1}p_{1}q_{1}r_{2}p_{1}q_{1}k}p_{3}] = \frac{a_- b_-}{a_+ b_+} [ijr_{3}p_{3}q_{3}r_{1}p_{3}\underline{q_{3}kp_{3}}] = \frac{b_-}{a_-} c_{-} [\underline{ijr_{3}p_{3}q_{3}r_{1}p_{3}q_{3}r_{1}}p_{3}] = 0 \right) \\
\hline 10 && \left( [ijr_{1}p_{1}q_{1}r_{2}p_{1}q_{1}\underline{kij}] = 0 \right)
\end{eqnarray*}
\normalsize

\noindent
Paths starting at vertex $k$:
\scriptsize

\normalsize

Let $u$ be the path $[ijr_{1}p_{1}jkp_{1}q_{1}ki]$, which is symmetric. This path is non-zero in $A$ since
$$\left( \frac{\sqrt{[3][4]}}{\sqrt{[2]}} + a_- \frac{1}{\mu} \right) [\underline{ijr_{1}p_{1}jk}p_{1}q_{1}ki] = - a_+ [\underline{ijr_{1}p_{1}q_{1}kp_{1}}q_{1}ki] = - b_- [ijr_{1}p_{1}q_{1}r_{2}p_{1}q_{1}ki].$$
We choose $u_{i\nu(i)} := u$. Then $v_{j\nu(i)} \nu(a_{ij})$ is given by
\begin{eqnarray*}
\lefteqn{ [\underline{kij}r_{1}p_{1}jkp_{1}q_{1}k] } \\
& \qquad = & -\frac{\sqrt{[2]}}{\sqrt{[4]}} [\underline{kp_{1}jr_{1}}p_{1}jkp_{1}q_{1}k] - \frac{\sqrt{[2]}}{\sqrt{[4]}} [kp_{2}jr_{1}p_{1}jkp_{1}q_{1}k] - \frac{\sqrt{[2]}}{\sqrt{[4]}} [\underline{kp_{3}jr_{1}}p_{1}jkp_{1}q_{1}k] \\
& \qquad = & -\frac{\sqrt{[2]}}{\sqrt{[4]}} \frac{b_-}{a_+} (1+c_{+}) [\underline{kp_{3}q_{3}r_{1}p_{1}}jkp_{1}q_{1}k] - \frac{\sqrt{[2]}}{\sqrt{[4]}} [\underline{kp_{2}jr_{1}p_{1}}jkp_{1}q_{1}k] \\
& \qquad = & -\frac{\sqrt{[2]}}{\sqrt{[4]}} \frac{a_- b_-}{a_+ b_+} (1+c_{+}) [\underline{kp_{3}q_{3}kp_{1}jkp_{1}}q_{1}k] - \frac{\sqrt{[2]}}{\sqrt{[4]}} \frac{a_+}{a_-} [kp_{2}q_{2}kp_{1}jkp_{1}q_{1}k] \\
&& \qquad + \frac{\sqrt{[2]}}{\sqrt{[4]}} \frac{a_+}{a_-} (1+c_{+}) [kp_{1}q_{1}kp_{1}jkp_{1}q_{1}k] \\
& \qquad = & \frac{\sqrt{[2]}}{\sqrt{[4]}} \frac{a_+}{a_-} (1+c_{+}) \left( 1-\frac{\lambda_{1}}{\lambda_{3}}c_{-} \right) [kp_{1}q_{1}kp_{1}jkp_{1}q_{1}k] \\
&& \qquad - \frac{\sqrt{[2]}}{\sqrt{[4]}} \frac{a_+}{a_-} \left( 1-\frac{\lambda_{2}}{\lambda_{3}}c_{-}(1+c_{+}) \right) [\underline{kp_{2}q_{2}kp_{1}jkp_{1}q_{1}}k] \\
& \qquad = & d_1 [kp_{1}q_{1}kp_{1}jkp_{1}q_{1}k], \\
\lefteqn{ \textrm{ where } d_1 = \frac{\sqrt{[2]}}{\sqrt{[4]}} \frac{a_+}{a_-} \left( (1+c_{+}) \left( 1-\frac{\lambda_{1}}{\lambda_{3}}c_{-} \right) - \mu_2 \left( 1-\frac{\lambda_{2}}{\lambda_{3}}c_{-}(1+c_{+}) \right) \right) \neq 0. } \\
\lefteqn{ \textrm{Now } \frac{\sqrt{[2]}}{\sqrt{[3][4]}} (1-\mu_1^2) a_- [\underline{kp_{1}q_{1}kp_{1}jr_{1}}p_{1}jk] } \\
& \quad = & \frac{\sqrt{[2]}}{\sqrt{[3][4]}} (1-\mu_1^2) \frac{a_+ b_+}{a_- b_-} [kp_{1}q_{1}kp_{3}\underline{q_{3}r_{1}p_{1}}jk] \\
& \quad = & -\frac{\sqrt{[2]}}{\sqrt{[3][4]}} (1-\mu_1^2) a_+ [\underline{kp_{1}q_{1}kp_{3}q_{3}k}p_{1}jk] \\
& \quad = & \mu_1^3 (1+c_{+}) [kp_{1}q_{1}kp_{1}jkp_{1}jk] - \mu_1 [kp_{2}q_{2}kp_{1}jkp_{1}jk] + \mu_1^2 (1+c_{-}) [\underline{kp_{3}q_{3}kp_{1}jkp_{1}}jk] \\
& \quad = & \mu_1^2 \left( \mu_1 (1+c_{+}) - \frac{\lambda_{1}}{\lambda_{3}} (1+c_{-}) \right) [kp_{1}q_{1}kp_{1}jkp_{1}jk] \\
&& \qquad - \mu_1 \left( 1 + \frac{\lambda_{1}}{\lambda_{3}} \mu_1 (1+c_{-}) \right) [\underline{kp_{2}q_{2}kp_{1}jkp_{1}j}k] \\
& \quad = & \mu_1 \left( \mu_1^2 (1+c_{+}) - \frac{\lambda_{1}}{\lambda_{3}} \mu_1 (1+c_{-}) - \lambda_{3}^2 \frac{\lambda_{3} + \lambda_{1} \mu_1 (1+c_{-})}{\lambda_{3}^3 + \lambda_{2}^3 (1+c_{-})^3} d_2 \right) [\underline{kp_{1}q_{1}kp_{1}jkp_{1}jk}] \\
& \quad = & d_4 [kp_{1}q_{1}kp_{1}jkp_{1}q_{1}k] \\
\lefteqn{ \textrm{where } d_4 = \mu_1 \left( \mu_1^2 (1+c_{+}) - \frac{\lambda_{1}}{\lambda_{3}} \mu_1 (1+c_{-}) - \lambda_{3}^2 \frac{\lambda_{3} + \lambda_{1} \mu_1 (1+c_{-})}{\lambda_{3}^3 + \lambda_{2}^3 (1+c_{-})^3} d_2 \right) d_3^{-1} \neq 0, } \\
\lefteqn{ \textrm{with } d_2 = \left( \frac{\lambda_{1}^{(2)}}{\lambda_{2}^{(2)}} - \frac{\lambda_{2}}{\lambda_{3}} (1+c_{+})(1+c_{-}) \right) \left( 1 + \frac{\lambda_{1}\lambda_{2}}{\lambda_{3}^2} c_{-} (1+c_{-}) \right) } \\
&& \qquad \qquad + \frac{\lambda_{2}^{(1)}\lambda_{1}^{(2)}-\lambda_{1}^{(1)}\lambda_{2}^{(2)}}{\lambda_{3}^{(1)}\lambda_{2}^{(2)}} (1+c_{+}) \left( \frac{\lambda_{2}^2}{\lambda_{3}^2} (1+c_{-})^2 - \frac{\lambda_{1}}{\lambda_{3}} c_{-} \right) \\
\lefteqn{ \textrm{and } d_3 = 1 - \mu_2 \frac{\lambda_{2}^{(1)}\lambda_{1}^{(2)}-\lambda_{1}^{(1)}\lambda_{2}^{(2)}}{\lambda_{3}^{(1)}\lambda_{2}^{(2)}} \frac{c_{-}}{1+c_{+}} + \frac{\mu_2^2}{(1+c_{+})^3} \frac{\lambda_{1}^{(2)}}{\lambda_{2}^{(2)}}. }
\end{eqnarray*}
Then $v_{j\nu(i)} \nu(a_{ij}) = d_1 d_4^{-1} [kp_{1}q_{1}kp_{1}jr_{1}p_{1}jk] = d_1 d_4^{-1} a_{jk} v'$, where $v'$ is symmetric, and we obtain $C=1$. \\

\subsection{The graph $\mathcal{E}_5^{(12)}$ for the twisted orbifold $\mathcal{D}^{(12)}$ invariant}

\begin{figure}[tb]
\begin{center}
\includegraphics[width=70mm]{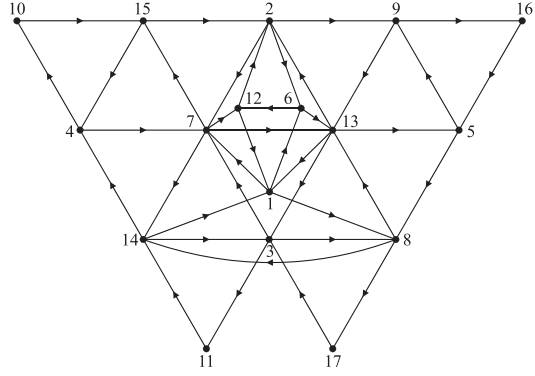}\\
 \caption{Labelled graph $\mathcal{E}_5^{(12)}$} \label{fig:E5(12)A}
\end{center}
\end{figure}

The Moore-Seiberg invariant $Z_{\mathcal{E}_{MS}^{(12)}}$ is realised by a braided subfactor which produces the nimrep $\mathcal{E}_5^{(12)}$ \cite[Section 5.4]{evans/pugh:2009ii}, illustrated in Figure \ref{fig:E5(12)A}. The unique cell system $W$ (up to equivalence) was computed in \cite[Theorem 13.1]{evans/pugh:2009i}.
For the graph $\mathcal{E}_5^{(12)}$ the automorphism $\nu$ is the identity. There are four possible choices for the vertices $j$, $k$: these are (3,8), (5,9), (14,3) and (15, 4). We see that it will not be possible to use the quicker method described above for the case where $\nu = \mathrm{id}$, since conjugating the graph will not interchange $j \leftrightarrow k$ for any of the choices of $j$, $k$.
We choose $i=10$, $j=15$ and $k=4$.
We will write out a basis for the space of paths which start from the vertices 10, 15, 4. Here $q=e^{2 \pi i/12}$, and we will write $[m]$ for the quantum integer $[m]_q$.

\noindent
Paths starting at vertex 10:
\scriptsize
\begin{eqnarray*}
l && \textrm{Paths} \\
\hline 2 && [10,15,2], \qquad ([10,15,4] = 0) \\
\hline 3 && [10,15,2,9], \qquad [10,15,2,6] \qquad \left( [10,\underline{15,2,7}] = -\frac{[2]}{[4]} [\underline{10,15,4},7] = 0 \right) \\
\hline 4 && [10,15,2,9,16], \qquad [10,15,\underline{2,9,13}] = \frac{[2]}{\sqrt{[3]}} [\underline{10,15,2,7},13] - \frac{\sqrt{[2]^3}}{\sqrt{[3][4]}} [10,15,2,6,13] = - \frac{\sqrt{[2]^3}}{\sqrt{[3][4]}} [10,15,2,6,13] \\
&& \left( [10,15,\underline{2,6,12}] = -\frac{\sqrt{[2]}}{\sqrt{[4]}} [\underline{10,15,2,7},12] = 0 \right) \\
\hline 5 && [10,15,2,\underline{9,16,5}] = - \frac{\sqrt{[4]}}{\sqrt{[2]}} [10,15,2,9,13,5], \qquad [10,15,2,9,13,3] \\
&& \left( [\underline{10,15,2,9,13},1] = -\frac{\sqrt{[2]^3}}{\sqrt{[3][4]}} [10,15,2,\underline{6,13,1}] = \frac{[2]}{[4]\sqrt{[3]}} [\underline{10,15,2,6,12},1] = 0 \right) \qquad \left( [10,15,2,\underline{9,13,2}] = 0 \right) \\
\hline 6 && [\underline{10,15,2,9,16,5},8] = - \frac{\sqrt{[4]}}{\sqrt{[2]}} [10,15,2,9,\underline{13,5,8}] = \frac{\sqrt{[4]}}{\sqrt{[2]}} [10,15,2,9,13,3,8], \qquad [10,15,2,9,13,3,11] \\
&& \left( [10,15,2,9,\underline{16,5,9}] = 0 \right) \qquad \left( [10,15,2,9,\underline{13,3,7}] = \frac{\sqrt{[3]}}{\sqrt{[2][6]}} [\underline{10,15,2,9,13,1},7] - \frac{\sqrt{[2]}}{\sqrt{[6]}} [\underline{10,15,2,9,13,2},7] = 0 \right) \\
\hline 7 && [\underline{10,15,2,9,16,5,8},14] = \frac{\sqrt{[4]}}{\sqrt{[2]}} [10,15,2,9,13,\underline{3,8,14}] = -\frac{\sqrt{[3][4]}}{\sqrt{[6]}} [10,15,2,9,13,3,11,14] \\
&& \left( [10,15,2,9,16,\underline{5,8,13}] = -\frac{\sqrt{[2]}}{\sqrt{[4]}} [\underline{10,15,2,9,16,5,9},13] = 0 \right) \\
&& \left( [\underline{10,15,2,9,16,5,8},17] = \frac{\sqrt{[4]}}{\sqrt{[2]}} [10,15,2,9,13,\underline{3,8,17}] = 0 \right) \\
\hline 8 && [10,15,2,9,16,5,8,14,4] \\
&& \left( [10,15,2,9,16,5,\underline{8,14,3}] = 0 \right) \qquad \left( [10,15,2,9,16,5,\underline{8,14,1}] = -\sqrt{[4]} [\underline{10,15,2,9,16,5,18,13},1] = 0 \right) \\
\hline 9 && [10,15,2,9,16,5,8,14,4,10] \qquad \left( [10,15,2,9,16,5,8,\underline{14,4,7}] = 0 \right) \\
\hline 10 && \left( [10,15,2,9,16,5,8,14,\underline{4,10,15}] = -\frac{\sqrt{[4]}}{\sqrt{[2]}} [\underline{10,15,2,9,16,5,8,14,4,7},15] = 0 \right)
\end{eqnarray*}
\normalsize

\noindent
Paths starting at vertex 15:
\scriptsize
\begin{eqnarray*}
l && \textrm{Paths} \\
\hline 2 && [15,2,6], \qquad [15,2,7] = -\frac{\sqrt{[2]}}{\sqrt{[4]}} [15,4,7], \qquad [15,2,9] \qquad ([10,15,4] = 0) \\
\hline 3 && [15,\underline{2,6,12}] = -\frac{\sqrt{[2]}}{\sqrt{[4]}} [15,2,7,12], \qquad [15,2,6,13], \qquad [15,2,7,13], \qquad [15,2,9,16] \\
&& \left( [\underline{15,2,7},14] = -\frac{\sqrt{[2]}}{\sqrt{[4]}} [15,\underline{4,7,14}] = 0 \right) \qquad \left( [15,\underline{2,7,15}] = 0 \right) \\
&& \left( [15,\underline{2,9,13}] = \frac{[2]}{\sqrt{[3]}} [15,2,7,13]^* - \frac{\sqrt{[2]^3}}{\sqrt{[3][4]}} [15,2,6,13]^* \right) \\
\hline 4 && [15,2,6,13,1] = -\frac{1}{\sqrt{[2][4]}} [15,2,6,12,1] = \frac{1}{[4]} [15,2,7,12,1] = -\frac{1}{\sqrt{[2]^3[4]}} [15,2,7,13,1], \\
&& [15,2,6,13,2] = -\frac{\sqrt{[4]}}{\sqrt{[2]}} [15,2,6,12,2] = [15,2,7,12,2] = -\frac{\sqrt{[4]}}{\sqrt{[2]}} [15,2,7,13,2], \\
&& [15,2,6,13,3], \qquad [15,2,6,13,5], \qquad [15,2,7,13,5] \qquad \left( [15,2,\underline{7,13,3}] = 0 \right) \\
&& \left( [15,2,\underline{9,16,5}] = -\frac{\sqrt{[4]}}{\sqrt{[2]}} [15,\underline{2,9,13},5] = \frac{[2]}{\sqrt{[3]}} [15,2,6,13,5]^* - \frac{\sqrt{[2][4]}}{\sqrt{[3]}} [15,2,7,13,5]^* \right) \\
\hline 5 && [15,2,\underline{6,13,2},7] = -\frac{\sqrt{[4]}}{\sqrt{[2]}} [15,2,6,\underline{12,2,7}] = \frac{\sqrt{[3][4]}}{\sqrt{[2]}} [15,2,\underline{6,12,1},7] = -[4]\sqrt{[3]} [15,2,6,\underline{13,1,7}] \\
&& \qquad = -[2][4] [15,2,6,13,2,7] - [4]\sqrt{[2][6]} [15,2,6,13,3,7], \quad \left( \Rightarrow [15,2,6,13,3,7] = - \frac{[3]^2}{[4]\sqrt{[2][6]}} [15,2,6,13,2,7]^* \right) \\
&& [15,2,6,12,1,8], \qquad [\underline{15,2,6,12,2},9] = -\frac{\sqrt{[2]}}{\sqrt{[4]}} [15,2,6,\underline{13,2,9}] = \frac{[2]}{[4]} [15,2,6,13,5,9] \\
&& [15,2,6,13,3,8], \qquad [15,2,6,13,3,11] \\
&& \left( [\underline{15,2,7,13,1},6] = [2] [15,2,6,\underline{12,1,6}] = -\frac{[2][4]}{\sqrt{[3]}} [15,2,\underline{6,12,2},6] = \frac{\sqrt{[2]^3[4]}}{\sqrt{[3]}} [\underline{15,2,6,13,2},6] \right. \\
&& \qquad = \frac{[2][4]}{\sqrt{[3]}} [15,2,7,\underline{13,2,6}] = -[2][4] [15,2,7,13,1,6] \qquad \Rightarrow [15,2,7,13,1,6] = 0 \Bigg) \\
&& \left( [15,2,6,\underline{13,5,8}] = -[15,2,6,13,3,8] - \frac{\sqrt{[6]}}{\sqrt{[4]}} [15,2,\underline{6,13,1},8] = -[15,2,6,13,3,8]^* + \frac{\sqrt{[6]}}{[4]\sqrt{[2]}} [15,2,6,12,1,8]^* \right) \\
&& \left( [15,2,7,\underline{13,5,8}] = -\frac{\sqrt{[6]}}{\sqrt{[4]}} [\underline{15,2,7,13,1},8] = -\frac{[2]\sqrt{[6]}}{\sqrt{[4]}} [15,2,6,12,1,8]^* \right) \\
\hline 6 && [15,2,6,12,\underline{1,7,13}] = -\frac{\sqrt{[6]}}{\sqrt{[4]}} [15,2,6,12,1,8,13], \qquad [15,2,6,12,\underline{1,7,14}] = -\sqrt{[4]} [15,2,6,12,1,8,14], \\
&& [15,2,6,12,1,8,17], \qquad [15,2,6,13,3,8,14] \qquad \left( [15,2,6,\underline{12,1,7},15] = -\frac{1}{\sqrt{[3]}} [15,2,6,12,\underline{2,7,15}] = 0 \right) \\
&& \left( [15,2,6,12,\underline{1,7,12}] = 0 \right) \qquad \left( [15,2,6,12,\underline{2,9,13}] = \frac{[2]}{\sqrt{[3]}} [15,2,6,\underline{12,2,7},13] = -[2] [15,2,6,12,1,7,13]^* \right) \\
&& \left( [\underline{15,2,6,12,2,9},16] = \frac{[2]}{[4]} [15,2,6,13,\underline{5,9,16}] = 0 \right) \qquad \left( [15,2,6,13,\underline{3,8,17}] = 0 \right) \\
&& \left( [15,2,6,13,\underline{3,8,13}] = \frac{\sqrt{[2][6]}}{\sqrt{[3]}} [\underline{15,2,6,13,3,7},13] = -\frac{[3]^2}{\sqrt{[2][4]}} [15,2,6,12,1,7,13]^* \right) \\
&& \left( [15,2,6,13,\underline{3,11,14}] = -[\underline{15,2,6,13,3,7},14] - \frac{\sqrt{[6]}}{\sqrt{[2][3]}} [15,2,6,13,3,8,14] \right. \\
&& \left. \qquad = \frac{\sqrt{[3]^5}}{\sqrt{[2]^3[4][6]}} [15,2,6,13,1,7,14]^* - \frac{\sqrt{[6]}}{\sqrt{[2][3]}} [15,2,6,13,3,8,14]^* \right) \\
\hline 7 && [15,2,6,12,1,\underline{7,13,1}] = -\frac{\sqrt{[6]}}{\sqrt{[4]}} [15,2,6,12,1,7,14,1] \qquad [15,2,6,12,1,\underline{7,13,3}] = \frac{\sqrt{[3]}}{\sqrt{[2][6]}} [15,2,6,12,1,7,14,3], \\
&& [15,2,6,13,3,8,14,4] \qquad \left( [15,2,6,12,1,\underline{7,13,2}] = 0 \right) \qquad \left( [15,2,6,12,1,\underline{7,14,4}] = 0 \right) \\
&& \left( [\underline{15,2,6,12,1,7,13},5] = -\frac{\sqrt{[6]}}{\sqrt{[4]}} [15,2,6,12,1,\underline{8,13,5}] = 0 \right) \\
&& \left( [15,2,6,12,1,\underline{8,17,3}] = -[\underline{15,2,6,12,1,8,13},3] - \frac{\sqrt{[6]}}{\sqrt{[2][3]}} [\underline{15,2,6,12,1,8,14},3] \right. \\
&& \left. \quad = \frac{\sqrt{[4]}}{\sqrt{[6]}} [15,2,6,12,1,7,13,3] + \frac{\sqrt{[4][6]}}{\sqrt{[2][3]}} [\underline{15,2,6,12,1,7,14,3}] = \frac{\sqrt{[4]}}{[3]\sqrt{[6]}} \left( [3]+\sqrt{[6]^3} \right) [15,2,6,12,1,7,13,3]^* \right) \\
&& \left( [15,2,6,13,3,\underline{8,14,3}] = -\frac{\sqrt{[2][3]}}{\sqrt{[6]}} [\underline{15,2,6,13,3,8,13},3] = \frac{\sqrt{[3]^5}}{\sqrt{[4][6]}} [15,2,6,12,1,7,13,3]^* \right) \\
\hline 8 && [15,2,6,12,1,7,\underline{13,1,7}] = \frac{\sqrt{[2][6]}}{\sqrt{[3]}} [15,2,6,12,1,7,13,3,7], \qquad [15,2,6,13,3,8,14,4,10] \\
&& \left( [15,2,6,12,1,7,\underline{13,1,6}] = 0 \right) \qquad \left( [\underline{15,2,6,12,1,7,13,3},11] = \frac{\sqrt{[3]}}{\sqrt{[2][6]}} [15,2,6,12,1,7,\underline{14,3,11}] = 0 \right) \\
&& \left( [15,2,6,12,1,7,\underline{13,1,8}] = -\frac{\sqrt{[4]}}{\sqrt{[6]}} [\underline{15,2,6,12,1,7,13,3},8] = -\frac{\sqrt{[3][4]}}{[6]\sqrt{[2]}} [15,2,6,12,1,7,\underline{14,3,8}] \right. \\
&& \left. \qquad = \frac{[3]\sqrt{[4]}}{[6]} [\underline{15,2,6,12,1,7,14,1},8] = -\frac{[3][4]}{\sqrt{[6]^3}} [15,2,6,12,1,7,13,1,8] \qquad \Rightarrow [15,2,6,12,1,7,13,1,8] = 0 \right) \\
&& \left( [15,2,6,13,3,8,\underline{14,4,7}] = -\frac{\sqrt{[6]}}{\sqrt{[4]}} [\underline{15,2,6,13,3,8,14,1},7] - [\underline{15,2,6,13,3,8,14,3},7] \right. \\
&& \qquad = -\frac{[3]^2\sqrt{[6]}}{\sqrt{[2][4]^3}} [15,2,6,12,1,7,13,1,7] - \frac{\sqrt{[3]^5}}{\sqrt{[4][6]}} [\underline{15,2,6,12,1,7,13,3,7}] \\
&& \left. \qquad = -\frac{[3]^2}{[6]\sqrt{[2][4]^3}} \left( \sqrt{[6]^3}+[3][4] \right) [15,2,6,12,1,7,13,1,7]^* \right) \\
\hline 9 && [15,2,6,12,1,7,13,1,7,15] \\
&& \left( [15,2,6,12,1,7,13,\underline{1,7,12}] = 0 \right) \qquad \left( [15,2,6,12,1,7,13,\underline{1,7,13}] = 0 \right) \qquad \left( [15,2,6,12,1,7,13,\underline{1,7,14}] = 0 \right) \\
&& \left( [15,2,6,13,3,8,14,\underline{4,10,15}] = -\frac{\sqrt{[4]}}{\sqrt{[2]}} [\underline{15,2,6,13,3,8,14,4,7},15] \right. \\
&& \left. \qquad = \frac{[3]^2}{[2][4][6]} \left( \sqrt{[6]^3}+[3][4] \right) [15,2,6,12,1,7,13,1,7,15]^* \right) \\
\hline 10 && \left( [15,2,6,12,1,7,13,1,\underline{7,15,2}] = 0 \right) \qquad \left( [15,2,6,12,1,7,13,1,\underline{7,15,4}] = 0 \right)
\end{eqnarray*}
\normalsize

\noindent
Paths starting at vertex 4:
\scriptsize
\begin{eqnarray*}
l && \textrm{Paths} \\
\hline 2 && [4,7,12], \qquad [4,7,13], \qquad [4,7,15] = -\frac{\sqrt{[2]}}{\sqrt{[4]}} [4,10,15] \qquad \left( [4,7,14] = 0 \right) \\
\hline 3 && [4,\underline{7,12,1}] = -\frac{\sqrt{[4]}}{\sqrt{[2]^3}} [4,7,13,1], \qquad [4,7,12,2], \qquad [4,7,13,2], \qquad [4,7,13,5] \\
&& \left( [4,\underline{7,13,3}] = 0 \right) \qquad \left( [4,\underline{7,15,2}] = \frac{\sqrt{[2]^3}}{\sqrt{[3][4]}} [4,7,12,2]^* + \frac{[2]}{\sqrt{[3]}} [4,7,13,2]^* \right) \qquad \left( [4,\underline{7,15,}4] = 0 \right) \\
\hline 4 && [4,7,\underline{13,2,6}] = -\sqrt{[3]} [\underline{4,7,13,1},6] = \frac{\sqrt{[2]^3[3]}}{\sqrt{[4]}} [4,7,\underline{12,1,6}] = -\sqrt{[2]^3} [4,\underline{7,12,2},6] \\
&& \quad = -[2][4] [4,7,13,2,6] + \sqrt{[2][3][4]} [4,7,15,2,6] \quad \left( \Rightarrow [4,7,15,2,6] = \frac{\sqrt{[3]}}{\sqrt{[2][4]}} [4,7,13,2,6] = \frac{[2][3]}{[4]} [4,7,12,1,6]^* \right) \\
&& [4,7,\underline{13,2,7}] = \frac{\sqrt{[3]}}{[2]} [\underline{4,7,13,1},7] = -\frac{\sqrt{[2][3]}}{\sqrt{[4]}} [4,7,\underline{12,1,7}] = \frac{\sqrt{[2]}}{\sqrt{[4]}} [4,\underline{7,12,2},7] \\
&& \qquad = [4,7,13,2,7] - \frac{\sqrt{[3]}}{\sqrt{[2][4]}} [4,7,15,2,7] \qquad \left( \Rightarrow [4,7,15,2,7] = 0 \right) \\
&& [4,7,12,1,8], \qquad [4,7,12,2,9], \qquad [4,7,\underline{13,2,9}] = -\frac{\sqrt{[2]}}{\sqrt{[4]}} [4,7,13,5,9] \\
&& \left( [4,7,\underline{13,2,7}] = \frac{\sqrt{[3]}}{[2]} [4,\underline{7,13,1},7] = -\frac{\sqrt{[2][3]}}{\sqrt{[4]}} [4,7,12,1,7]^* \right) \\
&& \left( [4,7,\underline{13,5,8}] = -\frac{\sqrt{[6]}}{\sqrt{[4]}} [\underline{4,7,13,1},8] = \sqrt{[2]^3[6]} [4,7,12,1,8]^* \right) \\
\hline 5 && [4,7,12,\underline{1,6,12}] = -\sqrt{[2][4]} [4,7,12,1,7,12], \qquad [4,7,12,1,6,13], \qquad [4,7,12,1,7,13], \\
&& [4,7,12,\underline{1,7,14}] = -\sqrt{[4]} [4,7,12,1,8,14], \qquad [4,7,12,1,8,17], \qquad [4,7,12,2,9,16] \\
&& \left( [\underline{4,7,12,1,7},15] = -\frac{1}{\sqrt{[3]}} [4,7,12,\underline{2,7,15}] = 0 \right) \qquad \left( [\underline{4,7,13,2,9},16] = -\frac{\sqrt{[2]}}{\sqrt{[4]}} [4,7,13,\underline{5,9,16}] = 0 \right) \\
&& \left( [4,7,12,\underline{1,8,13}] = -\frac{\sqrt{[2]^3}}{\sqrt{[6]}} [4,7,12,1,6,13]^* - \frac{\sqrt{[4]}}{\sqrt{[6]}} [4,7,12,1,7,13]^* \right) \\
&& \left( [4,7,12,\underline{2,9,13}] = -\frac{\sqrt{[2]^3}}{\sqrt{[3][4]}} [\underline{4,7,12,2,6},13] + \frac{[2]}{\sqrt{[3]}} [\underline{4,7,12,2,7},13] \right. \\
&& \left. \qquad = \frac{\sqrt{[2]^3}}{\sqrt{[4]^3}} [4,7,12,1,6,13]^* - [2] [4,7,12,1,7,13]^* \right) \\
&& \left( [4,7,13,\underline{2,9,13}] = -\frac{\sqrt{[2]^3}}{\sqrt{[3][4]}} [\underline{4,7,13,2,6},13] + \frac{[2]}{\sqrt{[3]}} [\underline{4,7,13,2,7},13] \right. \\
&& \left. \qquad = -\frac{[2]^3}{[4]} [4,7,12,1,6,13]^* - \frac{\sqrt{[2]^3}}{\sqrt{[4]}} [4,7,12,1,7,13]^* \right) \\
\hline 6 && [4,7,12,1,7,13,1], \qquad [4,7,12,1,6,13,3], \qquad [4,7,12,1,6,13,5], \qquad [4,7,12,1,7,13,3] \\
&& \left( \textrm{We know from the Hilbert series of Theorem \ref{thm:SU(3)Hilbert} that all paths of length 6 which start at 4 and end at 2 are zero.} \right) \\
&& \left( [\underline{4,7,12,1,6},12,1] = \frac{[4]}{[2][3]} [4,7,\underline{15,2,6,12,1}] = -\frac{\sqrt{[4]}}{[3]\sqrt{[2]}} [4,7,\underline{15,2,7},12,1] = \frac{1}{[3]} [\underline{4,7,15,4},7,12,1] = 0 \right) \\
&& \left( [4,7,12,1,\underline{7,14,1}] = -\frac{\sqrt{[2]^3}}{\sqrt{[6]}} [4,7,12,\underline{1,7,12},1] - \frac{\sqrt{[4]}}{\sqrt{[6]}} [4,7,12,1,7,13,1] \right. \\
&& \left. \qquad = \frac{[2]}{\sqrt{[4][6]}} [\underline{4,7,12,1,6,12,1}] - \frac{\sqrt{[4]}}{\sqrt{[6]}} [4,7,12,1,7,13,1] = -\frac{\sqrt{[4]}}{\sqrt{[6]}} [4,7,12,1,7,13,1] \right) \\
&& \left( [4,7,12,\underline{1,7,13},5] = -\frac{\sqrt{[2]^3}}{\sqrt{[4]}} [4,7,12,1,6,13,5]^* \right) \qquad \left( [4,7,12,1,\underline{7,14,3}] = \frac{\sqrt{[2][6]}}{\sqrt{[3]}} [4,7,12,1,7,13,3]^* \right) \\
&& \left( [4,7,12,1,\underline{8,17,3}] = -[4,7,12,\underline{1,8,13},3] - \frac{\sqrt{[6]}}{\sqrt{[2][3]}} [4,7,12,\underline{1,8,14},3] \right. \\
&& \qquad = \frac{\sqrt{[2]^3}}{\sqrt{[6]}} [4,7,12,1,6,13,3] + \frac{\sqrt{[4]}}{\sqrt{[6]}} [4,7,12,1,7,13,3] + \frac{\sqrt{[6]}}{\sqrt{[2][3][4]}} [\underline{4,7,12,1,7,14,3}] \\
&& \left. \qquad = \frac{\sqrt{[2]^3}}{\sqrt{[6]}} [4,7,12,1,6,13,3]^* + \left( \frac{\sqrt{[4]}}{\sqrt{[6]}} + \frac{[6]}{[3]\sqrt{[4]}} \right) [4,7,12,1,7,13,3]^* \right) \\
&& \left( [4,7,12,2,\underline{9,16,5}] = -\frac{\sqrt{[4]}}{\sqrt{[2]}} [\underline{4,7,12,2,9,13},5] = -\frac{[2]}{[4]} [4,7,12,1,6,13,5] + \sqrt{[2][4]} [\underline{4,7,12,1,7,13,5}] \right. \\
&& \left. \qquad = -\frac{[2][3]^2}{[4]} [4,7,12,1,6,13,5]^* \right) \qquad \left( [4,7,12,1,\underline{7,14,4}] = 0 \right) \\
\hline 7 && \textrm{We will let } \xi \textrm{ denote the path } [4,7,12,1] \textrm{ of length 3} \\
&& [\xi,7,13,1,7], \qquad [\xi,7,13,1,8], \qquad [\xi,6,13,3,11] \qquad \left( [\xi,6,\underline{13,3,7}] = \frac{\sqrt{[3]}}{\sqrt{[2][6]}} [\xi,\underline{6,13,1},7] = 0 \right) \\
&& \left( [\xi,7,\underline{13,1,6}] = 0 \right) \qquad \left( [\xi,\underline{7,13,3},8] = \frac{\sqrt{[3]}}{\sqrt{[2][6]}} [\xi,7,\underline{14,3,8}] = -\frac{[3]}{\sqrt{[6]}} [\underline{\xi,7,14,1},8] = \frac{[3]\sqrt{[4]}}{[6]} [\xi,7,13,1,8]^* \right) \\
&& \left( [\xi,6,\underline{13,3,8}] = -[4,7,12,\underline{1,6,13},5,8] = \frac{\sqrt{[4]}}{\sqrt{[2]^3}} [\xi,7,\underline{13,5,8}] = -\frac{\sqrt{[6]}}{\sqrt{[2]^3}} [\xi,7,13,1,8] - \frac{\sqrt{[4]}}{\sqrt{[2]^3}} [\underline{\xi,7,13,3,8}] \right. \\
&& \left. \qquad = -\frac{1}{[6]\sqrt{[2]^3}} \left( \sqrt{[6]^3}+[3][4] \right) [\xi,7,13,1,8]^* \right) \qquad \left( [\xi,6,\underline{13,5,9}] = 0 \right) \\
&& \left( [\xi,7,\underline{13,3,7}] = -\frac{\sqrt{[3]}}{\sqrt{[2][6]}} [\xi,7,13,1,7]^* \right) \qquad \left( [\xi,\underline{7,13,3},11] = \frac{\sqrt{[3]}}{\sqrt{[2][6]}} [\xi,7,\underline{14,3,11}] = 0 \right) \\
\hline 8 && [\xi,7,13,\underline{1,7,14}] = -\sqrt{[4]} [\xi,7,13,1,8,14], \qquad [\xi,7,13,1,7,15] \\
\lefteqn{ \left( \textrm{We know from the Hilbert series of Theorem \ref{thm:SU(3)Hilbert} that all paths of length 8 which start at 4 and end at 12, 13, 17 are zero.} \right) } \\
&& \left( [\xi,6,13,\underline{3,11,14}] = -[\xi,6,\underline{13,3,7},14] - \frac{\sqrt{[6]}}{\sqrt{[2][3]}} [\xi,6,\underline{13,3,8},14] = \frac{\sqrt{[6]}}{\sqrt{[2][3]}} [\underline{\xi,6,13,5,8},14] \right. \\
&& \left. \qquad = \frac{1}{[2]^2\sqrt{[3][6]}} \left( \sqrt{[6]^3}+[3][4] \right) [\xi,7,13,1,8,14]^* \right) \\
\hline 9 && [\xi,7,13,1,\underline{7,14,4}] = -\frac{\sqrt{[2]}}{\sqrt{[4]}} [\xi,7,13,1,7,15,4] \\
&& \left( [\xi,7,13,1,\underline{7,14,1}] = 0 \right) \qquad \left( [\xi,7,13,1,\underline{7,14,3}] = 0 \right) \qquad \left( [\xi,7,13,1,\underline{7,15,2}] = 0 \right) \\
\hline 10 && \left( [\xi,7,13,1,\underline{7,14,4},10] = -\frac{\sqrt{[2]}}{\sqrt{[4]}} [\xi,7,13,1,7,\underline{15,4,10}] = 0 \right) \qquad \left( [\xi,7,13,1,7,\underline{14,4,7}] = 0 \right)
\end{eqnarray*}
\normalsize

We choose $u_{i\nu(i)} = u_{10\nu(10)} = [10,15,2,9,16,5,8,14,4,10]$. Then
\begin{eqnarray*}
\lefteqn{ v_{15\nu(10)} \nu(a_{10,15}) \;\; = \;\; [\underline{15,2,9,16,5},8,14,4,10,15] } \\
& = & -\frac{\sqrt{[2][4]}}{\sqrt{[3]}} [\underline{15,2,7,13,5,8},14,4,10,15] + \frac{[2]}{\sqrt{[3]}} [\underline{15,2,6,13,5,8},14,4,10,15] \\
& = & \frac{\sqrt{[2]^3[6]}}{\sqrt{[3]}} [15,2,6,12,1,8,14,4,10,15] - \frac{[2]}{\sqrt{[3]}} [\underline{15,2,6,13,3,8,14,4,10,15}] \\
&& \qquad + \frac{\sqrt{[2][6]}}{[4]\sqrt{[3]}} [15,2,6,12,1,8,14,4,10,15] \\
& = & \frac{\sqrt{[2][6]}}{[4]\sqrt{[3]}} \left( [2]+[4] \right) [15,2,6,12,\underline{1,8,14},4,10,15] \\
&& \qquad -\frac{\sqrt{[3]^3}}{[4][6]} \left( \sqrt{[6]^3}+[3][4] \right) [15,\underline{2,6,12},1,7,13,1,7,15] \\
& = & -\frac{\sqrt{[2]^3[3][6]}}{\sqrt{[4]^3}} [\underline{15,2,6,12,1,7,14,4},10,15] \\
&& \qquad + \frac{\sqrt{[2][3]^3}}{[6]\sqrt{[4]^3}} \left( \sqrt{[6]^3}+[3][4] \right) [\underline{15,2,7},12,1,7,13,1,7,15] \\
& = & -\frac{[2]\sqrt{[3]^3}}{[4]^2[6]} \left( \sqrt{[6]^3}+[3][4] \right) [15,4,7,12,1,7,13,1,7,15] \;\; = \;\; c a_{15,4} v_{4\nu(15)},
\end{eqnarray*}
and
\begin{eqnarray*}
\lefteqn{ a_{4,10} v_{10\nu(4)} \;\; = \;\; [\underline{4,10,15},2,9,16,5,8,14,4] \;\; = \;\; -\frac{\sqrt{[4]}}{\sqrt{[2]}} [\underline{4,7,15,2},9,16,5,8,14,4] } \\
& = & \frac{[2]}{\sqrt{[3]}} [\underline{4,7,12,2,9,16,5},8,14,4] - \frac{\sqrt{[2][4]}}{\sqrt{[3]}} [\underline{4,7,13,2,9},16,5,8,14,4] \\
& = & -\frac{[2]^2\sqrt{[3]^3}}{[4]} [\underline{4,7,12,1,6,13,5,8},14,4] + \frac{[2]}{\sqrt{[3]}} [4,7,13,\underline{5,9,16},5,8,14,4] \\
& = & -\frac{\sqrt{[2][3]^3}}{[4][6]} \left( \sqrt{[6]^3}+[3][4] \right) [4,7,12,1,7,13,\underline{1,8,14},4] \\
& = & \frac{\sqrt{[2][3]^3}}{[6]\sqrt{[4]^3}} \left( \sqrt{[6]^3}+[3][4] \right) [4,7,12,1,7,13,1,\underline{7,14,4}] \\
& = & -\frac{[2]\sqrt{[3]^3}}{[4]^2[6]} \left( \sqrt{[6]^3}+[3][4] \right) [4,7,12,1,7,13,1,7,15,4] \;\; = \;\; c v_{4\nu(15)} \nu(a_{15,4}).
\end{eqnarray*}
Then $C=1$. \\

\subsection{$\mathcal{E}^{(24)}$ graph for the conformal embedding $SU(3)_{21} \subset (E_7)_1$}

\begin{figure}[tb]
\begin{center}
\includegraphics[width=80mm]{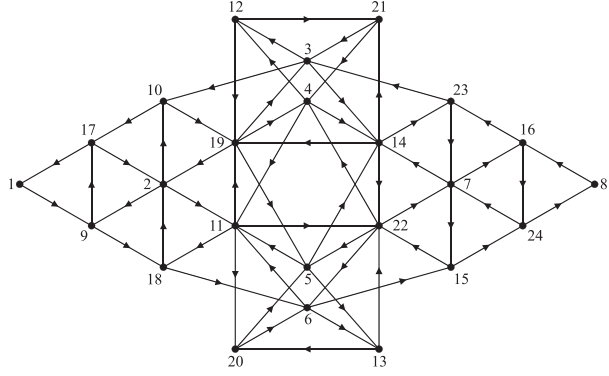}\\
 \caption{Labelled graph $\mathcal{E}^{(24)}$} \label{fig:E(24)A}
\end{center}
\end{figure}

For the graph $\mathcal{E}^{(24)}$, illustrated in Figure \ref{fig:E(24)A}, the automorphism $\nu$ is the identity.
The unique cell system $W$ (up to equivalence) was computed in \cite[Theorem 14.1]{evans/pugh:2009i}.
We choose $i=1$, $j=9$ and $k=17$.
We will write out a basis for the space of paths which start from the vertices 1, 9. \\

\noindent
Paths starting at vertex 1:
\scriptsize
\begin{eqnarray*}
l && \textrm{Paths} \\
\hline 2 && [1,9,18] \qquad ([1,9,17] = 0) \\
\hline 3 && [1,9,18,6] =: \xi_1 \qquad \left( [1,\underline{9,18,2}] = 0 \right) \\
\hline 4 && [1,9,18,6,13], \qquad [1,9,18,6,15] \qquad \left( [1,9,\underline{18,6,11}] = 0 \right) \\
\hline 5 && [1,9,18,\underline{6,13,22}] = -\frac{[4]}{\sqrt{[9]}} [1,9,18,6,15,22], \qquad [1,9,18,6,15,24] \qquad \left( [1,9,18,\underline{6,13,20}] = 0 \right) \\
\hline 6 && [\xi_1,13,22,4], \qquad [\underline{\xi_1,13,22},7] = -\frac{[4]}{\sqrt{[9]}} [\xi_1,\underline{15,22,7}] = \frac{[4]\sqrt{[3]}}{\sqrt{[5][9]}} [\xi_1,15,24,7], \qquad [\xi_1,15,24,8] \\
&& \left( [\xi_1,\underline{13,22,5}] = 0 \right) \qquad \left( [\xi_1,\underline{13,22,6}] = 0 \right) \\
\hline 7 && [\xi_1,13,22,4,12], \qquad [\xi_1,13,\underline{22,4,14}] = -\frac{\sqrt{[7]}}{[3]} [\xi_1,13,22,7,14], \\
&& [\underline{\xi_1,13,22,7},16] = \frac{[4]\sqrt{[3]}}{\sqrt{[5][9]}} [\xi_1,15,\underline{24,7,16}] = -\frac{\sqrt{[2][3][4]}}{\sqrt{[5][9]}} [\xi_1,15,24,8,16] \\
&& \left( [\xi_1,13,\underline{22,4,11}] = 0 \right) \qquad \left( [\underline{\xi_1,13,22,7},15] = \frac{[4]\sqrt{[3]}}{\sqrt{[5][9]}} [\xi_1,15,\underline{24,7,15}] = 0 \right) \\
\hline 8 && [\xi_1,13,22,\underline{4,12,19}] = -\frac{\sqrt{[5][7]}}{[3]\sqrt{[9]}} [\xi_1,13,22,4,14,19], \qquad [\xi_1,13,22,\underline{4,12,21}] = -\frac{[3]\sqrt{[5]}}{\sqrt{[7]}} [\xi_1,13,22,4,14,21], \\
&& [\underline{\xi_1,13,22,4,14},23] = -\frac{\sqrt{[7]}}{[3]} [\xi_1,13,22,\underline{7,14,23}] = \frac{\sqrt{[7]}}{\sqrt{[3][5]}} [\xi_1,13,22,7,16,23] \\
&& \left( [\xi_1,13,22,\underline{4,14,22}] = 0 \right) \qquad \left( [\xi_1,13,22,\underline{7,16,24}] = 0 \right) \\
\hline 9 && [\xi_1,13,22,4,\underline{12,19,4}] = -\frac{\sqrt{[7]}}{[3]\sqrt{[5]}} [\xi_1,13,22,4,12,21,4], \qquad [\xi_1,13,22,4,\underline{12,19,3}] = \frac{\sqrt{[5]}}{\sqrt{[3]}} [\xi_1,13,22,4,12,21,3], \\
&& [\xi_1,13,22,4,12,19,2] \qquad \left( [\underline{\xi_1,13,22,4,12,19},5] = \frac{[7]\sqrt{[5]}}{[3]^2\sqrt{[9]}} [\xi_1,13,22,7,\underline{14,19,5}] = 0 \right) \\
&& \left( [\xi_1,13,22,4,\underline{14,23,3}] = -\frac{\sqrt{[3][5]}}{[4]} [\underline{\xi_1,13,22,4,14,19},3] - \frac{\sqrt{[9]}}{[4]} [\underline{\xi_1,13,22,4,14},21,3] \right. \\
&& \qquad = \frac{\sqrt{[3]^3[9]}}{[4]\sqrt{[7]}} [\xi_1,13,22,4,12,19,3] + \frac{\sqrt{[7][9]}}{[3][4]} [\underline{\xi_1,13,22,7,14,21,3}] \\
&& \left. \qquad = \frac{\sqrt{[9]}}{[4][5]\sqrt{[3][7]}} \left( [3]^2[5]+[7] \right) [\xi_1,13,22,4,12,19,3]^* \right) \qquad \left( [\xi_1,13,22,4,\underline{14,23,7}] = 0 \right) \\
\hline 10 && \textrm{We will let } \xi_2 \textrm{ denote the path } [\xi_1,13,22,4,12] = [1,9,18,6,13,22,4,12] \textrm{ of length 7.} \\
&& [\xi_2,19,2,9], \qquad [\xi_2,\underline{19,2,10}] = -\frac{\sqrt{[4]}}{\sqrt{[2]}} [\xi_2,19,3,10], \qquad [\xi_2,\underline{19,2,11}] = -\frac{\sqrt{[7]}}{[3]} [\xi_2,19,4,11], \\
&& [\xi_2,\underline{19,3,12}] = -\sqrt{[3]} [\xi_2,19,4,12] \\
&& \left( [\xi_2,\underline{19,3,14}] = -\frac{\sqrt{[7]}}{[3]} [\underline{\xi_2,19,4},14] = \frac{[7]}{[3]^2\sqrt{[5]}} [\xi_2,\underline{21,4,14}] = -\frac{[7]}{\sqrt{[3]^5[5]}} [\underline{\xi_2,21,3},14] = -\frac{[7]}{[3]^2[5]} [\xi_2,19,3,14] \right. \\
&& \qquad \Rightarrow [\xi_2,19,3,14] = 0 \qquad \Rightarrow [\xi_2,19,4,14] = 0 \bigg) \\
\hline 11 && [\xi_2,19,\underline{2,9,17}] = -\sqrt{[3]} [\xi_2,19,2,10,17], \qquad [\xi_2,19,\underline{2,9,18}] = -\frac{\sqrt{[5]}}{\sqrt{[3]}} [\xi_2,19,2,11,18], \qquad [\xi_2,19,2,11,20], \\
&& [\xi_2,19,2,10,19] = -\sqrt{[3]} [\xi_2,19,2,11,19] = \frac{\sqrt{[7]}}{\sqrt{[3]}} [\xi_2,19,4,11,19] = -\frac{\sqrt{[3][9]}}{\sqrt{[2][4]}} [\xi_2,19,4,12,19] \\
&& \qquad = \frac{\sqrt{[9]}}{\sqrt{[2][4]}} [\xi_2,19,3,12,19] \qquad \left( [\underline{\xi_2,19,2,11},22] = -\frac{\sqrt{[7]}}{[3]} [\xi_2,19,\underline{4,11,22}] = 0 \right) \qquad \left( [\xi_2,19,\underline{3,12,21}] = 0 \right) \\
\hline 12 && [\xi_2,19,2,\underline{9,17,2}] = -\sqrt{[3]} [\underline{\xi_2,19,2,9,18},2] = \sqrt{[5]} [\xi_2,19,2,\underline{11,18,2}] = -\sqrt{[3][5]} [\xi_2,19,\underline{2,11,19},2] \\
&& \qquad = \sqrt{[5]} [\xi_2,19,2,10,19,2], \\
&& [\xi_2,19,2,\underline{11,18,6}] = -\frac{\sqrt{[9]}}{[4]} [\xi_2,19,2,11,20,6], \qquad [\xi_2,19,2,\underline{11,19,5}] = -\frac{[3]\sqrt{[9]}}{\sqrt{[2][4][7]}} [\xi_2,19,2,11,20,5] \\
&& \left( [\xi_2,19,2,\underline{9,17,1}] = 0 \right) \qquad \left( [\xi_2,19,2,\underline{10,19,3}] = 0 \right) \qquad \left( [\xi_2,19,2,\underline{11,19,4}] = 0 \right) \\
\hline 13 && [\xi_2,19,2,11,\underline{18,6,11}] = -\frac{\sqrt{[2]}}{\sqrt{[4]}} [\underline{\xi_2,19,2,11,18,2},11] = \frac{\sqrt{[2][3]}}{\sqrt{[4]}} [\xi_2,19,2,11,\underline{19,2,11}] \\
&& \qquad = - \frac{\sqrt{[2][7]}}{\sqrt{[3][4]}} [\xi_2,19,2,11,19,5,11], \\
&& [\xi_2,19,2,11,18,6,13] = -\frac{\sqrt{[9]}}{[4]} [\xi_2,19,2,11,20,6,13] = -\frac{\sqrt{[7][9]}}{[4][5]} [\xi_2,19,2,11,20,5,13] \\
&& \qquad = \frac{[7]\sqrt{[2]}}{[3][5]\sqrt{[4]}} [\xi_2,19,2,11,19,5,13], \qquad [\xi_2,19,2,11,18,6,15] \\
&& \left( [\xi_2,19,2,11,\underline{18,2,9}] = 0 \right) \qquad \left( [\xi_2,19,2,11,\underline{19,2,10}] = 0 \right) \qquad \left( [\xi_2,19,2,11,\underline{19,5,14}] = 0 \right) \\
\hline 14 && \textrm{We will let } \xi_3 \textrm{ denote the path } [\xi_2,19,2,11,18] = [1,9,18,6,13,22,4,12,19,2,11,18] \textrm{ of length 11.} \\
&& [\xi_3,6,11,20] = \frac{\sqrt{[5]}}{\sqrt{[3]}} [\xi_3,6,13,20], \qquad [\xi_3,6,15,24] \\
&& [\underline{\xi_3,6,11},22] = -\frac{\sqrt{[2][7]}}{\sqrt{[3][4]}} [\xi_2,19,2,11,19,\underline{5,11,22}] = \frac{\sqrt{[2][3][9]}}{\sqrt{[4][5]}} [\underline{\xi_2,19,2,11,19,5,13},22] = \frac{\sqrt{[3]^2[5][9]}}{[7]} [\xi_3,6,13,22] \\
&& \left( [\xi_3,\underline{6,11,18}] = 0 \right) \qquad \left( [\xi_3,\underline{6,15,22}] = -\frac{\sqrt{[3][5]}}{[4]} [\xi_3,6,11,22] - \frac{\sqrt{[9]}}{[4]} [\underline{\xi_3,6,13,22}] = - \frac{[2][7]}{\sqrt{[3]^3[5]}} [\xi_3,6,11,22]^* \right) \\
\hline 15 && [\xi_3,6,\underline{11,20,5}] = -\frac{\sqrt{[5][7]}}{[3]\sqrt{[9]}} [\xi_3,6,11,22,5], \qquad [\xi_3,6,11,22,7], \qquad [\xi_3,6,15,24,8] \\
&& \left( [\xi_3,6,\underline{11,20,6}] = -\frac{\sqrt{[3][5]}}{\sqrt{[9]}} [\underline{\xi_3,6,11,22},6] = -\frac{[3]^2[5]}{[7]} [\xi_3,6,\underline{13,22,6}] = -\frac{\sqrt{[3]^3[5]^3}}{[7]} [\underline{\xi_3,6,13,20},6] \right. \\
&& \left. \qquad = -\frac{[3]^2[5]}{[7]} [\xi_3,6,11,20,6] \qquad \Rightarrow [\xi_3,6,11,20,6] = 0 \qquad \Rightarrow [\xi_3,6,11,22,6] = 0 \right) \qquad \left( [\xi_3,6,\underline{11,22,4}] = 0 \right) \\
&& \left( [\xi_3,6,\underline{15,24,7}] = -\frac{\sqrt{[5]}}{\sqrt{[3]}} [\underline{\xi_3,6,15,22},7] = \frac{[2][7]}{[3]^2} [\xi_3,6,11,22,7]^* \right) \\
\hline 16 && [\xi_3,6,11,\underline{22,5,14}] = -\frac{\sqrt{[7]}}{[3]} [\xi_3,6,11,22,7,14], \qquad [\xi_3,6,11,22,7,16] \\
&& \left( [\xi_3,6,11,\underline{20,5,11}] = 0 \right) \qquad \left( [\xi_3,6,11,\underline{20,5,13}] = 0 \right) \qquad \left( [\xi_3,6,11,\underline{22,7,15}] = 0 \right) \\
&& \left( [\xi_3,6,15,24,8,16] = -\frac{\sqrt{[4]}}{\sqrt{[2]}} [\xi_3,6,15,24,7,16] = \frac{[7]\sqrt{[2][4]}}{[3]^2} [\xi_3,6,11,22,7,16]^* \right) \\
\hline 17 && [\xi_3,6,11,22,5,14,21], \qquad [\underline{\xi_3,6,11,22,5,14},23] = -\frac{\sqrt{[7]}}{[3]} [\xi_3,6,11,22,\underline{7,14,23}] = \frac{\sqrt{[7]}}{\sqrt{[3][5]}} [\xi_3,6,11,22,7,16,23] \\
&& \left( [\xi_3,6,11,22,\underline{5,14,19}] = 0 \right) \qquad \left( [\underline{\xi_3,6,11,22,5,14},22] = -\frac{\sqrt{[7]}}{[3]} [\xi_3,6,11,22,\underline{7,14,22}] = 0 \right) \\
&& \left( [\xi_3,6,11,22,\underline{7,16,24}] = 0 \right) \\
\hline 18 && [\xi_3,6,11,22,5,\underline{14,21,3}] = -\frac{[4]}{\sqrt{[9]}} [\xi_3,6,11,22,5,14,23,3] \qquad \left( [\xi_3,6,11,22,5,\underline{14,21,4}] = 0 \right) \\
&& \left( [\xi_3,6,11,22,5,\underline{14,23,7}] = 0 \right) \\
\hline 19 && [\xi_3,6,11,22,5,14,21,3,10] \qquad \left( [\xi_3,6,11,22,5,14,\underline{21,3,12}] = 0 \right) \qquad \left( [\xi_3,6,11,22,5,14,\underline{21,3,14}] = 0 \right) \\
\hline 20 && [\xi_3,6,11,22,5,14,21,3,10,17] \qquad \left( [\xi_3,6,11,22,5,14,21,\underline{3,10,19}] = 0 \right) \\
\hline 21 && [\xi_3,6,11,22,5,14,21,3,10,17,1] \qquad \left( [\xi_3,6,11,22,5,14,21,3,\underline{10,17,2}] = 0 \right) \\
\hline 22 && \left( [\xi_3,6,11,22,5,14,21,3,10,\underline{17,1,9}] = 0 \right)
\end{eqnarray*}
\normalsize

\noindent
Paths starting at vertex 9:
\scriptsize

\normalsize

Let $u$ be the path $[1,9,18,6,13,22,4,14,19,2,11,19,2,11,22,5,14,21,3,10,17,1]$, which is symmetric. This path is non-zero in $A$ since
\begin{eqnarray*}
\lefteqn{ [\underline{1,9,18,6,13,22,4,14,19,2,11,19,2},11,22,5,14,21,3,10,17,1] } \\
& = & -\frac{1}{\sqrt{[3]}} [\underline{1,9,18,6,13,22,4,14,19,2,11,18,2,11},22,5,14,21,3,10,17,1] \\
& = & -\frac{\sqrt{[4]}}{\sqrt{[2][3]}} [\underline{1,9,18,6,13,22,4,14,19},2,11,18,6,11,22,5,14,21,3,10,17,1] \\
& = & -\frac{\sqrt{[3][4][9]}}{\sqrt{[2][5][7]}} [1,9,18,6,13,22,4,12,19,2,11,18,6,11,22,5,14,21,3,10,17,1]. \\
\end{eqnarray*}
We choose $u_{i\nu(i)} = u_{1\nu(1)} := u$. Then $v_{j\nu(i)} \nu(a_{ij}) = v_{9\nu(1)} \nu(a_{1,9})$ is given by
\begin{eqnarray*}
\lefteqn{ [\underline{9,18,6,13,22,4,14,19},2,11,19,2,11,22,5,14,21,3,10,17,1,9] } \\
& = & -\frac{[3]\sqrt{[9]}}{\sqrt{[5][7]}} [\underline{9,18,6,13,22,4,12,19,2},11,19,2,11,22,5,14,21,3,10,17,1,9] \\
& = & \frac{\sqrt{[3][9]}}{\sqrt{[5][7]}} [\underline{\xi_3,10,19,2},11,22,5,14,21,3,10,17,1,9] \\
& = & \frac{\sqrt{[3][9]}}{[5]\sqrt{[7]}} [\underline{\xi_3,9,17,2,11,22},5,14,21,3,10,17,1,9] \\
& = & \frac{[9]\sqrt{[3]^5[4]}}{\sqrt{[2][5][7]^3}} [\xi_3,9,18,6,13,22,5,14,21,3,10,17,1,9] \\
&& + \frac{[4][6]\sqrt{[2][9][10]}}{[5]\sqrt{[3]}} [\xi_2,12,19,2,11,22,5,14,21,3,10,17,1,9] \\
&& + \frac{\sqrt{[2]^3[3][4]^3}}{[5]\sqrt{[9]}} ([3]^2+1) [\xi_2,12,19,3,14,22,5,14,21,3,10,17,1,9].
\end{eqnarray*}
Now
\begin{eqnarray*}
\lefteqn{ \frac{[4][6]\sqrt{[2][9][10]}}{[5]\sqrt{[3]}} [\underline{\xi_2,12,19,2,11,22,5},14,21,3,10,17,1,9] } \\
&& + \frac{\sqrt{[2]^3[3][4]^3}}{[5]\sqrt{[9]}} ([3]^2+1) [\underline{\xi_2,12,19,3,14,22,5},14,21,3,10,17,1,9] \\
& = & -\frac{[2]^2\sqrt{[4]^3[10]}}{\sqrt{[5]^3[9]}} [\underline{\xi_2,12,19,2,10,19,5,14},21,3,10,17,1,9] \;\; = \;\; 0,
\end{eqnarray*}
so we have
\begin{eqnarray*}
v_{j\nu(i)} \nu(a_{ij}) & = & \frac{[9]\sqrt{[3]^5[4]}}{\sqrt{[2][5][7]^3}} [\underline{\xi_3,9,18,6,13,22,5,14,21,3,10,17,1,9}] \\
& = & -\frac{[2][4][9][12]\sqrt{[3]^3[10]}}{\sqrt{[5][7]^3}} [\xi_2,12,19,2,9,18,6,13,22,4,12,19,2,9] \;\; \neq \;\; 0.
\end{eqnarray*}
Now
\begin{eqnarray*}
\lefteqn{ [9,17,2,11,20,5,14,21,3,10,17,2,9,18,6,13,22,4,12,19,2,9] } \\
& = & \frac{\sqrt{[4][5]}}{\sqrt{[2]}} [9,18,6,13,20,5,14,21,3,10,17,2,9,18,6,13,22,4,12,19,2,9] \\
& = & -\frac{[3][5]\sqrt{[4]}}{\sqrt{[2][7]}} [9,18,6,13,22,5,14,21,3,10,17,2,9,18,6,13,22,4,12,19,2,9] \\
& = & \frac{[5][12]\sqrt{[2][7]}}{[3]\sqrt{[10]}} [\xi_2,12,19,2,9,18,6,13,22,4,12,19,2,9] \;\; \neq \;\; 0,
\end{eqnarray*}
so that
\begin{eqnarray*}
\lefteqn{ v_{j\nu(i)} \nu(a_{ij}) } \\
& = & -\frac{[4][9][10]\sqrt{[2][3]^5}}{[7]^2\sqrt{[5]^3}} [9,17,2,11,20,5,14,21,3,10,17,2,9,18,6,136,22,4,12,19,2,9] \\
& = & a_{jk} v',
\end{eqnarray*}
where $v'$ is symmetric. Then $C=1$.

\end{appendix}

\end{document}